\documentclass[12pt]{article}
\usepackage{amsmath}
\usepackage{amstext}
\usepackage{amssymb}
\usepackage{amsfonts}
\usepackage{amsthm}
\usepackage{amsopn}
\usepackage{amsbsy}
\usepackage{layout}

\def\newtheorems{\newtheorem{theorem}{Theorem}[section]
                 \newtheorem{cor}[theorem]{Corollary}
                 \newtheorem{prop}[theorem]{Proposition}
                 \newtheorem{lemma}[theorem]{Lemma}
                 \newtheorem{defn}[theorem]{Definition}
                 \newtheorem{definition}[theorem]{Definition}

                 \newtheorem{question}[theorem]{Question}
                 }

\newtheorems

\def\fs#1{\mbox{\it #1\kern 1.3pt}}
\def\rfs#1{\mbox{\rm #1\kern 1.3pt}}
\def\bfs#1{\mbox{\bf #1\kern 1.3pt}}
\def\fss#1{\mbox{\scriptsize\it #1\kern 1.3pt}}
\def\fst#1{\mbox{\tiny\it #1\kern 1.1pt}}
\def\sifs#1{\mbox{\scriptsize\it #1\kern 1.3pt}}
\def\srfs#1{\mbox{\kern0.7pt\scriptsize\rm #1\kern 1.3pt}}
\def\trfs#1{\mbox{\kern0.7pt\tiny\rm #1\kern 1.3pt}}
\def\tfs#1{\mbox{\kern0.7pt\tiny #1\kern 1.3pt}}
\def\sbfs#1{\mbox{\kern0.7pt\srbf #1\kern -0.6pt}}
\def\srbfs#1{\mbox{\kern0.7pt\srbf #1\kern -0.6pt}}
\def\spfs#1{\mbox{\kern0.7pt\scmu #1\kern 1.3pt}}
\def\sspfs#1{\mbox{\kern0.5pt\sscmu #1\kern 1.1pt}}
\def\ssbfs#1{\mbox{\kern0.7pt\ssbf #1\kern 1.3pt}}
\def\fsm#1{\mbox{\tiny\it #1\kern 1.0pt}}

\newcommand{\comp}{\hbox{$<\kern -3pt >$}}
\newcommand{\ncomp}
		{\;\hbox{\hbox{/}\kern -9.5pt \hbox{$<\kern -3pt >$}}}

\newcommand{\meet}
	       {\hbox{$\wedge \kern -5.75pt \raise 1.5pt \hbox{$.$}\,$}}
\newcommand{\Meet}
	     {\hbox{$\bigwedge \kern -8pt \raise 0.75pt \hbox{$.$}\:$}}
\newcommand{\ld}
	       {\hbox{$< \kern -6pt \raise 2pt \hbox{$.$}\,$}}
\newcommand{\sss}{\: \hbox{$
\underline{\hbox{$\subset$}}\kern -4pt\raise -2pt \hbox{$\tiny |$}  
$}\: }
\newcommand{\almostcontained}{\mathrel{\hbox{$
\raise 1.5pt \hbox{\scriptsize $\subset$}\kern -6.3pt\raise -3.5pt
\hbox{\scriptsize $\sim$}  
$}\:}}
\newcommand{\rraro}[2]{\hbox{$\kern 3pt\raise 2pt \hbox{$\raro$}
 \kern -14pt \raise
-3.5pt\hbox{\tiny{$#1\raro #2$}}$}}

\newcommand{\ct}{\centerline}

\newcommand{\frc}{\hbox{$\parallel \kern -5.7pt \hbox{$-$}$}}

\newcommand{\nfrc}{\not \kern -5pt \frc}
\newcommand{\rest}{\vbox{\hbox{$\:\kern -2pt\mathbin{\vert\kern-3.1pt\lower-1pt
   \hbox{$\mathsurround=0pt\mathchar"0012$}\kern-4pt}\:$}}}

\newcommand{\dsum}
{\mathbin{\raise1pt\hbox{$-\kern-13pt-\kern-8.65pt
\raise1pt\hbox{\scriptsize$\vert$}\kern-1.1pt
\raise1pt\hbox{\scriptsize$\vert$}\kern3.7pt$}}}
\newcommand{\Dsum}
{\mathbin{\raise0.4pt\hbox{$-\kern-8.5pt-\kern-11.23pt
\raise0pt\hbox{\small$\vert$}\kern-1.5pt
\raise0pt\hbox{\small$\vert$}\kern4.7pt$}}}

\newcommand{\drest}{\vbox{\hbox{$\:\kern -2pt\mathbin{\rest\kern-4.6pt\lower1.7pt
   \hbox{$\mathsurround=0pt\mathchar"0012$}\kern-4pt}\:$}}}
\newcommand{\vdrest}{\vbox{\hbox{$\:\kern -2pt\mathbin{\rest\kern-3.95pt\lower1.7pt
   \hbox{$\mathsurround=0pt\mathchar"0012$}\kern-4pt}\:$}}}
   \newcommand{\smdrest}{\vbox{\hbox{$\:\kern -2pt\mathbin{\smrest\kern-2.7pt\lower1.7pt
      \hbox{$\mathsurround=0pt\mathchar"0012$}\kern-4pt}\:$\kern 1pt}}}

\newcommand{\rests}{\vbox{\hbox{\scriptsize$\:\kern
-1.4pt\mathbin{\vert\kern-2.4pt\lower-1pt
\hbox{\scriptsize$\mathsurround=0pt\mathchar"0012$}\kern-
3.0pt}\:$}}}

\newcommand{\nrestriction}{\kern-2.5pt\upharpoonright\kern-2.5pt}
\newcommand{\mrestriction}{\kern-0.1pt\upharpoonright\kern-2.5pt}

\newcommand{\itGamma}{{\it \Gamma}}
\newcommand{\itDelta}{{\it \Delta}}

\newcommand{\itOmega}{{\it \Omega}}

\newcommand{\bfOmega}{{\bf \Omega}}

\newcommand{\gammaprime}{\gamma\kern1pt'}

\newcommand{\calA}{{\cal A}}
\newcommand{\calB}{{\cal B}}
\newcommand{\calC}{{\cal C}}
\newcommand{\calD}{{\cal D}}
\newcommand{\calE}{{\cal E}}
\newcommand{\calF}{{\cal F}}

\newcommand{\calH}{{\cal H}}
\newcommand{\calI}{{\cal I}}

\newcommand{\calM}{{\cal M}}

\newcommand{\calO}{{\cal O}}
\newcommand{\calP}{{\cal P}}
\newcommand{\calQ}{{\cal Q}}
\newcommand{\calR}{{\cal R}}
\newcommand{\calS}{{\cal S}}
\newcommand{\calT}{{\cal T}}
\newcommand{\calU}{{\cal U}}
\newcommand{\calsubU}{{\cal U\kern0.4pt}}
\newcommand{\calV}{{\cal V}}
\newcommand{\calsubV}{{\cal V\kern0.4pt}}
\newcommand{\calW}{{\cal W}}
\newcommand{\calX}{{\cal X}}

\newcommand{\hatomega}{\hat{\omega}}

\newcommand{\whatF}{\widehat{F}}

\newcommand{\whatH}{\widehat{H}}

\newcommand{\whatX}{\widehat{X}}
\newcommand{\swhatX}{\kern6pt\widehat{\rule{0pt}{8pt}}\kern-5.9pt X}
\newcommand{\whatY}{\widehat{Y}}

\newcommand{\whatcalM}{\widehat{\calM}}
\newcommand{\whatcalX}{\widehat{\calX}}

\newcommand{\vecmu}{\vec{\mu}}
\newcommand{\vecnu}{\vec{\nu}}

\newcommand{\swtildeX}{\kern6pt\widetilde{\rule{0pt}{8pt}}\kern-5.9pt X}

\newcommand{\raro}{\rightarrow}

\newcommand{\itmb}[1]{\item[[\kern 1.5pt #1\kern -4pt]]}
\newcommand{\itms}[1]{\item[[#1\kern -5pt]]}

\newcommand{\II}{{\bf I\kern -1pt I}}

\newcommand{\eqdf}{\hbox{\bf \,:=\,}}

\newcommand{\surj}{\vbox{\hbox{$\longrightarrow $
                  \kern -22pt \hbox{\lower 2.5pt  \hbox{\tiny onto}}
                  \kern -16pt \hbox{\raise 5pt  \hbox{\tiny 1-1}}
                  \kern 3pt}}}
\newcommand{\uarrow}[2]{\vbox{\hbox{$\longrightarrow $
                  \kern -16pt \hbox{\raise 5pt  \hbox{\tiny $#1$}}
                  \kern 10pt }}}

\newcommand{\spri}{\vbox{\hbox{\raise 2pt \hbox{\tiny $\|$}}}}
\newcommand{\tlr}{\vbox{\hbox{\raise 2pt \hbox{\tiny $\leftarrow$}}}}
\newcommand{\trr}{\vbox{\hbox{\raise 2pt \hbox{\tiny $\rightarrow$}}}}
\newcommand{\spr}{\mathrel{\vbox{\hbox{\tlr \kern -1.7pt \spri \kern -1.7pt \trr}}}}

\newcommand{\mspri}{\vbox{\hbox{\raise 2pt \hbox{$\scriptscriptstyle \|$}}}}

\newcommand{\mtlr}{\vbox{\hbox{\raise 2pt \hbox{$\scriptscriptstyle \leftarrow$}}}}

\newcommand{\mtrr}{\vbox{\hbox{\raise 2pt \hbox{$\scriptscriptstyle \rightarrow$}}}}

\newcommand{\mspr}{\mathrel{\vbox{\hbox{\mtlr \kern -1.7pt \mspri \kern -1.7pt \mtrr}}}}

\newcommand{\Nm}{\hbox{\kern -1.3pt \em I\kern-.2300em N\,}}
\newcommand{\Na}{\hbox{\kern -1.3pt \it I\kern-.2300em N\,}}
\newcommand{\Aa}{\hbox{\it A\kern -6.8pt\lower 1.0pt\hbox{-}\,}}
\newcommand{\Ba}{\hbox{\kern -1.3pt \it I\kern-.2300em B\,}}
\newcommand{\Da}{\hbox{\kern -1.3pt \it I\kern-.2300em D\,}}
\newcommand{\Ka}{\hbox{\kern -1.3pt \it I\kern-.2300em K\,}}
\newcommand{\La}{\hbox{\kern -1.3pt \it I\kern-.2300em L\,}}
\newcommand{\Ta}{\hbox{\kern 0.5pt \it T\kern-.5550em T\,}}
\newcommand{\jTa}{\hbox{\kern 0.5pt \it T\kern-.5550em T\,}}
\newcommand{\nTa}{\hbox{\kern 0.5pt \it T\kern-.6300em T\,}}
\newcommand{\Nmt}{\hbox{\kern -1.3pt I\kern-.2300em N\,}}
\newcommand{\N}{\hbox{I\kern-.1500em \hbox{\sf N}}}
\newcommand{\As}{\hbox{\scriptsize\it A\kern -5.3pt\lower 0.7pt\hbox{\it -}\,}}
\newcommand{\Ns}{\hbox{{\scriptsize\it I\kern-.1500em N}}}
\newcommand{\Nss}{\hbox{{\tiny\it I\kern-.1500em N}}}
\newcommand{\Ls}{\hbox{{\scriptsize\it I\kern-.1700em L}}}
\newcommand{\Rm}{\hbox{\kern -1.3pt \em I\kern-.1950em R\,}}
\newcommand{\Ra}{\hbox{\kern -1.3pt \it I\kern-.1950em R\,}}
\newcommand{\R}{\hbox{I\kern-.1500em \hbox{\sf R}}}
\newcommand{\sR}{\hbox{\tiny \hbox{I\kern-.1500em \hbox{\sf R}}}}
\newcommand{\msR}{\hbox{\tiny \hbox{\it I\kern-.2200em \hbox{\it R}}}}
\newcommand{\Rs}{\hbox{\scriptsize\it \hbox{I\kern-.2100em \hbox{R}}}}

\newcommand{\Q}
   {\hbox{${\rm Q} \kern -7.5pt \raise 2pt \hbox{\tiny$|$}\kern 7.5pt$}}
\newcommand{\C}
   {\hbox{${\rm C} \kern -6.5pt \raise 2pt \hbox{\tiny$|$}\kern 6.5pt$}}
\newcommand{\sC}
   {\hbox{\tiny \hbox{${\rm C} \kern -5.0pt \raise 1pt \hbox{\tiny$|$}\kern 7.5pt$}}}
\newcommand{\Z}{\hbox{\sf Z\kern-0.720em\hbox{ Z}}}
\newcommand{\sZ}{\hbox{\tiny\hbox{ \sf Z\kern-0.720em\hbox{ Z}}}}

\newcommand{\myqed}{\kern 5pt\vrule height7.5pt width6.9pt depth0.2pt
\rule{1.3pt}{0pt}}
\newcommand{\proofend}
{
\hbox{
\rule{0.8pt}{7pt}
\kern-3.9pt
\raise6.3pt
\hbox{\rule{4.1pt}{0.8pt}}
\kern-3.8pt
\rule{0.8pt}{7pt}
\kern-9.0pt
\rule{4.9pt}{0.8pt}
}
\kern-06pt}
\newcommand{\proofendeol}
{
\hbox{
\rule{0.8pt}{7pt}
\kern-3.9pt
\raise6.3pt
\hbox{\rule{4.1pt}{0.8pt}}
\kern-3.8pt
\rule{0.8pt}{7pt}
\kern-9.0pt
\rule{4.9pt}{0.8pt}
\kern-5pt
}
\kern-8pt
}

\newcommand{\bd}{\begin{description}}
\newcommand{\ed}{\end{description}}

\newcommand{\ben}{\begin{enumerate}}
\newcommand{\een}{\end{enumerate}}

\newcommand{\sngltn}[1]{\{ #1 \}}
\newcommand{\dbltn}[2]{\{ #1, #2 \}}

\newcommand{\setm}[2]{\{#1\:|\:#2\}}
\newcommand{\setmlr}[2]{\left\{#1\:|\:#2\right\}}

\newcommand{\fsetn}[2]{\{\,#1,\ldots ,#2\}}

\newcommand{\boldnorm}[1]
{
\kern1pt\lower3pt\hbox{\rule{1.0pt}{12pt}\kern1.4pt}
#1
\lower3pt\hbox{\kern1.6pt\rule{1.0pt}{12pt}\kern2.0pt}
                                            }

\newcommand{\bboldnorm}[1]
{
\kern1pt\lower3pt\hbox{\rule{1.2pt}{12pt}\kern1.4pt}
#1
\lower3pt\hbox{\kern1.6pt\rule{1.2pt}{12pt}\kern2.0pt}
                                            }

\newcommand{\sboldnorm}[1]
{
\kern1pt\lower1.6pt\hbox{\rule{0.9pt}{07pt}\kern1.1pt}
#1
\lower1.6pt\hbox{\kern1.1pt\rule{0.9pt}{07pt}\kern1.8pt}
                                            }

\newcommand{\sbboldnorm}[1]
{
\kern1pt\lower0.5pt\hbox{\rule{1.0pt}{06pt}\kern1.1pt}
#1
\lower0.5pt\hbox{\kern1.1pt\rule{1.0pt}{06pt}\kern1.8pt}
                                            }
\newcommand{\abs}[1]{| #1 |}

\newcommand{\pair}[2]{\langle#1 ,#2\rangle}

\newcommand{\trpl}[3]{\langle#1 ,#2 ,#3 \rangle}

\newcommand{\change}[1]{\lower 0.7pt \hbox{\mbox{\large $#1$}}}

\newcommand{\fnn}[3]{#1:#2 \raro #3}
\newcommand{\tendin}[3]
{\lim_{#1\kern 1pt \ni\kern 1.5pt #2\kern 1pt \rightarrow\kern 1.5pt #3}}
\newcommand{\tend}[2]
{\lim_{#1\kern 1pt \rightarrow\kern 1.5pt #2}}

\newcommand{\iso}[3]{\hbox{$ #1 : #2 \cong #3 $}}

\newcommand{\titles}[1]{\bigskip\ct {\Large\bf #1}\bigskip \bigskip}

\newcommand{\authorss}[4]{\bigskip\ct {#1} \ct{#2}
\ct{#3} \ct{#4} }

\newcommand{\inverse}{^{-1}}

\newcommand{\tsubseteq}{\lower2.1pt\hbox{\tiny$\subseteq$}}
\newcommand{\tgeq}{\raise2.2pt\hbox{\tiny$\geq$}}
\newcommand{\tleq}{\raise2.2pt\hbox{\tiny$\leq$}}
\newcommand{\tngeq}{\raise0.5pt\hbox{\tiny$\geq$}\kern1pt}
\newcommand{\tsucceq}{\raise0.5pt\hbox{\tiny$\succeq$}\kern1pt}
\newcommand{\tnleq}{\raise0.5pt\hbox{\tiny$\leq$}\kern1pt}
\newcommand{\tpreceq}{\raise0.5pt\hbox{\tiny$\preceq$}\kern1pt}
\newcommand{\temptyset}{\lower1.6pt\hbox{\tiny$\emptyset$}}
\newcommand{\tSml}{\lower1.6pt\hbox{\tiny\it Sml}}
\newcommand{\tSpprtd}{\lower1.6pt\hbox{\tiny\it Spprtd}}
\newcommand{\tprec}{\lower1.6pt\hbox{\tiny$\prec$}}
\newcommand{\tcong}{\lower1.6pt\hbox{\tiny$\cong$}}
\newcommand{\tspr}{\lower3.5pt\hbox{\tiny$\spr$}}
\newcommand{\indexentry}[2]{#1 \hfill #2\newline}
\newcommand{\newprec}{\hbox{$\prec$\kern-2.6pt\raise6.3pt\hbox{\rule{2.5pt}{0.3pt}}\kern-2.5pt\raise-0.5pt\hbox{\rule{2.5pt}{0.3pt}}\kern3pt}}

\newcommand{\onetoone}{\hbox{1 \kern-3.2pt -- \kern-3.0pt 1}}
\newcommand{\onetoonen}
{\hbox
{1 \kern-3.2pt \rule[2.9pt]{4pt}{0.6pt} \kern-2.8pt 1}}

\newcommand{\num}[1]{$(#1)$}

\newcommand{\overbbA}
               {\kern3.1pt
	       \overline{\kern-1.1pt\bbA\kern-0.6pt}\kern0.6pt}
\newcommand{\oversbbA}
               {\kern2.6pt
	       \overline{\kern-2.6ptsboldbbA\kern-0.3pt}\kern0.3pt}
\newcommand{\overB}
               {\kern3.1pt\overline{\kern-3.1ptB\kern-0.6pt}\kern0.6pt}
\newcommand{\overE}
               {\kern3.1pt\overline{\kern-3.1ptE\kern-0.6pt}\kern0.6pt}
\newcommand{\oversE}
               {\kern2.6pt\overline{\kern-2.6ptE\kern-0.3pt}\kern0.3pt}
\newcommand{\overF}
               {\kern3.1pt\overline{\kern-3.1ptF\kern-0.6pt}\kern0.6pt}
\newcommand{\oversF}
               {\kern2.6pt\overline{\kern-2.6ptF\kern-0.3pt}\kern0.3pt}
\newcommand{\overG}
               {\kern3.1pt\overline{\kern-3.1ptG\kern-0.6pt}\kern0.6pt}
\newcommand{\overN}
               {\kern3.1pt\overline{\kern-0.9ptN\kern-0.3pt}\kern0.6pt}
\newcommand{\overV}
               {\kern3.1pt\overline{\kern-0.9ptV\kern-0.3pt}\kern0.6pt}
\newcommand{\overW}
               {\kern3.1pt\overline{\kern-0.9ptW\kern-0.3pt}\kern0.6pt}
\newcommand{\overX}
               {\kern3.1pt\overline{\kern-3.1ptX\kern-0.6pt}\kern0.6pt}
\newcommand{\oversG}
               {\kern2.6pt\overline{\kern-2.6ptG\kern-0.3pt}\kern0.3pt}
\newcommand{\oversX}
               {\kern2.6pt\overline{\kern-2.6ptX\kern-0.3pt}\kern0.3pt}

\newcommand{\overcalO}
               {\kern3.1pt\overline{\kern-3.1pt{\cal O}\kern-0.6pt}\kern0.6pt}
\newcommand{\arrowcalP}
               {\kern3.1pt\overline{\rule{0.80pt}{0pt}\rule{0pt}{8.8pt}
	       \kern-3.1pt\cal P\kern-0.6pt}
	       \kern-7pt\vec{\phantom{P}}\kern0.6pt}
\newcommand{\arrowcalQ}
               {\kern3.1pt\overline{\rule{0.80pt}{0pt}\rule{0pt}{8.8pt}
	       \kern-3.1pt\cal Q\kern-0.6pt}
	       \kern-7pt\vec{\phantom{Q}}\kern0.6pt}
\newcommand{\arrowcalR}
               {\kern3.1pt\overline{\rule{0.80pt}{0pt}\rule{0pt}{8.8pt}
	       \kern-3.1pt\cal R\kern-0.6pt}
	       \kern-7pt\vec{\phantom{R}}\kern0.6pt}
\newcommand{\overarrowP}
               {\kern3.1pt\overline{\rule{0.30pt}{0pt}\rule{0pt}{8.8pt}
	       \kern-3.1ptP\kern-0.6pt}
	       \kern-7pt\vec{\phantom{P}}\kern0.6pt}
\newcommand{\sarrowcalP}{\hbox{$\vec{\hbox{\scriptsize $\calP$}}$}}

\newcommand{\sarrowcalPX}{\sarrowcalP^{\lower3pt\hbox{\tiny$X$}}}

\newcommand{\sstar}{\hbox{\scriptsize \kern1pt$\star$\kern0pt}}
\newcommand{\rsstar}{\raise1pt\hbox{\scriptsize \kern1pt$\star$\kern1pt}}
\newcommand{\raisedsstar}{\kern0.5pt\raise1.0pt\hbox{\scriptsize$*$}}
\newcommand{\sperp}{\raise2pt\hbox{\kern0pt\tiny $\perp$}}
\newcommand{\ssperp}{\raise1.5pt\hbox{\kern0pt\tinier \symbol{63}}}

\newcommand{\nperp}
{
\hbox{\kern1.0pt
\hbox{\raise1.0pt \hbox{\rule{4.6pt}{0.3pt}}}\kern-5.1pt$\sperp$}
}
\newcommand{\snperp}
{
\hbox{\kern1.0pt
\hbox{\raise0.4pt \hbox{\rule{3.3pt}{0.2pt}}}\kern-3.65pt$\ssperp$}
}

\newcommand{\sor}{\raise2pt\hbox{\tiny\rm Or}}

\newcommand{\srs}[1]{\raise1.3pt\hbox{\tiny\rm #1}}
\newcommand{\ssrfs}[1]{\raise1.3pt\hbox{\tiny\rm #1}}
\newcommand{\ssrs}[1]{\raise1.3pt\hbox{\nmini #1}}
\newcommand{\scirc}
{\raise1pt\hbox{\scriptsize\kern1.5pt$\circ$\kern1.5pt}}
\newcommand{\sscirc}
{\hbox{\tiny\kern1.5pt$\circ$\kern0.3pt}}

\newcommand{\bcirc}
{\mathop{\lower1pt\hbox{\large\kern1.5pt$\circ$}}}
\newcommand{\bbcirc}
{\mathop{\lower1.8pt\hbox{\Large\kern1.5pt$\circ$}}}

\newcommand{\raisedstar}
{\raise-3.1pt\hbox{$\rule{0.5pt}{0pt}^*$}}
\newcommand{\rstar}
{\raise-3.1pt\hbox{$\rule{0.5pt}{0pt}^*$}}
\newcommand{\raisedtinystar}
{\kern0.7pt\raise1.5pt\hbox{\tiny $*$}}
\newcommand{\rtstar}
{\kern0.7pt\raise1.0pt\hbox{\tiny $*$}\kern-2.0pt}
\newcommand{\bfie}{\hbox{\blcmssifont \symbol{101}\kern1pt}}
\newcommand{\bfif}{\hbox{\blcmssifont \symbol{102}\kern1pt}}

\newcommand{\bfih}{\hbox{\blcmssifont \symbol{104}\kern1pt}}

\newcommand{\mcdot}{\kern-1.4pt\cdot\kern-1.4pt}
\newcommand{\lcdot}{\hbox{$\kern0.1pt\cdot\kern1.0pt$}}
\newcommand{\ncdot}{\hbox{$\kern-0.2pt\cdot\kern0.6pt$}}
\newcommand{\fr}{\kern1.0ptr}
\newcommand{\fzero}{\kern1.3pt0}
\newcommand{\fone}{\kern1.0pt1}
\newcommand{\ftwo}{\kern1.3pt2}
\newcommand{\fprime}{\kern1.5pt'}
\newcommand{\nprime}{\kern0.8pt'}
\newcommand{\sprime}{\raise1.3pt\hbox{\scriptsize\kern1.4pt$'$}}
\newcommand{\fprimesub}{\kern1.5pt'\kern-3.8pt}
\newcommand{\fprimew}{\kern0.7pt'}
\newcommand{\fprimei}{\kern1.5pt'\kern-3.0pt}
\newcommand{\fnprime}{\kern1.5pt'\kern-3.0pt}

\newcommand{\clubsign}{{\raise1.3pt\hbox{\tiny $\clubsuit$}}}
\newcommand{\sclub}{{\raise1.3pt\hbox{\tiny $\clubsuit$}}}
\newcommand{\srsp}[1]{\hbox{\kern1pt\tiny $(#1)$}}
\newcommand{\sminus}{\rule{5.6pt}{0.3pt}}
\newcommand{\ssim}{\hbox{\scriptsize $\sim$}}
\newcommand{\neweq}{
\mathrel{
\hbox{
\lower-0.0pt\hbox{\sminus}
\kern-9.95pt\ssim\kern-6.2pt\raise3.8pt\hbox{\sminus}
}}}

\newcommand{\dcup}{\mathbin{\hbox{$\cup$
\kern-10.43pt\raise1.4pt \hbox{\tiny $\cup$}}}}

\newcommand{\sprt}[1]
{\lower1.3pt\hbox{\kern1.5pt\rule{0.5pt}{10pt}}
\kern-0.5pt\underline{\kern1pt#1\kern1pt}
\lower1.3pt\hbox{\kern0.2pt\rule{0.5pt}{10pt}}\kern1pt}

\newcommand{\sprtd}[2]
{#1\lower1.3pt\hbox{\kern1.5pt\rule{0.5pt}{10pt}}
\kern-0.5pt\underline{\kern1pt#2\kern1pt}
\lower1.3pt\hbox{\kern0.2pt\rule{0.5pt}{10pt}}\kern1pt}

\newcommand{\sprtl}[1]
{\lower4.4pt\hbox{\kern1.5pt\rule{0.5pt}{14pt}}
\kern-0.5pt\underline{\kern1pt#1\kern1pt}
\lower4.4pt\hbox{\kern-0.5pt\rule{0.5pt}{14pt}}\kern1pt}

\newcommand{\sprtll}[1]
{\lower5.8pt\hbox{\kern1.5pt\rule{0.5pt}{15.4pt}}
\kern-0.5pt\underline{\kern1pt#1\kern1pt}
\lower5.8pt\hbox{\kern-0.5pt\rule{0.5pt}{15.4pt}}\kern1pt}

\newcommand{\sprtdl}[2]
{#1\lower4.0pt\hbox{\kern1.5pt\rule{0.5pt}{14pt}}
\kern-0.5pt\underline{\kern1pt#2\kern1pt}
\lower4.1pt\hbox{\kern0.0pt\rule{0.5pt}{14pt}}\kern1pt}

\newcommand{\sprtm}[1]
{\lower3.3pt\hbox{\kern1.5pt\rule{0.5pt}{14pt}}
\kern-0.5pt\underline{\kern1pt#1\kern1pt}
\lower3.4pt\hbox{\kern0.0pt\rule{0.5pt}{14pt}}\kern1pt}

\newcommand{\sprtdm}[2]
{#1\lower3.3pt\hbox{\kern1.5pt\rule{0.5pt}{14pt}}
\kern-0.5pt\underline{\kern1pt#2\kern1pt}
\lower3.4pt\hbox{\kern0.0pt\rule{0.5pt}{14pt}}\kern1pt}

\DeclareSymbolFont{AMSb}{U}{msb}{m}{n}
\DeclareMathSymbol{\bbA}{\mathbin}{AMSb}{"41}
\DeclareMathSymbol{\bbB}{\mathbin}{AMSb}{"42}
\DeclareMathSymbol{\bbC}{\mathbin}{AMSb}{"43}
\DeclareMathSymbol{\bbD}{\mathbin}{AMSb}{"44}
\DeclareMathSymbol{\bbE}{\mathbin}{AMSb}{"45}
\DeclareMathSymbol{\bbF}{\mathbin}{AMSb}{"46}
\DeclareMathSymbol{\bbG}{\mathbin}{AMSb}{"47}
\DeclareMathSymbol{\bbH}{\mathbin}{AMSb}{"48}
\DeclareMathSymbol{\bbI}{\mathbin}{AMSb}{"49}
\DeclareMathSymbol{\bbJ}{\mathbin}{AMSb}{"4A}
\DeclareMathSymbol{\bbK}{\mathbin}{AMSb}{"4B}
\DeclareMathSymbol{\bbL}{\mathbin}{AMSb}{"4C}
\DeclareMathSymbol{\bbM}{\mathbin}{AMSb}{"4D}
\DeclareMathSymbol{\bbN}{\mathbin}{AMSb}{"4E}
\DeclareMathSymbol{\bbO}{\mathbin}{AMSb}{"4F}
\DeclareMathSymbol{\bbP}{\mathbin}{AMSb}{"50}
\DeclareMathSymbol{\bbQ}{\mathbin}{AMSb}{"51}
\DeclareMathSymbol{\bbR}{\mathbin}{AMSb}{"52}
\DeclareMathSymbol{\bbS}{\mathbin}{AMSb}{"53}
\DeclareMathSymbol{\bbT}{\mathbin}{AMSb}{"54}
\DeclareMathSymbol{\bbU}{\mathbin}{AMSb}{"55}
\DeclareMathSymbol{\bbV}{\mathbin}{AMSb}{"56}
\DeclareMathSymbol{\bbW}{\mathbin}{AMSb}{"57}
\DeclareMathSymbol{\bbX}{\mathbin}{AMSb}{"58}
\DeclareMathSymbol{\bbY}{\mathbin}{AMSb}{"59}
\DeclareMathSymbol{\bbZ}{\mathbin}{AMSb}{"5A}
\newcommand{\boldbbA}{\mathbb{A}\kern-9pt\mathbb{A}}
\newcommand{\boldbbL}{\mathbb{L}\kern-8.3pt\mathbb{L}}

\newcommand{\sboldbbA}
                {\hbox{\kern0.7pt\scriptsize
		$\mathbb{A}\kern-6pt\mathbb{A}$}}
\newcommand{\sboldbbL}
                {\hbox{\kern0.7pt\scriptsize
		$\mathbb{L}\kern-5.68pt\mathbb{L}$}}
\newcommand{\sboldbbN}
                {\hbox{\kern0.7pt\scriptsize
		$\mathbb{N}\kern-5.60pt\mathbb{N}$}}
\newcommand{\sboldbbT}
                {\hbox{\kern0.7pt\scriptsize
		$\mathbb{T}\kern-5.60pt\mathbb{T}$}}

\newcommand{\wpm}{\kern0.7pt\pm}

\newcommand{\rmsup}[1]
{^{\raise1pt\hbox{\scriptsize\rm \kern0.5pt#1}}}

\newcommand{\rmsupsup}[1]
{^{\raise1pt\hbox{\tiny\rm \kern0.5pt#1}}}

\def\bfP{{\bf P}}
\def\bfGamma{{\bf\Gamma}}
\newcommand{\Erdos}{Erd\"os}

\newcommand{\beths}{\beth\kern1.7pt}
\newcommand{\sbeths}{\kern1pt\beth\kern1.7pt}

\newcommand{\symdiff}
{\hbox{$\mathbin{\hbox{\tiny \kern2pt$\triangle$\kern2pt}}$}}

\begin{document}
\baselineskip 18pt
\titles{A classification of CO spaces which are}
\vspace{-1cm}
\titles{continuous images of compact ordered spaces
\kern0.5pt\normalsize
\footnote{
AMS Subject Classification 2000:  Primary 06E05.
Secondary 54G12, 06A05
\newline
Keywords: scattered spaces, ordered spaces, superatomic Boolean algebras.
}}

    \bigskip

    \authorss
{\bf Robert Bonnet}
{D\'epartement de Math\'ematique}{Universit\'e de Savoie}{Chamb\'ery,
France}

    \bigskip

    \ct{ and}

    \authorss
   {\bf Matatyahu Rubin}
{Department of Mathematics}{Ben Gurion University of the
    Negev}{Beer Sheva, Israel\kern1pt\small
}

\kern2cm

\begin{abstract}
\noindent
A topological space $X$ is called a CO space, if every closed subset of $X$
is homeomorphic to some clopen subset of $X$.
Every ordinal with its order topology is a CO space.
This work gives a complete classification of CO spaces which are continuous images
of compact ordered spaces.
\end{abstract}

\newpage
\section{Introduction}\label{s1}
A topological space $X$ is called a CO space if every closed subset
of $X$ is homeomorphic to some clopen subset of $X$.
The simplest example of a compact Hausdorff CO space is
a successor ordinal with its order topology.

In this work we characterize the CO spaces which are continuous images
of compact interval spaces. There are such spaces
which are not ordinals, yet this class is not much bigger
than the class of successor ordinals.

So far there has been only one result concerning
compact Hausdorff CO spaces which are not continuous images
of compact interval spaces. It is due to Bonnet and Shelah \cite{BS}.
Assuming $\diamondsuit_{\aleph_1}$ they construct a thin tall CO space.
The significance of this result is that it indicates that there is no
explicit description of general compact Hausdorff CO spaces.

To state the main theorem of this work, we need the following terminology.
A space $\pair{X}{\tau^X}$ is an {\it interval space},
if there is a linear ordering $<$ of $X$ such that $\tau^X$
is the order topology of this linear ordering.
That is, a subbase for this topology is the family of sets
$$
\setm{\setm{x \in X}{x < a}}{a \in X}
\cup
\setm{\setm{x \in X}{x > a}}{a \in X}.
$$
An interval space $X$ is called an {\it ordinal space}
if there is a well ordering of $X$ such that $\tau^X$ is the order topology of
this well ordering.
For infinite cardinals $\lambda,\mu$,
let $\mu^*$ denote the reverse ordering of $\mu$
and $X_{\lambda,\mu}$ denote
the interval space of $\lambda + 1 + \mu^*$.
\index{N@$X_{\lambda,\mu}$. The interval space of $\lambda + 1 + \mu^*$}
\index{N@$\alpha(X_{\lambda,\mu}) = \max(\lambda,\mu) \cdot \omega$.}
\kern-5.5pt
Define $\alpha(X_{\lambda,\mu})$ to be the following ordinal:
$\alpha(X_{\lambda,\mu}) \eqdf \max(\lambda,\mu) \cdot \omega$.
For an infinite cardinal $\kappa$
let $X_{\kappa}$ denote the one point compactification of a discrete
space of cardinality $\kappa$
and set $\alpha(X_{\aleph_1}) \eqdf\omega^2$.
The notation $X \cong Y$ stands for the fact that
$X$ and $Y$ are homeomorphic,
and $\iso{f}{X}{Y}$ means that
$f$ is a homeomorphism between $X$ and $Y$.
The final result of this work is the following theorem.

\begin{theorem}\label{t1.1}
\num{a} Let $X$ be a Hausdorff space
which is a continuous image of a compact interval space,
and assume that $X$ is a CO space.
Then there is a partition $\fsetn{Z,Y_0}{Y_{k - 1}}$ of $X$
into open sets
such that
\begin{itemize}
\item[\num{1}]
For every $i < k$
either
$Y_i \cong X_{\aleph_1}$,
or 
$Y_i \cong X_{\lambda,\mu}$, where $\lambda,\mu$ are some infinite
regular cardinals and  $\mu > \aleph_0$.
\item[\num{2}]
$Z$ is an ordinal space homeomorphic to some successor ordinal $\beta$.
\item[\num{3}]
$\beta > \alpha(Y_i)$ for every $i \in I$.
\end{itemize}
Note that if $X_{\lambda,\mu} \cong X_{\lambda',\mu'}$,
then $\dbltn{\lambda}{\mu} = \dbltn{\lambda'}{\mu'}$.
So $\alpha(Y_i)$ is well-defined.

\num{b} If a space $X$ has the above form, then $X$ is a CO space,
and $X$ is a continuous image of a compact interval space.
\hfill\myqed
\end{theorem}

Part (b) of the above theorem is merely an observation.
It is Part (a) which is the real subject of this work.

\index{N@$K_{\srfs{CII}}$. The class of all Hausdorff spaces which are
       a continuous\\\indent
       image of a compact interval space}
A compact Hausdorff space $X$ is {\it scattered}
if every nonempty subset of $X$
has an isolated point in its relative topology.
Let $K_{\srfs{CII}}$ denote the class of all Hausdorff spaces which are
the continuous image of a compact interval space.
Section \ref{s2} deals with the following intermediate step in the proof
of Theorem~\ref{t1.1}.
\index{D@scattered space}

\begin{theorem}\label{scat-t1.1}
\label{t1.2}
If $X \in K_{\srfs{CII}}$ and $X$ is a CO space, then $X$ is scattered.
\end{theorem}

Three main questions arise.
\begin{question}\label{metr-bldr-q5.4}
\begin{rm}
(a) Is there a non-scattered compact Hausdorff CO space?
It is even not known whether it is consistent with ZFC that such a space exists.

(b) The construction of \cite{BS} works only for $\aleph_1$.
So we ask whether there is a compact Hausdorff CO space of cardinality $> \aleph_1$
which is not a finite direct sum of a member of the class defined in
Theorem~\ref{t1.1} and a CO space with cardinality $\aleph_1$?
It is even not known whether this statement is consistent.

(c) Does it follow from ZFC that there is
a compact Hausdorff CO space
which does not belong the the class defined in Theorem~\ref{t1.1}?
\end{rm}
\end{question}

Let $K_{\srfs{IVL}}$
be the class of $0$-dimensional compact interval spaces.
The classification those CO spaces which belong to $K_{\srfs{IVL}}$
was dealt with in \cite{BBR}. The classification theorem proved in
\cite{BBR} is of course a special case of Theorem~\ref{t1.1}.

After the authors had proved Theorem~\ref{t1.2},
Shelah proved a theorem which turned out to be almost equivalent to
\ref{t1.2}.
The statement of this theorem appears in \cite{S} p.355.
However, a proof of that theorem has never appeared.
That Shelah's statement is equivalent to a statement
about continuous images of interval spaces follows from~\cite{H}.

\kern1mm

{\bf The main steps in the proof of Theorem \ref{t1.1}.}

In Section \ref{s2} we prove that a CO space
which is a continuous image of a compact interval space
must be scattered (Theorem~\ref{t1.2}).
The rest of the sections deals with scattered spaces which are a continuous image
of a compact interval space.

Section \ref{s3} deals with the question: when a scattered
continuous image of a compact interval space is itself an interval space.
The characterization uses ``obstructions''.
We define a class $\calO$ of topological spaces,
and prove that every space which is
a scattered continuous image of a compact interval space,
and which does not embed any member of $\calO$ must be an interval space.
This statement appears in Theorem~\ref{t2.1}.

Let $X$ be a scattered continuous image of a compact interval space,
and suppose that $X$ is a CO space.
We shall show that $X$ is the sum of finitely many copies of $X_{\aleph_1}$
and a space $Y$ which omits all members of $\calO$. (See above).
Then we use the characterization of CO compact interval spaces from \cite{BBR}
to get a description of $Y$.

Section \ref{s4} contains the main technical lemma
needed in the proof that the obstructions are omitted (Theorem~\ref{t3.2}).
It says that if $X$ is a scattered CO space then there are no
subsets $M,L,K \subseteq X$ such that $M \prec L \prec^{\srfs{w}} K$.
(See Definition~\ref{d3.1}).
In Theorem~\ref{t3.2}, the CO space $X$ is assumed to have a very strong
Hausdorff property. Because of this assumption we are able to deal only with
continuous images of compact interval spaces and not with general compact spaces.

In section \ref{s5} we show that the obstructios are omitted
and in \ref{s6} we obtain the desired characterization.

As a matter of fact,
using Theorem~\ref{t3.2} there is a short clean proof of the characterization
of CO scattered compact interval spaces.
This has already been done in \cite{BBR}, but in a less elegant way.
So in Section~\ref{s7} we prove this characterization.
By doing so, this work becomes self-contained and easier to read.

\section{Scatteredness of CO spaces which are
continuous images of compact interval spaces}\label{s2}

\kern 2mm

In this section we prove Theorem~\ref{t1.2}.
The proof is by way of contradiction, but it takes
till Theorem~\ref{t2.28} to reach this contradiction.
In two of the intermediate lemmas -- Proposition~\ref{p2.6}
and Corollary~\ref{c2.9}, a space $X$ is given,
and it is assumed that $X$ is a non-scattered CO space.
Since it is not known whether non-scattered CO spaces exist,
these lemmas have little or no use once Theorem~\ref{t1.2} is proved.
In addition,
Propositions \ref{tightly-haus-p1.2}, \ref{tightly-haus-p1.3}
and \ref{tightly-haus-l1.4} assume the existence of a CO space
which has some extra propreties.
These assumptions too are likely to be contradictory.
See especially~\ref{tightly-haus-l1.4}.

We do not prove directly that every CO space
which belongs to $K_{\srfs{CII}}$ is scattered.
This turns out to be too tedious.
Rather, we find certain topological properties of members
of $K_{\srfs{CII}}$ which serve as interpolants.
For example, in Part 1.\ below,
we prove that every member $X$ of $K_{\srfs{CII}}$
is tightly Hausdorff,
and later we use {\it this} property of $X$ rather than assuming that
$X \in K_{\srfs{CII}}$.
There are four other properties of members of $K_{\srfs{CII}}$
which are used as interpolants,
and we prove them just before they are used.
The class of spaces with these five properties is denoted by $K$.
In Theorem~\ref{t2.28} we prove that every
member of $K$ which is a CO space is scattered.
Hence the same is true for members of $K_{\srfs{CII}}$.

\kern2mm

\noindent
{\bf 1. Some Hausdorff-type properties
of members of ${\bf K_{\srfs{CII}}}$.}
\smallskip

\noindent
We start by defining the notion of a tightly Hausdorff space.
We shall show that members of $K_{\srfs{CII}}$
are tightly Hausdorff. This property and some of its weaker variants
will be used extensively.

\begin{defn}\label{d2.2}
\begin{rm}
Let $X$ be a topological space.

(a) We denote by $\tau^X$ the topology of $X$.
If $A \subseteq X$, then $\rfs{cl}^X(A)$, $\rfs{int}^X(A)$ and
$\rfs{acc}^X(A)$ denote respectively the closure, interior
and the set of accumulaton points of $A$ in $X$.
If $x \in X$,
then the set of open neighborhoods of $x$ in $X$ is denoted by
$\rfs{Nbr}^X(x)$.
Similarly,
$\rfs{Nbr}^X_{\srfs{cl}}(x)$ and $\rfs{Nbr}^X_{\srfs{clp}}(x)$
denote respectively the
the set of closed neighborhoods of $x$ in $X$
and the set of clopen neighborhoods of $x$ in $X$.
Superscript~$^X$ is omitted when the indended space $X$
can be understood from the context.
\index{N@$\rfs{cl}^X(A)$. Closure of $A$ in $X$}
\index{N@$\rfs{int}^X(A)$. Interior of $A$ in $X$}
\index{N@$\rfs{acc}^X(A)$. Set of accumulation points of $A$ in $X$}
\index{N@$\rfs{Nbr}^X(x)$. Set of open neighborhoods of $x$ in $X$}
\index{N@$\rfs{Nbr}^X_{\srfs{cl}}(x)$.
       Set of closed neighborhoods of $x$ in $X$}
\index{N@$\rfs{Nbr}^X_{\srfs{clp}}(x)$.
       Set of clopen neighborhoods of $x$ in $X$}

(b) A family $\calA$ of subsets of $X$ is called
a {\it pairwise disjoint family},
if $A \cap B = \emptyset$ for any distinct $A,B \in \calA$.
Let $\calA$ a be pairwise disjoint family of subsets of $X$
and $x \in X$. We say that $x$ is an
{\it accumulation point of $\calA$},
if every neighborhood of $x$ intersects infinitely many members
of $\calA$.
The set of accumulation points of $\calA$ is denoted by
$\rfs{acc}(\calA)$.
Suppose that $\calA$ is a pairwise disjoint family of subsets of $X$,
such that for every $B,C \subseteq \bigcup \calA$,
if
$$
\setm{A \in \calA}{B \cap A \neq \emptyset} =
\setm{A \in \calA}{C \cap A \neq \emptyset},
$$
then
$$
\rfs{acc}(\setm{B \cap A}{A \in \calA}) =
\rfs{acc}(\setm{C \cap A}{A \in \calA}).
$$
Then $\calA$ is called a {\it tight family}.
\index{D@pairwise disjoint family of subsets of $X$}
\index{D@accumulation point of $\calA$}
\index{N@$\rfs{acc}(\calA)$. The set of accumulation points of a
       family of sets $\calA$}
\index{D@tight family of subsets of $X$}

(c) A subset $A \subseteq X$ is {\it relatively discrete} if
$A$ together with its relative topology is a discrete space.
So $A$ is relatively discrete iff
$A \cap \rfs{acc}^X(A) = \emptyset$. 
\index{D@relatively discrete.
       $A$ is relatively discrete if $A \cap \rfs{acc}(A) = \emptyset$}

(d) Let $A \subseteq X$.
A family $\calU = \setm{U_x}{x \in A}$ is called a
{\it Hausdorff system} for $A$, if $\calU$ is a pairwise disjoint
family and for every $x \in A$, $U_x \in \rfs{Nbr}(x)$.
\index{D@Hausdorff system@@Hausdorff system}

(e) We say that $\calU$ is a {\it strong Hausdorff system} for $A$,
if $\calU$ is a Hausdorff system for $A$
and $\rfs{acc}(\calU) = \rfs{acc}(A)$.
\index{D@strong Hausdorff system}

(f) Let $X$ be a Hausdorff space.
If every relatively discrete subset of $X$ has a Hausdorff system,
then $X$ is called a {\it collectionwise Hausdorff space}.
If every relatively discrete subset of $X$
has a strong Hausdorff system,
then $X$ is said to be a {\it strongly Hausdorff space},
and if every relatively discrete subset of $X$
has a tight Hausdorff system,
then we call $X$ a {\it tightly Hausdorff space}.
\rule{0pt}{0pt}\hfill\myqed
\index{D@collectionwise Hausdorff space.
       $X$ is collectionwise Hausdorff if every\\\indent
       relatively discrete subset of $X$ has a Hausdorff system}
\index{D@strongly Hausdorff space.
       $X$ is strongly Hausdorff if every relatively\\\indent
       discrete subset of $X$ has a strong Hausdorff system}
\index{D@tightly Hausdorff space.
       $X$ is tightly Hausdorff if every relatively\\\indent
       discrete subset of $X$ has a tight Hausdorff system}
\end{rm}
\end{defn}

Note that
\newline\centerline{
``tightly Hausdorff'' $\Rightarrow$
``strongly Hausdorff'' $\Rightarrow$
``collectionwise Hausdorff''.
}

\begin{lemma}\label{l2.3}
If $X \in K_{\srfs{CII}}$, then $X$ is tightly Hausdorff.
\end{lemma}

\noindent{\bf Proof }
Let $N$ be a subset of a chain $\pair{L}{<}$
and $I \subseteq N$ be a convex subset of $L$.
We say that $I$ is a convex component of $N$ in $L$ if there is
no convex set $I'$ such that $I' \subseteq N$ and $I'$
properly contains $I$.
Denote the family of convex components of $N$ in $L$ by $\calI(N)$.
Clearly, $\calI(N)$ is a partition of $N$,
and if $N$ is open in the order topology of $L$,
then every member of $\calI(N)$ is open.
\index{N@$\calI(N)$. The family of convex components of $N$
       in a linear ordering $L$}

Let $\pair{L}{<}$ be a compact chain
and $\fnn{f}{L}{X}$ be a continuous surjective function.
Denote the order topology of $\pair{L}{<}$ by $\tau^L$
and the topology of $X$ by $\tau^X$.
Suppose that $A \subseteq X$ is relatively discrete.
For every $x \in A$ we define $L_x \in \tau^L$
and $U_x \in \tau^X$.
Let $\setm{x_i}{i < \alpha}$ be an enumeration of $A$.
We define $L_{x_i}$ and $U_{x_i}$ by induction on $i$.
Suppose that $L_{x_j}$ and $U_{x_j}$
have been defined for every $j < i$,
set $A_0 = \setm{x_j}{j < i}$,
and assume the following induction hypotheses.
\begin{itemize}
\addtolength{\parskip}{-11pt}
\addtolength{\itemsep}{06pt}
\item[(1)]
For every $x \in A_0$ and $I \in \calI(L_x)$ there is $s_I \in I$
such that $f(s_I) = x$.
\item[(2)]
$L_x \cap L_y = \emptyset$ for every distinct $x,y \in A_0$.
\item[(3)]
$f\inverse(U_x) \subseteq L_x$ for every $x \in A_0$.
\item[(4)]
$f(\rfs{cl}(L_x)) \cap A = \sngltn{x}$ for every $x \in A_0$.
\end{itemize}

{\bf Claim 1 }
(i) If $s \in \rfs{acc}(\setm{L_x}{x \in A_0})$,
then $f(s) \in \rfs{acc}(A_0)$.
(ii)
If $s \in \rfs{cl}(\bigcup \setm{L_x}{x \in A_0})$,
then either $f(s) \in \rfs{acc}(A_0)$
or for some $x \in A_0$, $s \in \rfs{cl}(L_x)$.
\newline
{\bf Proof }
Statement (ii) follows trivially from (i).
Let $s$ be as in the (i) and $J$ be an open interval containing~$s$.
Then for every finite set $\sigma \subseteq A_0$
there are distinct $x,y,z \in A_0 - \sigma$ such that
$L_x,L_y,L_z$ intersect $J$.
Then there is $I \in \calI(L_x) \cup \calI(L_y) \cup \calI(L_z)$
such that $I \subseteq J$.
Assume that $I \in \calI(L_x)$.
Then $f(s_I) = x$.
This implies that for every neighborhood $N$ of $s$,
$f(N) \cap A_0$ is infinite.
So if $U \in \rfs{Nbr}(f(s))$, then $f\inverse(U) \in \rfs{Nbr}(s)$,
so $f(f\inverse(U))$ contains an infinite subset of $A_0$.
Now, $f(f\inverse(U)) = U$.
Hence $U$ contains an infinite subset of $A_0$.
So $f(s) \in \rfs{acc}(A_0)$.
Claim~1 is proved.

Denote $x_i$ by $y$,
and set $K = \rfs{cl}(\bigcup \setm{L_x}{x \in A_0})$.
Then $y \not\in f(K)$.
This relies on the following three facts.
\begin{itemize}
\addtolength{\parskip}{-11pt}
\addtolength{\itemsep}{06pt}
\item
$A \cap \bigcup_{x \in A_0} f(\rfs{cl}(L_x)) = A_0$.
\item
If $s \in \rfs{cl}(\bigcup \setm{L_x}{x \in A_0}) -
\bigcup_{x \in A_0} f(\rfs{cl}(L_x))$,
then $f(s) = \rfs{acc}(A_0)$.
\item
$A$ is relatively discrete.
\end{itemize}
Hence $V_y \eqdf X - f(K) \in \rfs{Nbr}(y)$.
Choose $W_y \in \rfs{Nbr}(y)$ such that 
$\rfs{cl}(W_y) \cap A = \sngltn{y}$
and define $M_y = f\inverse(V_y \cap W_y)$
and $L_y = \bigcup \setm{I \in \calI(M_y)}{y \in f(I)}$.
Clearly,
$\calI(L_y) = \setm{I \in \calI(M_y)}{y \in f(I)}$
and $f\inverse(y) \subseteq L_y$.

{\bf Claim 2 }
There is $U_y \in \rfs{Nbr}(y)$ such that
$f\inverse(U_y) \subseteq L_y$.
\newline
{\bf Proof }
Suppose that Claim 2 is false.
Then for every $F \in \rfs{Nbr}_{\srfs{cl}}(y)$,\break
$f\inverse(F) \cap (L - L_y) \neq \emptyset$.
So
$H \eqdf \bigcap \setm{f\inverse(F) \cap
(L - L_y)}{F \in \rfs{Nbr}_{\srfs{cl}}(y)} \neq
\emptyset$.
Let $a \in H$. Then for every $F \in \rfs{Nbr}_{\srfs{cl}}(y)$,
$f(a) \in F$. So $f(a) = y$. But $a \not\in L_y$.
A contradiction.
This proves Claim 2.

Let $U_y \in \rfs{Nbr}(y)$ be such that $f\inverse(U_y) \subseteq L_y$.
We check that the induction hypotheses (1)\,-\,(4) hold for $L_y$
and $U_y$.
The definition of $L_y$ implies that (1) holds,
and the definition of $U_y$ implies that (3) holds.
$$
\mbox{$
L_y \subseteq M_y \subseteq f\inverse(X - f(K)) \subseteq
f\inverse(X - f(\bigcup_{x \in A_0} L_x)) \subseteq
L - \bigcup_{x \in A_0} L_x.
$}
$$
So (2) holds.

We prove (4).
Certainly, $y \in f(L_y)$.
Recall that $L_y \subseteq M_y \subseteq f\inverse(W_y)$.
So $\rfs{cl}(L_y) \subseteq \rfs{cl}(f\inverse(W_y))$.
Also, $\rfs{cl}(f\inverse(W_y)) \subseteq f\inverse(\rfs{cl}(W_y))$.
So $\rfs{cl}(L_y) \subseteq f\inverse(\rfs{cl}(W_y))$
and hence
$$
f(\rfs{cl}(L_y)) \cap A \subseteq f(f\inverse(\rfs{cl}(W_y)) \cap A =
\rfs{cl}(W_y) \cap A = \sngltn{y}.
$$
The first equality follows from the surjectiveness of $f$.
This shows that (4) is fulfilled.
We have completed the inductive construction.

We show that $\calU \eqdf \setm{U_x}{x \in A}$
is a tight Hausdorff system for $A$.
Let $x,y \in A$ be distinct. Then $f\inverse(U_x) \subseteq L_x$
and $f\inverse(U_y) \subseteq L_y$. Since $L_x,L_y$ are disjoint,
so are $U_x$ and $U_y$.
Observe the following fact.
\begin{list}{}
{\setlength{\leftmargin}{23pt}
\setlength{\labelsep}{08pt}
\setlength{\labelwidth}{20pt}
\setlength{\itemindent}{-00pt}
\addtolength{\topsep}{-01pt}
\addtolength{\parskip}{-12pt}
\addtolength{\itemsep}{-05pt}
}
\item[$(*)$] 
Let $W,Z$ be compact Hausdorff spaces,
\hbox{$\fnn{h}{W}{Z}$} be continuous
and $\calC \subseteq \calP(W)$.
If $\setm{h(C)}{C \in \calC}$ is a pairwise disjoint family,
then $h(\rfs{acc}(\calC)) = \rfs{acc}(h(\calC))$.
\end{list}

\kern-1.5mm

\noindent
To see this, note that the fact
$h(\rfs{acc}(\calC)) \subseteq \rfs{acc}(h(\calC))$
holds even without assuming that $W$ and $Z$ are compact.
Now, the sets $h(\rfs{acc}(\calC))$ and $\rfs{acc}(h(\calC))$
are closed, and it is easy to see that
$h(\rfs{acc}(\calC))$ is dense in $\rfs{acc}(h(\calC))$.
So these sets must be equal and hence $(*)$ holds.

Let $A' \subseteq A$ and suppose that $B = \setm{y_x}{x \in A'}$,
where $y_x \in U_x$ for every $x \in A'$.
We show that $\rfs{acc}(B) = \rfs{acc}(A')$.
For every $x \in A'$ let $w_x \in L_x \cap f\inverse(y_x)$.
Such a choice is possible since $y_x \in U_x \subseteq f(L_x)$.
Let $I_x \in \calI(L_x)$ be such that $w_x \in I_x$.
Then by the definition of $L_x$
there is $z_x \in I_x$ such that $f(z_x) = x$.
Set $M = \setm{z_x}{x \in A'}$ and $N = \setm{w_x}{x \in A'}$.
Then $f(M) = A'$ and $f(N) = B$.
Applying  $(*)$ to $ \calC \eqdf \setm{\sngltn{m}}{m \in M}$,
we conclude that
(i) $\rfs{acc}(A') = f(\rfs{acc}(M))$.
Similarly,
(ii) $\rfs{acc}(B) = f(\rfs{acc}(N))$.
It is also clear that (iii) $\rfs{acc}(M) = \rfs{acc}(N)$.
To see this let $z \in \rfs{acc}(M)$. Then, without loss of generality,
there is a strictly increasing sequence
$\setm{z_{x_i}}{i < \mu} \subseteq M$ which converges to $z$.
This implies that $\setm{I_{x_i}}{i < \mu}$
is a strictly increasing sequence converging to $z$,
and so $\setm{w_{x_i}}{i < \mu}$ converges to $z$.
That is, $z \in \rfs{acc}(N)$.
We have shown that $\rfs{acc}(M) \subseteq \rfs{acc}(N)$,
and the same argument proves that
$\rfs{acc}(N) \subseteq \rfs{acc}(M)$.
So (iii) holds.

From (i) \,-\, (iii) it follows that $\rfs{acc}(B) = \rfs{acc}(A')$.

\smallskip

We prove that if $A' \subseteq A$,
then $\rfs{acc}(\setm{U_x}{x \in A'}) \subseteq \rfs{acc}(A')$.
For $x \in A'$ let $V_x = f\inverse(U_x)$.
By $(*)$,
$f(\rfs{acc}(\setm{V_x}{x \in A'})) = \rfs{acc}(\setm{U_x}{x \in A'})$.
Let $y \in \rfs{acc}(\setm{U_x}{x \in A'})$.
So there is $z \in \rfs{acc}(\setm{V_x}{x \in A'})$
such that $y = f(z)$.
There are a $\onetoone$ sequence $\setm{x_i}{i < \mu} \subseteq A'$
and a strictly monotonic sequence
$\setm{z_i}{i < \mu}$ such that $z_i \in V_{x_i}$ and
$\lim_{i < \mu} z_i = z$.
By the construction, $V_x = f\inverse(U_x) \subseteq L_x$,
so for every $i < \mu$ there is $I_i \in \calI(L_{x_i})$ such that
$z_i \in I_i$. Let $w_i \in f\inverse(x_i) \cap I_i$.
(The definition of the $L_x$'s assures the existence of~$w_i$).
Since the $I_i$'s are pairwise disjoint
and since $z_i,w_i \in I_i$,
it follows that $\lim_{i < \mu} w_i = \lim_{i < \mu} z_i$.
Hence $\lim_{i < \mu} x_i = \lim_{i < \mu} f(w_i) =
\lim_{i < \mu} f(z_i)  = f(z) = y$.
So $y \in \rfs{acc}(A')$.
We have proved the following facts.
\begin{itemize}
\addtolength{\parskip}{-11pt}
\addtolength{\itemsep}{06pt}
\item[(1)] 
If
$A' \subseteq A$ and $\setm{y_x}{x \in A'}$ is such that $y_x \in U_x$
for every $x \in A'$,
then $\rfs{acc}(A') = \rfs{acc}(\setm{y_x}{x \in A'})$.
\item[(2)] 
For every $A' \subseteq A$,
$\rfs{acc}(\setm{U_x}{x \in A'}) \subseteq \rfs{acc}(A')$.
\end{itemize}

Facts (1) and (2) imply that $\setm{U_x}{x \in A}$ is a tight family.
So $A$ has a tight Hausdorff system.
\smallskip\hfill\proofend

For a Hausdorff space $X$ denote by $\rfs{Is}(X)$ the set of
isolated points of $X$ and set $\rfs{D}(X) = X - \rfs{Is}(X)$.
\index{N@$\rfs{Is}(X)$. The set of isolated points of $X$}
\index{N@$\rfs{D}(X) = X - \rfs{Is}(X)$}
\index{N@$\rfs{D}_{\alpha}(X)$. The $\alpha$'s derivative of $X$}
Now define the {\it $\alpha$'s derivative} of $X$ as follows.
$\rfs{D}_0(X) = X$,
$\rfs{D}_{\alpha + 1}(X) = \rfs{D}(\rfs{D}_{\alpha}(X))$
and
$\rfs{D}_{\delta}(X) = \bigcap_{\alpha < \delta} \rfs{D}_{\alpha}(X)$
when $\delta$ is a limit ordinal.
Suppose now that {\em $X$ is a compact Hausdorff space.}
The {\it rank} of $X$ is
the first ordinal $\alpha$ such that
$\rfs{D}_{\alpha}(X)$ is finite or perfect.
(A set is {\it perfect}
if it does not have isolated points in its relative topology).
Denote the rank of $X$ by $\rfs{rk}(X)$.
Define $\rfs{ker}(X) = \rfs{D}_{\srfs{rk}(X) + 1}(X)$
and call $\rfs{ker}(X)$ the {\it perfect kernel} of $X$.
Hence $\rfs{ker}(X)$ is either the empty set
or an infinite perfect set.
It is easy to check that $X$ is scattered
iff $\rfs{ker}(X) = \emptyset$.
\index{D@perfect set. A set which does not have isolated points in its
       \\\indent relative topology}
\index{D@derivative.}
\index{D@rank. The rank of $X$, the first ordinal $\alpha$ such that
       $\rfs{D}_{\alpha}(X)$ is finite\\\indent
       or perfect}
\index{N@$\rfs{rk}(X)$. The rank of $X$}
\index{N@$\rfs{ker}(X) = \rfs{D}_{\srfs{rk}(X)}(X)$. The maximal perfect
       subset of $X$}
\index{D@perfect kernel}

Let $\rfs{Clop}(X)$ and $\rfs{Clsd}(X)$
denote respectively
the set of clopen subsets,
and the set of closed subsets
of a general Hausdorff space $X$
and set\break
$\rfs{Po}(X) =
\setm{x \in X}{\mbox{there is } U \in \rfs{Nbr}(x)
\mbox{ such that } \rfs{Is}(U) = \emptyset}$.
\index{N@$\rfs{Clop}(X)$. The set of clopen subsets of $X$}
\index{N@$\rfs{Clsd}(X)$. The set of closed subsets of $X$}
\index{N@$\rfs{Po}(X) = \setm{x \in X}{\mbox{there is }
       U \in \rfs{Nbr}(x) \mbox{ such that $\rfs{cl}(U)$ is perfect}}$}
{\em For a compact Hausdorff space} define
$\calS(X) = \setm{F \in \rfs{Clsd}(\rfs{ker}(X))}{
F \mbox{ is scattered}}$
and $\itOmega(X) = \rfs{sup}(\setm{\rfs{rk}(F)}{F \in \calS(X)})$.
\smallskip
\index{N@$\calS(X) =
       \setm{F \in \rfs{Clsd}(\rfs{ker}(X))}{F \mbox{ is scattered}}$}
\index{N@$\itOmega(X) = \rfs{sup}(\setm{\rfs{rk}(F)}{F \in \calS(X)})$}

The proof of Theorem \ref{t1.2} is by way of contradiction.
In the end of this section we assume that $X$ is a counter-example
to the theorem, and conclude that
$(2^{\abs{\itOmega(X)}})^+ < \abs{\rfs{ker}(X)}$,
which turns out to be a contradiction.
The proof is divided to a series of subclaims,
the first of which is the following statement
about members $X$ of $K_{\srfs{CII}}$.

\kern2mm

\noindent
{\bf 2. If ${\bf X}$ is a CO space, then
$\bfOmega(\bf X)$ is not attained by any
member\break
\rule{14pt}{0pt} of~$\calS(\bf X)$.}
\smallskip

\noindent
{\thickmuskip=2mu \medmuskip=1mu \thinmuskip=1mu 
Let $F \in \calS(X)$ and $A \subseteq \rfs{Is}(F)$ be such that
$\rfs{D}(F) = \rfs{acc}(A) = \rfs{acc}(\rfs{Is}(F) - A)$.}
Let $\setm{U_x}{x \in \rfs{Is}(F)}$ be a Hausdorff system
for $\rfs{Is}(F)$ and define
$F_1 = \rfs{cl}(F \cup \bigcup \setm{U_x \cap \rfs{ker}(X)}{x \in A})$.
The set $F_1$ is called a fattening of $F$.
The precise definition of a fattening is given below.

\begin{defn}\label{d2.4}
\begin{rm}
Let $X$ be a Hausdorff compact space,
$F \in \calS(X)$,
$F_1 \in \rfs{Clsd}(\rfs{ker}(X))$ and $F \subseteq F_1$.
We call $F_1$ a {\it fattening} of $F$ if the following holds.
\begin{itemize}
\addtolength{\parskip}{-11pt}
\addtolength{\itemsep}{06pt}
\item[(F1)] 
$\rfs{Is}(F) = \rfs{Is}(F_1) \cup (\rfs{Is}(F) \cap \rfs{Po}(F_1))$.
\item[(F2)] 
$\rfs{D}(F) = \rfs{acc}(\rfs{Is}(F_1))
= \rfs{acc}(\rfs{Is}(F) \cap \rfs{Po}(F_1))$.
\end{itemize}
\end{rm}
\end{defn}

\begin{prop}\label{p2.5}
If $X$ is collectionwise Hausdorff and compact and $F \in \calS(X)$,
then $F$ has a fattening.
\end{prop}

\noindent{\bf Proof }
Let $\setm{U_x}{x \in \rfs{Is}(F)}$
and
$\setm{V_y}{y \in \rfs{Is}(\rfs{D}(F))}$
be Hausdorff systems for $\rfs{Is}(F)$
and $\rfs{Is}(\rfs{D}(F))$ respectively.
For every $y \in \rfs{Is}(\rfs{D}(F))$
let $A_y$ be an infinite subset of $\rfs{Is}(F) \cap V_y$ such that
$(\rfs{Is}(F) \cap V_y) - A_y$ is also infinite.
Then
$\rfs{acc}(A_y) =
\rfs{acc}((\rfs{Is}(F) \cap V_y) - A_y) = \sngltn{y}$.
Let
$$
F_1 =
\rfs{cl}
\left(F \cup \bigcup \setmlr{U_x \cap \rfs{ker}(X)}
{x \in \mbox{$\bigcup_{y \in \srfs{Is}(\srfs{D}(F))}$} A_y}\right).
$$
It is left the reader to check that $F_1$ is as required.
\smallskip\hfill\proofend

\begin{prop}\label{p2.6}
Let $X$ be a collectionwise Hausdorff CO space
and $E \in \calS(X) - \sngltn{\emptyset}$.
Suppose that $\rfs{rk}(E) = \alpha$.
Then for every $n \in \omega$ there is $F \in \calS(X)$ such that
$\abs{\rfs{D}_{\alpha}(F)} = n$.
\end{prop}

\noindent{\bf Proof }
For $\alpha = 0$ the claim of the proposition follows from the fact
that $\rfs{ker}(X)$ is infinite,
so we assume that $\alpha > 0$.
The proof is by induction on $n$.
We may assume that $\abs{\rfs{D}_{\alpha}(E)} = 1$.
Suppose that $F \in \calS(X)$ and $\abs{\rfs{D}_{\alpha}(F)} = n$.
We show that $(*)$ there is $G \in \calS(X)$ such that
$\abs{\rfs{D}_{\alpha}(G)} = 2n$.

Let $\whatH$ be a fattening of $F$,
let $H \in \rfs{Clop}(X)$ and $\varphi$
be such that $\iso{\varphi}{\whatH}{H}$.
Set $\whatH^0 = \rfs{cl}(\rfs{Is}(F) \cap \rfs{Po}(\whatH))$
and $H^0 = \varphi(\whatH^0)$.
Clearly, 
$$
D(\whatH^0) = D(\rfs{cl}(\rfs{Is}(F) \cap \rfs{Po}(\whatH)))
= \rfs{acc}(\rfs{Is}(F) \cap \rfs{Po}(\whatH)),
$$
and by (F2),
$$
\rfs{acc}(\rfs{Is}(F) \cap \rfs{Po}(\whatH)) =
\rfs{acc}(\rfs{Is}(\whatH)).
$$
So
$$
D(\whatH^0) = \rfs{acc}(\rfs{Is}(\whatH)).
$$
The same holds for $H^0$ and $H$, namely,
$$
D(H^0) = \rfs{acc}(\rfs{Is}(H)).
$$
Since $H$ is clopen in $X$, it follows that
$\rfs{Is}(H) = \rfs{Is}(X) \cap H$. So
$$
D(H^0) = \rfs{acc}(\rfs{Is}(X) \cap H),
$$
and hence
$$
D(H^0) \cap \rfs{Po}(X) = \emptyset.
$$

Let $K_0$ be a clopen subset of $X$
homeomorphic to $\rfs{ker}(X)$.
Then there is $F^0 \subseteq K_0$ such that $F^0 \cong F$.
We shall show that $F^0 \cup H^0$ is the set $G$ required in $(*)$.
Since $F^0 \subseteq K_0$ and $K_0$ is open and perfect,
we have that $F^0 \subseteq \rfs{Po}(X)$.
So $\rfs{D}(H^0) \cap F^0 = \emptyset$.
This implies that
$H^0 \cap F^0$ is a finite subset of $\rfs{Is}(H^0)$,
and hence
\begin{itemize}
\addtolength{\parskip}{-11pt}
\addtolength{\itemsep}{06pt}
\item[(1)] 
$\rfs{D}(F^0 \cup H^0)$ is the disjoint union of $\rfs{D}(F^0)$
and $\rfs{D}(H^0)$.
\end{itemize}

We next show that
\begin{itemize}
\addtolength{\parskip}{-11pt}
\addtolength{\itemsep}{06pt}
\item[(2)] 
$H^0 \subseteq \rfs{ker}(X)$.
\end{itemize}
Recall that $\whatH^0 = \rfs{cl}(\rfs{Is}(F) \cap \rfs{Po}(\whatH))$.
So $\rfs{Is}(\whatH^0) = \rfs{Is}(F) \cap \rfs{Po}(\whatH)$
and hence $\rfs{Is}(\whatH^0) \subseteq \rfs{Po}(\whatH)$.
It follows that $\rfs{Is}(\whatH^0) \subseteq \rfs{ker}(\whatH)$,
and this implies that
$\rfs{cl}(\rfs{Is}(\whatH^0)) \subseteq \rfs{ker}(\whatH)$.
But $\rfs{cl}(\rfs{Is}(\whatH^0)) = \whatH^0$,
so $\whatH^0 \subseteq \rfs{ker}(\whatH)$.
Since $\varphi$ takes $\whatH$ to $H$ and 
$\whatH^0$ to $H^0$, it follows that
$H^0 \subseteq \rfs{ker}(H)$.
So $H^0 \subseteq \rfs{ker}(X)$.

Now we show that
\begin{itemize}
\addtolength{\parskip}{-11pt}
\addtolength{\itemsep}{06pt}
\item[(3)] 
$\abs{\rfs{D}_{\alpha}(H^0)} = n$.
\vspace{-2mm}
\end{itemize}

Clearly,
$\rfs{D}(\whatH^0) = \rfs{acc}(\rfs{Is}(\whatH^0)) =
\rfs{acc}(\rfs{Po}(\whatH) \cap \rfs{Is}(F))$.
By (F2),\break
$\rfs{acc}(\rfs{Po}(\whatH) \cap \rfs{Is}(F)) = \rfs{D}(F)$.
So $\rfs{D}(\whatH^0) = \rfs{D}(F)$.
Since $\alpha > 0$,
it follows that $\rfs{D}_{\alpha}(\whatH^0) = \rfs{D}_{\alpha}(F)$.
So
$\abs{\rfs{D}_{\alpha}(H^0)} = \abs{\rfs{D}_{\alpha}(\whatH^0)} =
\abs{\rfs{D}_{\alpha}(F)} = n$.
We have proved (3).

Recall that $F^0 \subseteq \rfs{Po}(X)$. So
\begin{itemize}
\addtolength{\parskip}{-11pt}
\addtolength{\itemsep}{06pt}
\item[(4)] 
$F^0 \subseteq \rfs{ker}(X)$.
\vspace{-2mm}
\end{itemize}

From (2) and (4) it follows that, $F^0 \cup H^0 \in \calS(X)$.
Since $F^0 \cong F$, $\abs{D_{\alpha}(F^0)} = n$.
Hence by Facts (1) and (3),
$\abs{\rfs{D}_{\alpha}(F^0 \cup H^0)} = 2n$.
\smallskip\hfill\proofend

Any continuous image of a sequentially compact space
is sequentially compact. So we have the following fact.

\begin{prop}\label{p2.7}
Every member of $K_{\srfs{CII}}$ is sequentially compact.
\end{prop}

For scattered spaces $F$ and $G$ define $F \newprec G$,
if either $\rfs{rk}(F) < \rfs{rk}(G)$
or $\rfs{rk}(F) = \rfs{rk}(G)$ and
$\abs{\rfs{D}_{\srfs{rk}(F)}(F)} < \abs{\rfs{D}_{\srfs{rk}(F)}(G)}$.
Let $X$ be a Hausdorff space. We say that $X$ is
{\it strongly Hausdorff for convergent sequences\,}
if every $\onetoone$ convergent sequence in $X$ has
a strong Hausdorff system.
\index{D@strongly Hausdorff for convergent sequences}

\begin{prop}\label{p2.8}
Let $X$ be a sequentially compact space which is
strongly Hausdorff for convergent sequences.
Suppose that $F_0 \newprec F_1 \newprec \ldots$ is a sequence of members
of $\calS(X)$. Then there is $F \in \calS(X)$ such that
$F_i \newprec F$ for every $i \in \omega$.
\end{prop}

\noindent{\bf Proof }
We may assume that $\abs{D_{\srfs{rk}(F_i)}(F_i)} \geq i$.
More precisely, there is a subsequence $\setm{m_i}{i \in \omega}$
and for every $i$ there is a closed subset $\whatF_i \subseteq F_{m_i}$
with the property that $F_i \newprec \whatF_{i + 1}$
and $\abs{D_{\srfs{rk}(\whatF_i)}(\whatF_i)} \geq i$
for every $i \in \omega$.
To see this we distinguish between the cases:
(i) $\setm{\rfs{rk}(F_i)}{i \in \omega}$ is eventually constant,
and
(ii) $\setm{\rfs{rk}(F_i)}{i \in \omega}$ is not eventually constant.
If (i) happens we define $\whatF_i = F_{i + n_0}$,
where $n_0$ is such that for every $i,j \geq n_0$,
$\rfs{rk}(F_i) = \rfs{rk}(F_j)$.
Suppose that (ii) happens. Then we take a subsequence
$\sngltn{m_i}_{i \in \omega}$ such that for every~$i$,
$\rfs{rk}(F_{m_{i + 1}}) > \rfs{rk}(F_{m_i}) + 1$.
Let $\setm{\alpha_i}{i \in \omega}$ be such that
$\lim_i \alpha_i = \lim_i \rfs{rk}(F_i)$ and for every $i$,
$\rfs{rk}(F_{m_i}) < \alpha_i < \rfs{rk}(F_{m_{i + 1}})$.
Let $\whatF_i$ be a closed subset of $F_{m_i}$ such that
$\abs{D_{\alpha_i)}(\whatF_i)} = i$.
Hence $\setm{\whatF_i}{i \in \omega}$ is as desired.

It follows that there is a $\onetoone$ sequence
$\setm{x_i}{i \in \omega}$ such that
$x_i \in D_{\srfs{rk}(\whatF_i)}(\whatF_i)$ for every $i \in \omega$.
We may assume that $\sngltn{x_i}$ is a convergent sequence.
So by the sequential strong Hausdorff property of $X$,
$\sngltn{x_i}_{i \in \omega}$ has a strong Hausdorff system
$\sngltn{U_i}_{i \in \omega}$.
Let $F = \rfs{cl}(\bigcup_{i \in \omega} \whatF_i \cap U_i)$.
Let $x = \rfs{lim}_{i \in \omega} x_i$
and $\alpha = \rfs{Sup}_{i \in \omega}\rfs{rk}(\whatF_i)$.
It easy to see that $\rfs{D}_{\alpha}(F) = \sngltn{x}$
and clearly, $F \subseteq \calS(X)$.
So for every $i \in \omega$, $F_i \newprec \whatF_{i + 1} \newprec F$.
\smallskip\hfill\proofend

\begin{cor}\label{c2.9}
Let $X$ be a CO space with a nonempty kernel.
Suppose also that $X$ is collectionwise Hausdorff,
strongly Hausdorff for convergent sequences
and sequentially compact.
Then for every $F \in \calS(X)$, $\rfs{rk}(F) < \itOmega(X)$.
\end{cor}

\noindent{\bf Proof }
Suppose by contradiction that
$F \in \calS(X)$ and $\rfs{rk}(F) = \itOmega(X)$.
By Proposition \ref{p2.6},
there is a sequence $F_1,F_2,\ldots$ of members of $\calS(X)$
such that $\abs{\rfs{D}_{\itOmega(X)}(F_i)} = i$.
By Proposition \ref{p2.8}, there is $H \in \calS(X)$ such that
$\rfs{rk}(H) > \itOmega(X)$.
A contradiction, so $\rfs{rk}(F) < \itOmega(X)$
for every $F \in \calS(X)$.
\smallskip\hfill\proofend

If $\rfs{ker}(X) \neq \emptyset$
and $\rfs{rk}(F) < \itOmega(X)$ for every $F \in \calS(X)$,
then we say that {\it $\itOmega(X)$ is not attained in $X$}.
\index{D@attained: $\itOmega(X)$ is not attained in $X$}

We next define the notion of a good point
and prove the following statement for non-scattered members $X$ of
$K_{\srfs{CII}}$.

\kern2mm

\noindent
{\bf 3. If ${\bf\Omega(X)}$ is not attained in ${\bf X}$,
then the set of good points of ${\bf X}$ is\break
\rule{12.5pt}{0pt} perfect.}
\smallskip

\noindent
Let $X$ be a compact Hausdorff space.
A member $x \in X$ is called a {\it good point} of $X$
if for every $\alpha < \itOmega(X)$ and $U \in \rfs{Nbr}(x)$
there is $F \in \calS(X)$
such that $\rfs{D}_{\alpha}(F) \cap U \neq \emptyset$.
\index{D@good point}
Note that if $\itOmega(X)$ is not attained, then it is a limit ordinal.

We shall show that if $\itOmega(X)$ is not attained,
then the set of good points is a nonempty perfect set.
The existence of a good point is a trivial consequence of
the compactness of $X$.
It is also trivial that the set of good points is closed.
So we have the following fact.

\begin{prop}\label{p2.11}
Let $X$ be a compact Hausdorff space with a nonempty kernel.
Then the set of good points of $X$ is closed and nonempty.
\end{prop}

The following proposition is well-known and easy to prove.
Recall that according to our definition of scatterednes,
a scattered space is compact Hausdorff.

\begin{prop}\label{p2.12}
Let $Y$ be a scattered space, $X$ be a Hausdorff space
and $\fnn{g}{Y}{X}$ be a continuous surjective function.
Then $X$ is scattered.
\end{prop}

There is another property of members of $K_{\srfs{CII}}$
that we now establish.
Let $\lambda$ be an infinite cardinal, $A \subseteq X$ and $x \in X$.
Call $x$ a {\it $\lambda$-accumulation point} of $A$ if
$\abs{U \cap A} = \lambda$ for every $U \in \rfs{Nbr}(x)$.
\index{D@accumulation point: $\lambda$-accumulation point.
       A point $x$ is a\\\indent
       $\lambda$-accumulation point of $A$ if $\abs{U \cap A} = \lambda$
       for every $U \in \rfs{Nbr}(x)$}
\index{D@dense@@dense: $\lambda$-dense linear ordering}
A linear ordering $\pair{L}{<}$ is {\it $\lambda$-dense},
if $\abs{L} > 1$ and for every $a < b$ in $L$, $\abs{(a,b)} = \lambda$.

\begin{prop}\label{p2.13}
\num{a} Let $\lambda$ be an infinite regular cardinal
and $\pair{L}{<}$ be a linear ordering of power $\lambda$.
Then either $L$ has a subset of order type $\lambda$ or $\lambda^*$,
or $L$ has a $\lambda$-dense subset.

\num{b}
Let $\alpha$ be a successor ordinal equipped with its order topology
and\break
$\fnn{g}{\alpha}{X}$ be a continuous surjective function.
Then $\abs{\rfs{Is}(X)} = \abs{X}$.

\num{c} Let $X \in K_{\srfs{CII}}$.
Suppose that $A \subseteq X$ and $\lambda \eqdf \abs{A}$
is an infinite regular cardinal. Then either
there is $B \subseteq A$ such that $\abs{B} = \lambda$,
$B$ is relatively discrete and $\rfs{cl}(B)$ is scattered,
or $A$ has at least two $\lambda$-accumulation points.
\end{prop}

\noindent{\bf Proof }
(a) Define an equivalence relation on $L$ as follows: $a \sim b$ if
the open interval whose endpoints are $a$ and $b$ has cardinality
$< \lambda$.
If there is an equivalence class of cardinality $\lambda$,
then that equivalence class contains
an increasing or decreasing sequence of type $\lambda$.
If every equivalence class has cardinality $< \lambda$,
then the chain of equivalence classes is $\lambda$-dense.
\smallskip

(b) For every
$x \in \rfs{Is}(X)$
let $\beta_x \in g\inverse(x)$.
Define $B = \setm{\beta_x}{x \in \rfs{Is}(X)}$ and $C = \rfs{cl}(B)$.
Then
\newline
(1)
\centerline{
$g(\rfs{cl}(B)) = \rfs{cl}(g(B))$.
}
As $X$ is a continuous image of a scattered space,
$X$ must be scattered.
In particular, $\rfs{Is}(X)$ is dense in $X$.
So from (1) we conclude that $g(C) = X$.
Either $C$ or $C$ minus its maximum have the same order type as $B$,
so in particular, $\abs{C} = \abs{B}$.
We thus have that
$\abs{X} \leq \abs{C} = \abs{B} = \abs{\rfs{Is}(X)}$.
\smallskip

(c) Suppose that $X$, $\lambda$ and $A$ are as in Part (c),
and let $\pair{L}{<}$ be a compact linear ordering and
$\fnn{h}{L}{X}$ be continuous and surjective.
For every $a \in A$ let $\ell_a \in h\inverse(a)$ and let
$M = \setm{\ell_a}{a \in A}$.
By Part (a), either (i) $M$ contains
an increasing or decreasing sequence of type $\lambda$
or (ii) $M$ contains a $\lambda$-dense subset.
Assume first that (i) happens. We may then assume that $M$ is an
increasing sequence of type $\lambda$. Let $N = \rfs{cl}^L(M)$.
Then $N$ is a compact interval space, and $N$ is scattered.
Let $C = h(N)$.
By Proposition~\ref{p2.12}, $C$ is a scattered space.
If $b \in C - A$, then there is $n \in \rfs{cl}(M) - M$ such that
$h(n) = b$. Hence since $h {\restriction} M$ is $\onetoone$,
$b \in \rfs{acc}(A)$.
It follows that $\rfs{Is}(C) \subseteq A$.
By Part (b), $\abs{\rfs{Is}(C)} = \lambda$.
Hence $\rfs{Is}(C)$ is a relatively discrete subset of $A$
of cardinality $\lambda$ and its closure is scattered.
That is, $B \eqdf \rfs{Is}(C)$ fulfills the requirements
of the proposition.

Suppose next that (ii) happens.
We may assume that $M$ is $\lambda$-dense.
It is trivial that if $n \in \rfs{acc}(M)$,
then $n$ is a $\lambda$-accumulation point of $M$.
Since $h {\restriction} M$ is $\onetoone$,
$h(n)$ is then a $\lambda$-accumulation point of $h(M) = A$.
We have thus verified that
\begin{itemize}
\addtolength{\parskip}{-11pt}
\addtolength{\itemsep}{06pt}
\item[(1)] 
for every $n \in \rfs{acc}(M)$,
$h(n)$ is a $\lambda$-accumulation point of $A$.
\end{itemize}

\kern-1mm

\noindent
Next we notice that
\begin{itemize}
\addtolength{\parskip}{-11pt}
\addtolength{\itemsep}{06pt}
\item[(2)] 
$\rfs{acc}(A) \subseteq h(\rfs{acc}(M))$.
\end{itemize}

\kern-1mm

{\thickmuskip=3.5mu \medmuskip=2mu \thinmuskip=1mu 
\noindent
Let $a \in \rfs{acc}(A)$.
Set $A' = A - \sngltn{a}$ and $M' = M - h\inverse(\sngltn{a})$.
Then $h(M') = A'$}
and hence $\rfs{cl}(A') = h(\rfs{cl}(M'))$.
Clearly, $a \in \rfs{acc}(A') \subseteq \rfs{cl}(A')$
and thus $a \in h(\rfs{cl}(M'))$.
But $a \not\in h(M')$.
So $a \in h(\rfs{acc}(M')) \subseteq h(\rfs{acc}(M))$.

It follows from (1) and (2) that every accumulation point of $A$
is a $\lambda$-accumulation point.
If $A$ has at least two accumulation points,
then the requirements of the proposition are fulfilled.
Otherwise $A$ has exactly one accumulation point, which means that
$\rfs{cl}(A)$ is homeomorphic to the one point compactification
of a discrete space of cardinality $\lambda$.
If this happens, then $A - \rfs{acc}(A)$ is relatively discrete
and $\rfs{cl}(A - \rfs{acc}(A))$ is scattered.
So the requirements of the proposition are again fulfilled.
\smallskip\hfill\proofend

For a compact Hausdorff space $X$,
let $\rfs{Good}(X)$ denote the set of good points of $X$.
\index{N@$\rfs{Good}(X)$. The set of good points of $X$}
If $X$ is a scattered space and $x \in X$,
then the {\it rank of $x$ in $X$} is defined to be
$\rfs{max}(\setm{\alpha}{x \in \rfs{D}_{\alpha}(X)})$.
The rank of $x$ in $X$ is denoted by $\rfs{rk}^{X}(x)$.
Note that if $F \in \rfs{Nbr}^X_{\srfs{cl}}(x)$,
then $\rfs{rk}^{F}(x) = \rfs{rk}^{X}(x)$.
\index{D@rank of $x$ in $X$}
\index{N@$\rfs{rk}^{X}(x) =
       \rfs{max}(\setm{\alpha}{x \in \rfs{D}_{\alpha}(X)})$}

\begin{prop}\label{p2.14}
Let $X$ be a strongly Hausdorff compact space.
Suppose also that
\begin{itemize}
\addtolength{\parskip}{-11pt}
\addtolength{\itemsep}{06pt}
\item[$(*)$] 
for every $A \subseteq X$,
if $\lambda \eqdf \abs{A}$ is an infinite regular cardinal,
then either $A$ has at least two $\lambda$-accumulation points,
or there is $B \subseteq A$ such that $\abs{B} = \lambda$,
$B$ is relatively discrete and $\rfs{cl}(B)$ is scattered.
\vspace{-2mm}
\end{itemize}
Assume further that $\rfs{ker}(X) \neq \emptyset$
and that $\itOmega(X)$ is not attained.
\underline{Then} $\rfs{Good}(X)$ is a nonempty perfect set.
\end{prop}

\noindent{\bf Proof }
By Proposition~\ref{p2.11}, $\rfs{Good}(X) \neq \emptyset$.
So suppose by contradiction that $x$ is a good point of $X$,
$U' \in \rfs{Nbr}(x)$
and $\rfs{Good}(X) \cap U' = \sngltn{x}$.
Let $U \in  \rfs{Nbr}(x)$ be such that $\rfs{cl}(U) \subseteq U'$.
Since $\itOmega(X)$ is not attained,
$\itOmega(X)$ is a limit ordinal.
Let $\lambda = \rfs{cf}(\itOmega(X))$
and $\setm{\alpha_i}{i < \lambda}$ be a strictly increasing sequence
of ordinals converging to $\itOmega(X)$.
For every $i < \lambda$ let $F_i \in \calS(X)$ be such that
$\rfs{rk}(F_i) = \alpha_i$ and $F_i \subseteq U$ and choose
$x_i \in \rfs{D}_{\alpha_i}(F_i)$.
It may happen that $\abs{\setm{x_i}{i < \lambda}} < \lambda$.
In that case we show that $\setm{x_i}{i < \lambda}$ can be replaced
by another sequence $\setm{x'_i}{i < \lambda}$ such that
$\abs{\setm{x'_i}{i < \lambda}} = \lambda$.
So suppose that $\abs{\setm{x_i}{i < \lambda}} < \lambda$.
Then by taking a subsequence we may assume that $x_i = x_j$ for every
$i,j < \lambda$. Then $x_0$ is a good point,
and since $x$ is the only good
point in $U$, $x_0 = x$.
For every $i < \lambda$ let
$x'_i \in
\rfs{D}_{\alpha_i}(F_{i + 1}) - \rfs{D}_{\alpha_i + 1}(F_{i + 1})$.
So $x'_i \neq x_{i + 1} = x$.
It therefore follows from the above argument that for every $i$,
$\abs{\setm{j}{x'_j = x'_i}} < \lambda$.
So $\abs{\setm{x'_i}{i < \lambda}} = \lambda$.
We may thus assume that $\setm{x_i}{i < \lambda}$ is $\onetoone$.
We apply $(*)$ to $A \eqdf \setm{x_i}{i < \lambda}$.
Every $\lambda$-accumulation point of $A$
is a good point and it belongs to $U'$.
So since $\abs{U \cap \rfs{Good}(X)} = 1$,
it follows that $A$ has at most one $\lambda$-accumulation point.
So by $(*)$ there is $B \subseteq A$ such that
$\abs{B} = \lambda$, $B$ is relatively discrete
and $\rfs{cl}(B)$ is scattered.
Let $B = \setm{x_{i(j)}}{j < \lambda}$
and $ \calU = \setm{U_j}{j < \lambda}$
be a strong Hausdorff system for $B$.
We may assume that $\rfs{cl}(U_j) \cap \rfs{cl}(U_{j'}) = \emptyset$
for every $j \neq j'$.
For every $j < \lambda$ define $\whatF_j = F_{i(j)} \cap \rfs{cl}(U_j)$
and define $F = \rfs{cl}(\bigcup_{j < \lambda} \whatF_j)$.
Clearly $F \subseteq \rfs{ker}(X)$.
Also, since $\rfs{cl}(B)$ is scattered
and $\calU$ is a strong Hausdorff system,
it follows that $F$ is scattered. So $F \in \calS(X)$.
For every $j < \lambda$,
$\rfs{rk}^F(x_{i(j)}) = \rfs{rk}^{\whatF_j}(x_{i(j)}) =
\rfs{rk}^{F_{i(j)}}(x_{i(j)}) = \alpha_{i(j)}$.
It follows that $\rfs{rk}(F) \geq \itOmega(X)$.
A contradiction to the non-attainment of $\itOmega(X)$.
\hfill\proofend

\kern2mm

\noindent
{\bf 4. Sets which code ordinals.}
\smallskip

\noindent
We shall define the notion of a code of an ordinal.
Codes are certain compact Hausdorff spaces which code ordinals.
Two other notions are to be defined:
the ``perfect end'' of a compact Hausdorff space
$F$ \,--\, this is a certain nonempty closed subset of $F$,
and the notion of a ``demonstrative subset'' of a compact Hausdorff
space $X$.
Two facts about codes are important:
\begin{itemize}
\addtolength{\parskip}{-11pt}
\addtolength{\itemsep}{06pt}
\item 
If $F$ and $H$ are codes for two different ordinals, $x$ belongs
to the perfect end of $F$, and $U \in \rfs{Nbr}(x)$,
then $U$ is not homeomorphic to any open subset of $H$,
(Proposition~\ref{p2.21}).
\item 
If $X \in K_{\srfs{CII}}$,
$\rfs{ker}(X) \neq \emptyset$, $\itOmega(X)$
is not attained and $0 < \alpha < \itOmega(X)$,
then any neighborhood of a member of $\rfs{Good}(X)$
contains a demonstrative set which is an $(\alpha + 1)$-code,
(Lemma~\ref{l2.24}).
\vspace{-2mm}
\end{itemize}

We first verify the following property of members of $K_{\srfs{CII}}$.

\begin{prop}\label{p2.16}
Let $X \in K_{\srfs{CII}}$ and $A \subseteq X$ be relatively discrete.
Suppose that $\rfs{cl}(A)$ is not scattered.
Then there is $B \subseteq A$ such that

\kern-5mm

$$
\rfs{acc}(B) = \rfs{ker}(\rfs{cl}(A)).
$$
\end{prop}

\noindent{\bf Proof }
Let $\pair{L}{<}$ be a linear ordering and $\fnn{g}{L}{X}$
be continuous and surjective.
For every $a \in A$ let $\ell_a \in g\inverse(a)$ and define
$L_0 = \rfs{cl}(\setm{\ell_a}{a \in A})$ and $g_0 = g {\restriction} L_0$.
So $\fnn{g_0}{L_0}{\rfs{cl}(A)}$ and $g_0\inverse(A)$ is topologically dense
in $L_0$.
We may thus assume that $\fnn{g}{L}{\rfs{cl}(A)}$
and $g\inverse(A)$ is topologically dense in $L$.

Let $U = L - g\inverse(\rfs{ker}(\rfs{cl}(A)))$.
Then $U$ is an open subset of $L$.
Let $\calI$ be the partition of $U$ into maximal convex subsets of $L$.
Then every member of $\calI$ is an open interval of $L$.
Since $g\inverse(A)$ is topologically dense in $L$,
$g\inverse(A) \cap I \neq \emptyset$ for every $I\in \calI$.
Choose $m_I \in I \cap g\inverse(A)$ and set
$M = \setm{m_I}{I \in \calI}$ and $B = g(M)$.
Clearly, $\rfs{acc}(M) \cap U = \emptyset$.
So $\rfs{cl}(M) \subseteq M \cup g\inverse(\rfs{ker}(\rfs{cl}(A)))$.
Since $L$ is compact, $g(\rfs{cl}(M)) = \rfs{cl}(B)$.
Hence $\rfs{cl}(B) \subseteq B \cup \rfs{ker}(\rfs{cl}(A))$.
Since $\rfs{acc}(B) \cap B = \emptyset$, it follows that
$\rfs{acc}(B) \subseteq \rfs{ker}(\rfs{cl}(A))$.

Let $x \in \rfs{ker}(\rfs{cl}(A))$, and assume by contradiction that
$x \not\in \rfs{acc}(B)$.
Since $A$ is relatively discrete, we have $x \not\in A$,
and in particular, $x \not\in B$. So $x \not\in \rfs{cl}(B)$.
This implies that $g\inverse(x) \cap \rfs{acc}(M) = \emptyset$.
Then every $y \in g\inverse(x)$ has an open neighborhood $W_y$
such that $\setm{I \in \calI}{I \cap W_y \neq \emptyset}$ is finite.
Note that $g\inverse(A) \subseteq U = \bigcup \calI$.
So $\bigcup \calI$ is dense in $L$.
Using the facts that $\bigcup \calI$ is dense in $L$ and that
$y \not\in \bigcup \calI$ we conclude that there are $I_y,J_y \in \calI$
such that $y$ is the right endpoint of $I_y$ and the left endpoint
of $J_y$.
Let $V_y = I_y \cup \sngltn{y} \cup J_y$.
Define $V = \bigcup \setm{V_y}{y \in g\inverse(x)}$.
Since $g\inverse(x) \subseteq V$ and $L$ is compact, it follows that
$g(V)$ is a neighborhood of $x$.
But
$V \subseteq (L - g\inverse(\rfs{ker}(\rfs{cl}(A)))) \cup g\inverse(x)$,
so $g(V) \cap \rfs{ker}(\rfs{cl}(A)) = \sngltn{x}$.
This means that $x$ is an isolated point of $\rfs{ker}(\rfs{cl}(A))$,
a contradiction.
It follows that $x \in \rfs{acc}(B)$,
so $\rfs{ker}(\rfs{cl}(A)) \subseteq \rfs{acc}(B)$.
\smallskip\hfill\proofend

\begin{prop}\label{p2.17}
Let $X$ be a compact Hausdorff perfect space.
Then there is a relatively discrete subset $A \subseteq X$ such that
$\rfs{cl}(A)$ contains a nonempty perfect set.
\end{prop}

\noindent{\bf Proof }
We show that $[0,1]$ is a continuous image of $X$.
If $X$ is not\break
$0$-dimensional let $x,y$ be two distinct points
in the same connected component of $X$
and $\fnn{g}{X}{[0,1]}$ be a continuous function such that
$g(x) = 0$ and $g(y) = 1$. Then $\rfs{Rng}(g) = [0,1]$.
If $X$ is $0$-dimensional and perfect,
then the Cantor set is a continuous image of $X$
and $[0,1]$ is a continuous image of the Cantor set.
So there is a continuous surjective $g$ from $X$ to $[0,1]$.

Let $B \subseteq [0,1]$ be a relatively discrete set such that
$\rfs{cl}(B) - B$ is perfect.
For every $b \in B$ choose $x_b \in g\inverse(b)$ and define
$A = \setm{x_b}{b \in B}$.
It follows from the relative discreteness of $B$ that $A$ is also
relatively discrete.
We have shown  in Proposition \ref{p2.12}
that a continuous image of a scattered space is scattered.
As $\rfs{cl}(B)$ is not scattered, $\rfs{cl}(A)$ cannot be scattered.
\smallskip\hfill\proofend

In order to define codes, we introduce the notion of the
{\it perfect derivative} of a compact Hausdorff space $X$.
For a compact Hausdorff space $X$ define

$$
\rfs{PD}(X) \eqdf X - \rfs{Is}(X) - \rfs{Po}(X),
$$

\kern-5mm

$$
\rfs{PD}_0(X) \eqdf X
\mbox{, \ }
\rfs{PD}_{\alpha + 1}(X) \eqdf \rfs{PD}(\rfs{PD}_{\alpha}(X))
$$

\noindent
and if $\delta$ is a limit ordinal, then define

\kern-3.1mm

$$
\mbox{$
\rfs{PD}_{\delta}(X) \eqdf
\bigcap_{\alpha < \delta} \rfs{PD}_{\alpha}(X).
$}
$$
The {\it perfect rank} of $X$ is defined by
$\rfs{prk}(X) =
\rfs{max}(\setm{\alpha}{\rfs{PD}_{\alpha}(X) \neq \emptyset})$
and the {\it perfect end} is defined by
$\rfs{Pend}(X) = \rfs{PD}_{\srfs{prk}(X)}(X)$.
Note that $\rfs{Pend}(X)$ is the union of a finite set
of isolated points and a perfect set. Each may be empty but not both.
\index{D@perfect derivative}
\index{N@$\rfs{PD}(X) = X - \rfs{Is}(X) - \rfs{Po}(X)$.
       The perfect derivative of $X$}
\index{N@$\rfs{PD}_{\alpha}(X)$. The $\alpha$'s perfect derivative
       of $X$}
\index{D@perfect rank}
\index{D@perfect end}
\index{N@$\rfs{prk}(X) =
       \rfs{max}(\setm{\alpha}{\rfs{PD}_{\alpha}(X) \neq \emptyset})$.
       The perfect rank of $X$}
\index{N@$\rfs{Pend}(X) = \rfs{PD}_{\srfs{prk}(X)}(X)$.
       The perfect end of $X$.}

For $x \in X$, the property:
``$x$ belongs to $\rfs{PD}_{\alpha}(X)$''
is a local property. This is expressed in the next observation
which is trivial and is not proved.

\begin{prop}\label{p2.19}
Let $Z$ be a compact Hausdorff space.
Suppose that\break
$U \subseteq G \subseteq Z$, $U$ is open and $G$ is closed.
Then $U \cap \rfs{PD}_{\alpha}(G) = U \cap \rfs{PD}_{\alpha}(Z)$
for every ordinal $\alpha$.
\end{prop}

\begin{defn}\label{d2.20}
\begin{rm}
Let $\alpha \geq 2$.
A compact Hausdorff space $F$ is called an {\it $\alpha$-code} if
\begin{itemize}
\addtolength{\parskip}{-11pt}
\addtolength{\itemsep}{06pt}
\item[(C1)] 
$\rfs{prk}(F) = \alpha$,
\item[(C2)] 
$\rfs{Po}(\rfs{PD}_{\beta}(F)) = \emptyset$
for every $0 < \beta < \alpha$,
\item[(C3)] 
$\rfs{Pend}(F)$ is perfect.
\end{itemize}
A set which is an $\alpha$-code for some ordinal $\alpha$,
is called a {\it code}.
\index{D@code: $\alpha$-code}
\index{D@code}
\hfill\myqed
\end{rm}
\end{defn}

\begin{prop}\label{p2.21}
Suppose that $F$ is an $\alpha$-code, $H$ is a $\beta$-code
and $\alpha \neq \beta$.
If $x \in \rfs{Pend}(F)$, then there are no $U \in \rfs{Nbr}^F(x)$
and $V \in \tau^H$ such that $U \cong V$.
\end{prop}

\noindent{\bf Proof }
Suppose by contradiction that $F$ is an $\alpha$-code,
$H$ is a $\beta$-code, $\alpha \neq \beta$ and there are
$x \in \rfs{Pend}(F)$, $U' \in \rfs{Nbr}^F(x)$,
$V' \in \tau^H$ and $\varphi$
such that $\iso{\varphi}{U'}{V'}$.
Note that by the definition of codes, $\alpha > 1$.
Let $U$ be an open subset  of $U'$
such that $F_1 \eqdf \rfs{cl}(U) \subseteq U'$,
and set $V = \varphi(U)$ and $H_1 = \varphi(F_1)$.
Note that $\rfs{Pend}(F) = \rfs{PD}_{\alpha}(F)$
and hence $\rfs{PD}_{\alpha}(F)$ is a nonempty perfect set.
By Proposition~\ref{p2.19},
$\rfs{PD}_{\alpha}(F_1) \cap U = \rfs{PD}_{\alpha}(F) \cap U$.
So
$\rfs{PD}_{\alpha}(F_1) \cap U \neq \emptyset$
and $\rfs{PD}_{\alpha}(F_1) \cap U$ has no isloated points.
Now, since $\varphi \nrestriction F_1$
is a homeomorphism between $F_1$ and $H_1$ which takes $U$ to $V$,
we have that
$\rfs{PD}_{\alpha}(H_1) \cap V = 
\varphi(\rfs{PD}_{\alpha}(F_1) \cap U)$.
It follows that
\begin{itemize}
\addtolength{\parskip}{-11pt}
\addtolength{\itemsep}{06pt}
\item[$(\dagger)$] 
$\rfs{Po}(\rfs{PD}_{\alpha}(H_1) \cap V) \neq \emptyset$.
\vspace{-2mm}
\end{itemize}
By Proposition~\ref{p2.19},
$\rfs{PD}_{\alpha}(H_1) \cap V = \rfs{PD}_{\alpha}(H) \cap V$.
The only ordinal $\gamma > 1$ for which
$\rfs{Po}(\rfs{PD}_{\gamma}(H))$ is nonempty
is $\beta$ but $\alpha > 1$, and is different from $\beta$.
So $\rfs{Po}(\rfs{PD}_{\alpha}(H_1) \cap V) = \emptyset$.
This contradicts $(\dagger)$, hence the Proposition is proved.
\smallskip\hfill\proofend

\begin{defn}\label{d2.22}
\begin{rm}

Let $\itOmega$ be a limit ordinal and $F$ be a compact Hausdorff space.
We say that $F$ is {\it $\itOmega$-demonstrative}
if\smallskip
\begin{list}{}
{\setlength{\leftmargin}{33pt}
\setlength{\labelsep}{08pt}
\setlength{\labelwidth}{20pt}
\setlength{\itemindent}{-00pt}
\addtolength{\topsep}{-11pt}
\addtolength{\parskip}{-12pt}
\addtolength{\itemsep}{-05pt}
}
\item[(D1)]
$\itOmega(F) = \itOmega$ and $\itOmega$ is not attained in $F$,
\item[(D2)]
$\rfs{Pend}(F) \subseteq \rfs{Good}(F)$.
\end{list}
\index{D@demonstrative set: $\itOmega$-demonstrative set}
\hfill\myqed
\end{rm}
\end{defn}

In the next lemma we use the following properties of members of
$K_{\srfs{CII}}$.

\kern2mm

\begin{list}{}
{\setlength{\leftmargin}{41pt}
\setlength{\labelsep}{08pt}
\setlength{\labelwidth}{32pt}
\setlength{\itemindent}{00pt}
\addtolength{\topsep}{-11pt}
\addtolength{\parskip}{-12pt}
\addtolength{\itemsep}{-05pt}
}
\item[(TH1)] 
$X$ is tightly Hausdorff.
\item[(TH2)] 
For every relatively discrete subset $A$ of $X$,
if $\rfs{cl}(A)$ is not scattered, then there is $B \subseteq A$
such that $\rfs{acc}(B)$ is a nonempty perfect set.
\item[(TH3)] 
For every $A \subseteq X$,
if $\lambda \eqdf \abs{A}$ is an infinite regular cardinal,
then either $A$ has at least two $\lambda$-accumulation points,
or there is $B \subseteq A$ such that $\abs{B} = \lambda$,
$B$ is relatively discrete and $\rfs{cl}(B)$ is scattered.
\end{list}
Observe the following fact.

\begin{prop}\label{p2.20n}
For any of the properties (TH1)\,-\,(TH3),
if $\kern2ptY$ is a closed subspace of a space having the property,
then $Y$ has the same property.
\end{prop}

We shall later use two other properties of members of $K_{\srfs{CII}}$.
\begin{itemize}
\addtolength{\parskip}{-11pt}
\addtolength{\itemsep}{06pt}
\item
$X$ is sequentially compact.
\end{itemize}
This property is used in showing that
$\itOmega(X)$ is not attained -- a fact which is assumed
in the next lemma.
See Proposition~\ref{p2.8} and Corollary~\ref{c2.9}.
Another (and last) property to be used is
\begin{itemize}
\addtolength{\parskip}{-11pt}
\addtolength{\itemsep}{06pt}
\item
For every infinite cardinal $\lambda$ and a closed subset
$F \subseteq X$:
if $\abs{F} = 2^{\lambda^+}$,
then there is a scattered subspace $H \subseteq F$
such that $\abs{H} = \lambda^+$.
\vspace{-2mm}
\end{itemize}
This is proved in Proposition~\ref{p2.31},
and it is used at the end of the proof of Theorem~\ref{t1.2},
where it is shown that $\abs{\rfs{Good}(X)}$
cannot be much larger than $\abs{\itOmega(X)}$.

Let $K_{\srfs{TH}}$ be the class of all compact Hausdorff spaces
that have Properties~(TH1)\,-\,(TH3).
\index{N@$K_{\srfs{TH}}$. The class of all compact Hausdorff spaces
       that have\\\indent Properties~(TH1)\,-\,(TH3)}
Note that by Lemma~\ref{l2.3}, Proposition~\ref{p2.16}
and Proposition~\ref{p2.14}, $K_{\srfs{CII}} \subseteq K_{\srfs{TH}}$.

\begin{lemma}\label{l2.24}
Let $X \in K_{\srfs{TH}}$.
Suppose that $\rfs{ker}(X) \neq \emptyset$ and that $\itOmega(X)$
is not attained in $X$.
Let $g \in \rfs{Good}(X)$ and $V \in \rfs{Nbr}(g)$.
Then for every\break
$\alpha \in \itOmega(X) - \sngltn{0}$,
$V$ contains an $\itOmega(X)$-demonstrative $(\alpha + 1)$-code.
\end{lemma}

\noindent{\bf Proof }
Let $H$ be a closed neighborhood of $g$ such that $H \subseteq V$.
Clearly, we may replace $V$ by $H$.
Also, $H \in K_{\srfs{TH}}$ (this follows from~\ref{p2.20n}),
$\rfs{ker}(H) \neq \emptyset$,
$\itOmega(H) = \itOmega(X)$,
and thus $\itOmega(H)$ is not attained in $H$.
So we may replace $H$ by $X$ and prove that $X$ contains
an $\itOmega(X)$-demonstrative $(\alpha + 1)$-code.

Property (TH3) is just $(*)$ of Proposition~\ref{p2.14},
and (TH1) is stronger than being strongly Hausdorff.
So \ref{p2.14} implies that $\rfs{Good}(X)$ is a nonempty perfect set.

By~Proposition~\ref{p2.17},
there is a relatively discrete subset $A \subseteq \rfs{Good}(X)$
such that $\rfs{cl}(A)$ contains a nonempty perfect set,
and by (TH2) we may assume that
$\rfs{acc}(A)$ is perfect and nonempty.

By (TH1),
$A$ has a tight Hausdorff system
$\calU = \setm{U_a}{a \in A}$.
For every $a \in A$ let $F_a$ be a subset of $U_a$
such that $F_a$ is compact and scattered,
$a \not\in F_a$ and $\rfs{rk}(F_a) = \alpha$.
The existence of $F_a$ is assured by the goodness of $a$.
Also, choose a closed neighborhood $H'_a$ of $a$
such that $H'_a \subseteq U_a$ and $H'_a \cap F_a = \emptyset$,
and define $H_a = H'_a \cap \rfs{ker}(X)$.
Set $S = \bigcup_{a \in A} F_a$ and $T = \bigcup_{a \in A} H_a$
and define $C = \rfs{cl}(S \cup T)$.
We shall show that $C$ is an
$\itOmega(X)$-demonstrative $(\alpha + 1)$-code.

We start with the fact that $C$ is a code.
Clearly,
\begin{itemize}
\addtolength{\parskip}{-05pt}
\addtolength{\itemsep}{06pt}
\item[(1)] 
$F_a,H_a \in \rfs{Clop}(C)$ for every $a \in A$.
\end{itemize}
Since $\calU$ is tight and $F_a \cup H_a \subseteq U_a$
for every $a \in A$, it follows that
\begin{itemize}
\addtolength{\parskip}{-05pt}
\addtolength{\itemsep}{06pt}
\item[(2)] 
$\rfs{acc}(\setm{F_a \cup H_a}{a \in A}) =
\rfs{acc}(\setm{F_a}{a \in A}) =
\rfs{acc}(\setm{H_a}{a \in A}) = \rfs{acc}(A)$.
\end{itemize}
Also recall that
\begin{itemize}
\addtolength{\parskip}{-05pt}
\addtolength{\itemsep}{06pt}
\item[(3)] 
$\rfs{acc}(A)$ is perfect.
\end{itemize}
and
\begin{itemize}
\addtolength{\parskip}{-05pt}
\addtolength{\itemsep}{06pt}
\item[(4)] 
$F_a$ is scattered with rank $\alpha$ and $H_a$ is perfect.
\end{itemize}
From (2) and the tightness of $\calU$ it follows that
$$
\mbox{$
C = (\bigcup_{a \in A} F_a) \cup (\bigcup_{a \in A} H_a) \cup
\rfs{acc}(A),$}
$$
{\thickmuskip=2mu \medmuskip=1mu \thinmuskip=1mu 
and Facts (1)\,-\,(4) imply that
$\rfs{Po}(C) = \bigcup_{a \in A} H_a$
and that $\rfs{Is}(C) = \bigcup_{a \in A} \rfs{Is}(F_a)$.}
It thus follows that
$$
\mbox{$
\rfs{PD}(C) = (\bigcup_{a \in A} \rfs{D}(F_a)) \cup \rfs{acc}(A).
$}
$$
From the tightness of $\calU$ and Fact (4) it now follows that
\begin{itemize}
\addtolength{\parskip}{-11pt}
\addtolength{\itemsep}{06pt}
\item[(5)] 
for every $0 < \beta \leq \alpha$,
$\rfs{PD}_{\beta}(C) =
(\bigcup_{a \in A} \rfs{D}_{\beta}(F_a)) \cup \rfs{acc}(A)$
and\break
$\rfs{PD}_{\alpha + 1}(C) = \rfs{acc}(A)$.
\end{itemize}
Since $\rfs{acc}(A)$ is perfect, we conclude that
$\rfs{PD}_{\alpha + 2}(C) = \emptyset$
and hence\break
(i) $\rfs{prk}(C) = \alpha + 1$.
From the first part of (5) and from (4) and (1) it follows that
(ii) for every $0 <\beta < \alpha + 1$,
$\rfs{Po}(\rfs{PD}_{\beta}(C)) = \emptyset$.
Finally, from the second part of (5) and from (i) we have
(iii) $\rfs{Pend}(C) = \rfs{PD}_{\alpha + 1}(C) = \rfs{acc}(A)$,
and from (3) we conclude that (iv) $\rfs{Pend}(C)$ is perfect.
Facts (i), (ii) and (iv) are Clauses (C1)\,-\,(C3) in the definition of
an $(\alpha + 1)$-code.
We have thus proved that $C$ is an $(\alpha + 1)$-code.

We show that $C$ is $\itOmega(X)$-demonstrative.
Let $a \in A$, then $H_a$ is the intersection of a closed neighborhood
of $a$ with $\rfs{ker}(X)$. Since $a \in \rfs{Good}(X)$,
it follows that $\itOmega(H_a) = \itOmega(X)$,
and that $a \in \rfs{Good}(H_a)$.
Since $H_a \subseteq C \subseteq X$,
we also have that $\itOmega(H_a) \leq \itOmega(C) \leq \itOmega(X)$.
Hence $\itOmega(C) = \itOmega(X)$, and since $\itOmega(X)$
is not attained in $X$, it is not attained in $C$.
We have shown that
$$\hbox{
$\itOmega(C) = \itOmega(X)$ and $\itOmega(X)$ is not attained in $C$.}
$$
That is, $C$ fulfills
Clause (D1) in the definition of demonstrative sets.

We have also shown above that for every $a \in A$,
$a \in \rfs{Good}(C)$. That is, $A \subseteq \rfs{Good}(C)$.
In Fact (iii) we proved that $\rfs{Pend}(C) = \rfs{acc}(A)$.
So $\rfs{Pend}(C) \subseteq \rfs{acc}(\rfs{Good}(C))$.
Since $\rfs{Good}(C)$ is closed,
$\rfs{Pend}(C) \subseteq \rfs{Good}(C)$.
So Clause (D2) holds.
We have shown that $C$ is $\itOmega(X)$-demonstrative.
\hfill\proofend

\kern2mm

\noindent
{\bf 5. Coding subsets of $\bf \Omega(X)$ and proliferation systems.}
\smallskip

\noindent
The assumption ``$X$ is a non-scattered CO space'' is contradictory.
To reach this contradiction, we show that
$\abs{\rfs{Good}(X)}$ is much larger than $\abs{\itOmega(X)}$.
First we code subsets of $\itOmega(X)$ by subsets of $\rfs{Good}(X)$.
This coding implies that
$\abs{\rfs{Good}(X)} \geq 2^{\abs{\itOmega(X)}}$.
Next we code sets of subsets of $\itOmega(X)$ by subsets of
$\rfs{Good}(X)$.
This leads to the conclusion
that $\abs{\rfs{Good}(X)} \geq 2^{2^{\abs{\itOmega(X)}}}$.
We repeat this procedure twice more and then reach a contradiction.
The above three steps use an identical argument
which in the first case is applied to the set of $\alpha$-codes,
and in the second, to the set of codes of subsets of $\itOmega(X)$.
The notion which provides the unified argument is called
a ``proliferation system'',
and the conclusion of the iterated use of this argument
is stated in Corollary \ref{tightly-haus-c1.5}.
It will be evident that for an any limit ordinal $\itOmega$ 
the set of $\itOmega$-demonstrative codes is a proliferation system,
and that if $X$ is a non-scattered CO space, then all these codes are
realized in $X$.
This makes Corollary \ref{tightly-haus-c1.5} applicable to $X$.

A pair $\whatX = \pair{X}{e}$ consisting of a topological space $X$
and a point $e \in X$
is called a {\it pointed space}.
\index{D@pointed space. A pair $\pair{X}{x}$,
       where $X$ is a topological space and $x \in X$}
Let $\calP = \setm{\whatcalX_t}{t \in T}$
be an indexed family such that for every $t \in T$,
$\whatcalX_t$ is a class of pointed spaces closed under homeomorphisms.
Then $\calP$ is called a {\it type system}.
By ``closed under homeomorphisms'' we mean that if
$\pair{X}{e} \in \whatcalX_t$ and $\varphi$ is a homeomorphism
between $X$ and $Y$,
then $\pair{Y}{\varphi(e)} \in \whatcalX_t$.
Denote $T$ by $T_{\calP}$
and set $\whatcalX_{\calP} = \bigcup_{t \in T} \whatcalX_t$.
For $t \in T$ define
$\calX_t \eqdf \setm{X}{\mbox{there is } e \in X \mbox{ such that }
\pair{X}{e} \in \whatcalX_t}$
and $\calX_{\calP} \eqdf \bigcup_{t \in T} \calX_t$.
\index{D@type system}
\index{N@$\whatcalX_t$. The family of pointed spaces associated with $t$
       in a P-system $\calP$}
\index{N@$T_{\calP}$. The index set of a type system $\calP$}
\index{N@$\calX_t$. The family of spaces associated with $t$
       in a P-system $\calP$}
\index{N@$\whatcalX_{\calP}$. The class of all pointed spaces of
       a type system $\calP$}
\index{N@$\calX_{\calP}$. The class of all spaces of
       a type system $\calP$}

Let $\calX$ be a class of topological spaces
and $Y$ be a topological space.
We say that $\calX$ {\it occurs} in $Y$ if there is $X \subseteq Y$
such that $X \in \calX$.
We say that a class $\whatcalX$ of pointed spaces {\it occurs} in $Y$
if there are $X \subseteq Y$ and $e \in X$ such that
$\pair{X}{e} \in \whatcalX$.
\index{D@occurs: $\calX$ occurs in $Y$}

\begin{defn}\label{tightly-haus-d1.1}
\begin{rm}
(a)
A type system $\calP = \setm{\whatcalX_t}{t \in T}$
is called a {\it proliferation system
(P-system)} if the following hold.
\begin{list}{}
{\setlength{\leftmargin}{40pt}
\setlength{\labelsep}{08pt}
\setlength{\labelwidth}{32pt}
\setlength{\itemindent}{-00pt}
\addtolength{\topsep}{06pt}
\addtolength{\parskip}{-12pt}
\addtolength{\itemsep}{-05pt}
}
\item[(P1)] 
$T$ is infinite and every member of $\calX_{\calP}$
is compact Hausdorff.
\item[(P2)] 
Suppose that $s,t \in T$ are distinct,
$\pair{F}{d} \in \whatcalX_s$, $H \in \calX_t$
and\break
$V \in \rfs{Nbr}^F(d)$.
Then there is no $U \in \tau^H$ such that $V \cong U$.
\item[(P3)] 
For every $\pair{F}{d} \in \whatcalX_{\calP}$, $V \in \rfs{Nbr}^F(d)$
and $t \in T$
there is $Y \in \calX_t$ such that $Y \subseteq V$.
\end{list}
\index{D@proliferation system}
\index{D@p-system@@P-system. Abbreviation of a proliferation system}
Note the definition of a proliferation system does not exclude the
possibility that for every $t \in T_{\calP}$, $\whatcalX_t = \emptyset$.
However, by (P3),
if for some $t \in T_{\calP}$, $\whatcalX_t \neq \emptyset$,
then for all $t \in T_{\calP}$, $\whatcalX_t \neq \emptyset$.

(b) Let $\calP$ be a P-system and $X$ be a compact Hausdorff space.
We say that $\calP$ {\it occurs} in $X$ if $\calX_{\calP}$
occurs in $X$.
\index{D@occurs: A P-system occurs in $Y$}

(c)
Let $\calP = \setm{\whatcalX_t}{t \in T}$ be a P-system
and $\emptyset \neq \itGamma \subseteq T_{\calP}$.
A pointed compact Hausdorff space $\pair{F}{d}$
is called a {\it $\itGamma$-marker}
if there is a family $\calF$ of subsets of $F$
and $\setm{d_{F'}}{F' \in \calF}$ such that
\begin{list}{}
{\setlength{\leftmargin}{40pt}
\setlength{\labelsep}{08pt}
\setlength{\labelwidth}{32pt}
\setlength{\itemindent}{-00pt}
\addtolength{\topsep}{06pt}
\addtolength{\parskip}{-12pt}
\addtolength{\itemsep}{-05pt}
}
\item[(M1)] 
$\calF \subseteq \rfs{Clop}(F)$,
\item[(M2)] 
$\calF$ is a tight family,
\item[(M3)] 
for every $F' \in \calF$,
$\pair{F'}{d_{F'}} \in \bigcup_{t \in \itGamma} \whatcalX_t$,
\item[(M4)] 
for every $V \in \rfs{Nbr}(d)$ and $t \in \itGamma$
there is $F' \in \calF \cap \calX_t$ such that $d_{F'} \in V$,
\item[(M5)] $F = \rfs{cl}(\bigcup \calF)$.
\vspace{-2.5mm}
\end{list}
$\calF$ is called a {\it filler} for $\pair{F}{d}$.
\index{D@marker: $\itGamma$-marker}
\index{D@filler}

(c)
Denote the powerset of a set $A$ by $\bfP(A)$.
Suppose that $\calP = \setm{\whatcalX_t}{t \in T}$ is a P-system.
For every $\itGamma \in \bfP(T) - \sngltn{\emptyset}$ define
$$
\whatcalM_{\itGamma}^{\calP} = \setm{\pair{X}{d}}{\pair{X}{d}
\mbox{ is a $\itGamma$-marker}}.
$$
For ${\bfGamma} \subseteq \bfP(T) - \sngltn{\emptyset}$ define
$$
\calQ_{\bfGamma}^{\calP} =
\setm{\whatcalM_{\itGamma}^{\calP}}{\itGamma \in \mathbf \Gamma}.
$$
\index{N@$\bfP(A)$. Powerset of $A$}
\index{N@$\calQ_{\bfGamma}^{\calP}$}
\index{N@$\whatcalM_{\itGamma}^{\calP} = \setm{\pair{X}{d}}{\pair{X}{d}
       \mbox{ is a $\itGamma$-marker}}$}
\index{N@$\calQ_{\bfGamma}^{\calP} =
       \setm{\whatcalM_{\itGamma}^{\calP}}
       {\itGamma \in \mathbf \Gamma}$}
\rule{0pt}{0pt}\hfill\myqed
\end{rm}
\end{defn}

In general $\calQ_{\bfGamma}^{\calP}$, need not be a P-system,
but we shall see that if $\bfGamma$ is an infinite set
of pairwise incomparable subsets of $T$ with the same cardinality,
then $\calQ_{\bfGamma}^{\calP}$ is a P-system.

\begin{prop}\label{tightly-haus-p1.2}
\num{a} Let $\calP$ be a P-system and
$\itGamma,\itDelta \in \bfP(T_{\calP}) - \sngltn{\emptyset}$,
and assume that $\itGamma \not\subseteq \itDelta$.
Suppose that $\pair{F}{d}$ and $\pair{H}{e}$
are respectively a $\itGamma$-marker and a $\itDelta$-marker.
Then there do not exist $U \in \rfs{Nbr}^F(d)$ and $V \in \tau^H$
such that $U \cong V$.

\num{b}
Let $\calP$ be a P-system
and $\bfGamma \subseteq \bfP(T_{\calP}) - \sngltn{\emptyset}$
be infinite.
Suppose that for every distinct $\itGamma,\itDelta \in \bfGamma$,
$\abs{\itGamma} = \abs{\itDelta}$
and
$\itGamma \not\subseteq \itDelta$.
Then $\calQ_{\bfGamma}^{\calP}$ is a P-system.
\end{prop}

\noindent{\bf Proof }
(a)
Let $\itGamma$, $\itDelta$, $\pair{F}{d}$ and $\pair{H}{e}$
be as specified in (a).
Suppose by way of contradiction that there are $U \in \rfs{Nbr}^F(d)$
and $V \in \tau^H$ such that $U \cong V$.
We may assume that $U = V$.
So for every $A \subseteq U$,
$\tau^F \nrestriction A = \tau^U \nrestriction A =
\tau^H \nrestriction A$,
and $A$ is open in $F$ iff $A$ is open in $U$ iff $A$ is open in $H$.

Let $\calF$ and $\calH$ be fillers for $\pair{F}{d}$ and $\pair{H}{e}$
respectively.
Let $t \in \itGamma - \itDelta$.
Let $F' \in \calX_t \cap \calF$
be such that $d_{F'} \in U$ and $F_0 = F' \cap U$.
So $d_{F'} \in F_0$,
and since $F'$ is open in $F$, it follows that $F_0$ is open in $U$,
and hence in $H$. That is, $F_0 \in \rfs{Nbr}^H(d_{F'})$.
Suppose by contradiction that $d_{F'} \in \bigcup \calH$.
Let $H' \in \calH$ be such that $d_{F'} \in H'$.
There is $s \in \itDelta$ such that $H' \in \calX_s$.
Then $s \neq t$. Since $H'$ is open in $H$, we have that
$H' \cap F_0$ is open in $H$, and hence it is open in $F$.
This implies that $H' \cap F_0$ is open in $F'$.
Hence
$$H' \cap F_0 \in \rfs{Nbr}^{F'}(d_{F'}).$$
On the other hand,
$$
\hbox{$H' \cap F_0$ is open in $H'$.}
$$
Recall that $\pair{F'}{d_{F'}} \in \whatcalX_t$
and that $H' \in \calX_s$.
The last four mentioned facts contradict Property (P2) of~$\calP$.

It follows that $d_{F'} \in H - \bigcup \calH$.
Since every member of $\calH$ is clopen in $H$,
$d_{F'} \in H - \bigcup \setm{\rfs{cl}^H(H')}{H' \in \calH}$.
Now, $\rfs{cl}^H(\bigcup \calH) = H$,
hence $d_{F'} \in \rfs{acc}^H(\calH)$.
From the tightness of $\calH$ it follows that
$d_{F'} \in \rfs{acc}^H(\setm{d_{H'}}{H' \in \calH}$.
Recall that $F_0 \in \rfs{Nbr}^H(d_{F'})$.
So there is $H' \in \calH$ such that $d_{H'} \in F_0$.
Clearly, $F_0 \cap H'$ is open in $H'$ and hence
$$
H' \cap F_0 \in \rfs{Nbr}^{H'}(d_{H'}).
$$
$H' \cap F_0$ is open in $H$ and it is a subset of $U$.
So it is open in $F$. It follows that
$$
\hbox{$H' \cap F_0$ is open in $F'$.}
$$
There is $s \neq t$, such that $\pair{H'}{d_{H'}} \in \whatcalX_s$,
and on the other hand, $F' \in \calX_t$.
These facts contradict Property (P2) of $\calP$.
This proves (a).

(b)
Denote $\calQ^{\calP}_{\bfGamma}$ by $\calQ$.
So $T_{\calQ} = \bfGamma$ and for every $\itGamma \in \bfGamma$,
$\whatcalX_{\itGamma} = \whatcalM^{\calP}_{\itGamma}$.
By definition, every $\itGamma$-marker is compact Hausdorff.
So every member of $\calX_{\calQ}$ is compact Hausdorff. Since also,
$\bfGamma$ is infinite, $\calQ$ fulfills (P1).

That
$\calQ$ fulfills~(P2), was indeed proved in Part (a).

We prove (P3). Let $\pair{F}{d} \in \whatcalX_{\calQ}$
and $\itDelta \in \bfGamma$.
There is $\itGamma \in \bfGamma$
such that $\pair{F}{d} \in \whatcalX_{\itGamma}$.
Suppose that $\calF$ is a filler for $\pair{F}{d}$.
Let $\fnn{f}{\itGamma}{\itDelta}$ be a bijection.
If $F' \in \calF \cap \calX_t$, choose $H_{F'} \subseteq F'$ such that
$H_{F'} \in \calX_{f(t)}$
and define $\calH = \setm{H_{F'}}{F' \in \calF}$
and $H = \rfs{cl}(\bigcup \calH)$.
That $H_{F'}$ exists follows from (P3) applied to $\calP$.
The tightness of $\calF$ implies that $d \in \rfs{cl}(\bigcup \calH)$.
It now follows trivially from the definition of $H$ and $\calH$ that
$\pair{H}{d} \in \whatcalX_{\itDelta}$.
So $\calQ$ fulfills (P3).
\smallskip\hfill\proofend

Note that if $\calP$ is a P-system and $\calP$ occurs in $X$,
then for every $t \in T_{\calP}$, $\calX_t$ occurs in $X$.
This follows from (P3).
Suppose that $X$ is a tightly Hausdorff compact CO space
and $\calP$ is a P-system occurring in $X$.
We shall show that $\calM_{\itGamma}^{\calP}$ occurs in $X$
for every $\itGamma \in \bfP(T_{\calP}) - \sngltn{\emptyset}$.
In order to show this, we first establish the existence of
so-called $\mu$-special $\sngltn{t}$-markers.

Let $X$ be a topological space, $A \subseteq X$, $x \in X$
and $\mu$ be an infinite cardinal.
Denote the set of $\mu$-accumulation points of $A$ in $X$ by
$\rfs{acc}^X_{\mu}(A)$. We use $\rfs{acc}_{\mu}(A)$ as
an abbreviation of the above.
\index{N@$\rfs{acc}_{\mu}(A)$.
       The set of $\mu$-accumulation points of $A$}
Let $\calP$ be a P-system, $t \in T_{\calP}$
and $\pair{F}{e}$ be a $\sngltn{t}$-marker with a filler $\calF$.
For every $F' \in \calF$ choose $e_{F'}$
such that $\pair{F'}{e_{F'}} \in \whatcalX_t$.
We call $\pair{F}{e}$ a {\it $\mu$-special $\sngltn{t}$-marker}
if\break
$e \in \rfs{acc}_{\mu}(\setm{e_{F'}}{F' \in \calF})$.
\index{D@marker: $\mu$-special $\sngltn{t}$-marker}

\begin{prop}\label{tightly-haus-p1.3}
Let $X$ be a compact Hausdorff CO space and $\calP$ be a\break
P-system such that $\calP$ occurs in $X$.

\num{a} There is a set $\setm{\pair{G_t}{g_t}}{t \in T_{\calP}}$
such that
\begin{itemize}
\addtolength{\parskip}{-06pt}
\addtolength{\itemsep}{01pt}
\item[\num{1}] 
for every distinct $s,t \in T_{\calP}$, $g_s \neq g_t$,
\item[\num{2}] 
for every $t \in T_{\calP}$,
$\pair{G_t}{g_t} \in \whatcalX_t$ and
$G_t \subseteq X$,
\item[\num{3}] 
$\setm{g_t}{t \in T_{\calP}}$ is relatively discrete.
\vspace{-2mm}
\end{itemize}

\num{b} Suppose that in addition to the above,
$X$ is tightly Hausdorff.
Then for every $t \in T_{\calP}$
there is $K \subseteq X$ and $c \in K$
such that $\pair{K}{c}$
is a $\abs{T_{\calP}}$-special $\sngltn{t}$-marker.
\end{prop}

\noindent{\bf Proof }
Suppose that $\calP = \setm{\whatcalX_t}{t \in T}$.

(a)
Let $\pair{F}{d} \in \whatcalX_{\calP}$
be such that
$F \subseteq X$.
By (P3), for every $t \in T$ there is $\pair{H_t}{e_t} \in \whatcalX_t$
such that $H_t \subseteq F$.
Let $G_t$ be a clopen subset of $X$ homeomorphic to $H_t$
and let $g_t \in G_t$ be the image of $e_t$
under the homeomorphism between $H_t$ and $G_t$.
Suppose by way of contradiction that for some distinct $s,t \in T$,
$g_s \in G_t$. Then $G_s \cap G_t \in \rfs{Nbr}^{G_s}(g_s)$
and $G_s \cap G_t$ is open in $G_t$.
This contradicts (P2), so for every distinct $s,t \in T$,
$g_s \not\in G_t$.
We thus have that for every $t \in T$, $G_t$ is open and
$G_t \cap \setm{g_s}{s \in T} = \sngltn{g_t}$.
This means that $\setm{g_t}{t \in T}$ is relatively discrete
and that for every distinct $s,t \in T$, $g_s \neq g_t$.
So $\setm{\pair{G_t}{g_t}}{t \in T}$ is as required.

(b) Denote $\mu = \abs{T}$
and let $A = \setm{g_t}{t \in T}$ be as assured in (a).
So $A$ is relatively discrete, $\abs{A} = \mu$
and for every $a \in A$ there is $G \subseteq X$
such that $\pair{G}{a} \in \whatcalX_{\calP}$.
It is trivial that in a compact space every set of cardinality
$\mu$ has a $\mu$-accumulation point.
So let $c$ be a $\mu$-accumulation point of $A$.

Fix $t \in T$. We construct a set $K \subseteq X$ such that
$\pair{K}{c}$ is a
$\sngltn{t}$-marker.
Let $\setm{U_a}{a \in A}$ be a tight Hausdorff system for $A$.
If $a \in A$, then for some $s \in T$, $a = g_s$.
Since $g_s \in G_s$ and $\pair{G_s}{g_s} \in \whatcalX_s$,
we may apply (P3) to $\pair{G_s}{g_s}$.
Now, $U_a \in \rfs{Nbr}(g_s)$,
so there is $\pair{E_a}{e_a} \in \whatcalX_t$
such that $E_a \subseteq U_a \cap G_s$.
Let $\calE = \setm{E_a}{a \in A}$.
Define $K = \rfs{cl}(\bigcup \calE)$.
In the definition of a $\itGamma$-marker
we need to have a choice function $\setm{d_{F'}}{F' \in \calF}$.
So define $d_{E_a}$ to be $e_a$ for every $a \in A$.
It follows trivially from the construction that (M1)\,-\,(M5) hold for
$\pair{K}{c}$, $\calE$, $\setm{e_a}{a \in A}$ and $\sngltn{t}$.
That is, $\pair{K}{c}$ is a $\sngltn{t}$-marker.
Recall that $c$ was chosen to be a $\mu$-accumulation point of $A$.
Since $\setm{U_a}{a \in A}$ is a tight family, and
$a,e_a \in U_a$ for every $a \in A$, it follows that
$c$ is a $\mu$-accumulation point of $\setm{e_a}{a \in A}$.
This assures that $\pair{K}{c}$ is $\mu$-special.
\rule{10pt}{0pt}\smallskip\hfill\proofend
 
\begin{lemma}\label{tightly-haus-l1.4}
Let $X$ be a tightly Hausdorff compact CO space and $\calP$ be a
P-system such that $\calP$ occurs in $X$.
Then for every
$\itGamma \in \bfP(T_{\calP}) - \sngltn{\emptyset}$,
there are $H \subseteq X$ and $c \in H$
such that $\pair{H}{c}$ is a $\itGamma$-marker.
\end{lemma}

\noindent{\bf Proof }
Let $\calP = \setm{\calX_t}{t \in T}$ and denote $\abs{T}$ by $\mu$.
Choose a countable subset $T_0$ of $T$,
and for every $t \in T_0$
let $\pair{K_t}{c_t}$ be a $\mu$-special $\sngltn{t}$-marker
such that $K_t$ is a clopen subset of $X$.
The existence of a $\mu$-special $\pair{K_t}{c_t}$
was proved in Proposition~\ref{tightly-haus-p1.3}(b),
and that $K_t$ may be a clopen set follows from the fact that
$X$ is a CO space.
 
Let $t \in T_0$. Since $\pair{K_t}{c_t}$ is $\mu$-special,
there are a filler $\calF_t$ for $\pair{K_t}{c_t}$
and a set $\setm{e_{F'}}{F' \in \calF_t}$ such that
\begin{itemize}
\addtolength{\parskip}{-11pt}
\addtolength{\itemsep}{06pt}
\item[(1)] 
for every $F' \in \calF_t$, $\pair{F'}{e_{F'}} \in \whatcalX_t$,
\item[(2)] 
$c_t \in \rfs{acc}_{\mu}(\setm{e_{F'}}{F' \in \calF_t})$.
\vspace{-2mm}
\end{itemize}

By Proposition \ref{tightly-haus-p1.2}(a),
for every distinct $s,t \in T_0$,
$c_s \not\in K_t$.
This implies that $\setm{c_t}{t \in T_0}$
is infinite and relatively discrete.
For every $t \in T_0$ choose $V_t \in \rfs{Nbr}(c_t)$
in such a way that
\begin{itemize}
\addtolength{\parskip}{-11pt}
\addtolength{\itemsep}{06pt}
\item[(3)] 
$\setm{V_t}{t \in T_0}$ is a tight Hausdorff system
for $\setm{c_t}{t \in T_0}$.
\vspace{-2mm}
\end{itemize}
Define $\calF'_t = \setm{F' \in \calF_t}{e_{F'} \in V_t}$
and $E_t = \setm{e_{F'}}{F' \in \calF'_t}$.
Also choose $c \in \rfs{acc}(\setm{c_t}{t \in T_0})$. So
\begin{itemize}
\addtolength{\parskip}{-06pt}
\addtolength{\itemsep}{01pt}
\item[(4)] 
$c \in \rfs{acc}(\setm{V_t}{t \in T_0})$.
\end{itemize}

Let $\itGamma \in \bfP(T_{\calP}) - \sngltn{\emptyset}$.
We construct $H$ such that $\pair{H}{c}$ is a $\itGamma$-marker.
For every $t \in T_0$ let $\fnn{f_t}{\calF'_t}{\itGamma}$
be a surjection.
Let $t \in T_0$ and $F' \in \calF'_t$.
Then $V_t \in \rfs{Nbr}(e_{F'})$.
We use (P3) and the fact that
$\pair{F'}{e_{F'}} \in \whatcalX_{\calP}$
in order to conclude that there are $H_{F'}$ and $d_{F'}$ such that
\begin{itemize}
\addtolength{\parskip}{-11pt}
\addtolength{\itemsep}{06pt}
\item[(5)] 
$\pair{H_{F'}}{d_{F'}} \in \whatcalX_{f_t(F')}$
and $H_{F'} \subseteq V_t \cap F'$.
\vspace{-2mm}
\end{itemize}

Let $\calH = \setm{H_{F'}}{t \in T_0 \mbox{ and } F' \in \calF'_t}$.
For $H' = H_{F'} \in \calH$ denote $d_{F'}$ by~$b_{H'}$.
Define $H = \rfs{cl}(\bigcup \calH)$.
We verify that $\pair{H}{c}$ is a $\itGamma$-marker,
that $\calH$ is a filler for $\pair{H}{c}$
and that $\setm{b_{H'}}{H' \in \calH}$ is the choice function required
in the definition of a $\itGamma$-marker.

That (M1), (M3) and (M5) hold is trivial.
We check that (M2) holds.
We have to show that $\calH$ is a tight family.
For $t \in T_0$ set $\calH_t = \setm{H_{F'}}{F' \in \calF'_t}$.
Then
\begin{itemize}
\addtolength{\parskip}{-11pt}
\addtolength{\itemsep}{06pt}
\item[(i)] 
$\calH = \bigcup_{t \in T_0} \calH_t$.
\item[(ii)] 
For every $t \in T_0$, $\calH_t$ is a tight family.
\vspace{-2mm}
\end{itemize}
That $\calH_t$ is tight follows from the facts:
$\calF'_t$ is tight and $H_{F'} \subseteq F'$
for every $F' \in \calF'_t$. By (5),
\begin{itemize}
\addtolength{\parskip}{-11pt}
\addtolength{\itemsep}{06pt}
\item[(iii)] 
$H' \subseteq V_t$ for every $t \in T_0$ and $H' \in \calH_t$,
\item[(iv)] 
$\setm{V_t}{t \in T_0}$ is tight.
\vspace{-2mm}
\end{itemize}
Facts (i)\,-\,(iv) easily imply that $\calH$ is tight.
So $\calH$ satisfies (M2).

We next verify (M4).
Let $s \in \itGamma$ and $W \in \rfs{Nbr}(c)$.
For every $t \in T_0$ choose $F'_t \in \calF'_t$
such that $f_t(F'_t) = s$.
Set $H^0_t = H_{F'_t}$ and $a_t = b_{H^0_t}$.
So $a_t \in V_t$ and $\pair{H^0_t}{a_t} \in \whatcalX_s$.
From the facts:
\begin{itemize}
\addtolength{\parskip}{-11pt}
\addtolength{\itemsep}{06pt}
\item 
$c$ is an accumulation point of $\setm{V_t}{t \in T_0}$,
\item 
$a_t \in V_t$ for every $t \in T_0$,
\item 
$\setm{V_t}{t \in T_0}$ is tight,
\vspace{-2mm}
\end{itemize}
we conclude that $c \in \rfs{acc}(\setm{a_t}{t \in T_0})$.
So there is $t_0 \in T_0$ such that $a_{t_0} \in W$.
Recall that $\pair{H^0_{t_0}}{a_{t_0}} \in \whatcalX_s$.
Also $\pair{H^0_{t_0}}{a_{t_0}} \in \calH$.
We have thus found $H' \in \calH$
such that $\pair{H'}{b_{H'}} \in \whatcalX_s$
and $b_{H'} \in W$.
This shows that
$H$, $\calH$ and $\setm{b_{H'}}{H' \in \calH}$ fulfill (M4).
\smallskip\hfill\proofend

Let $X$ be a compact Hausdorff space and $\calP$ be a P-system.
Define
$$
\rfs{End}_{\calP}(X)
=
\setm{e \in X}{\mbox{There is } F \subseteq X \mbox{ such that }
\pair{F}{e} \in \whatcalX_{\calP}},\ 
$$
and 
$\rfs{Good}_{\calP}(X) = \rfs{cl}(\rfs{End}_{\calP}(X))$.
For an infinite cardinal $\mu$ set
$\beths_0(\mu) = \mu$, for every $n \in \omega$,
$\beths_{n + 1}(\mu) = 2^{\sbeths_n(\mu)}$
and $\beths_{\omega}(\mu) = \bigcup_{n \in \omega} \beths_n(\mu)$.
\index{N@$\rfs{End}_{\calP}(X)$}
\index{N@$\rfs{Good}_{\calP}(X)$}

\begin{cor}\label{tightly-haus-c1.5}
Let $X$ be a tightly Hausdorff compact CO space,
$\calP$ be
a P-system and $\mu = \abs{T_{\calP}}$.
If $\calP$ occurs in $X$,
then $\abs{\rfs{Good}_{\calP}(X)} \geq \beths_4(\mu)$.
\end{cor}

\noindent{\bf Proof }
We prove that $\abs{\rfs{Good}_{\calP}(X)} \geq \beths_{\omega}(\mu)$.
Denote $\calP$ and $T_{\calP}$ by $\calP_0$ and $T_0$ respectively.
We define by induction a P-system $\calP_n$.
Suppose that $\calP_n$ has been defined.
For simplicity, denote $\calP_n$ by $\calR$
and $T_{\calP_n}$ by $T$. We assume by induction that
$\calR$ occurs in $X$, that $\abs{T} = \beths_n(\mu)$
and that $\rfs{Good}_{\calP_n}(X) \subseteq \rfs{Good}_{\calP}(X)$.
Let $\bfGamma \subseteq \bfP(T)$ be such that
\begin{itemize}
\addtolength{\parskip}{-06pt}
\addtolength{\itemsep}{01pt}
\item[(1)] 
$\abs{\bfGamma} = 2^{\abs{T}}$,
\item[(2)] 
for any distinct $\itGamma,\itDelta \in \bfGamma$,\,
$\abs{\itGamma} = \abs{T}$ and $\itGamma \not\subseteq \itDelta$.
\vspace{-2mm}
\end{itemize}
Then by Proposition \ref{tightly-haus-p1.2}(b),
$\calP_{n + 1} \eqdf \calQ_{\bfGamma_n}^{\calP_n}$ is a P-system.
Denote $\calP_{n + 1}$ by $\calS$.
By the induction Hypothesis, $\calR$ occurs in $X$, and hence
by Lemma~\ref{tightly-haus-l1.4}, $\calS$ occurs in $X$.
Note that $T_{\calS} = \bfGamma$.
So $\abs{T_{\calS}} = 2^{\abs{T}}$.
By Proposition~\ref{tightly-haus-p1.3}(a),
$\abs{\rfs{End}_{\calS}(X)} \geq \abs{T_{\calS}} = 2^{\abs{T}} =
\beths_{n + 1}(\mu)$ and since
$\rfs{Good}_{\calS}(X) \supseteq \rfs{End}_{\calS}(X)$,
$$
\abs{\rfs{Good}_{\calS}(X)} \geq \beths_{n + 1}(\mu).
$$
By the definition of markers,
$\rfs{End}_{\calS}(X) \subseteq \rfs{Good}_{\calR}(X)$.
Since $\rfs{Good}_{\calS}(X) =$
{\thickmuskip=2mu \medmuskip=1mu \thinmuskip=1mu 
$\rfs{cl}(\rfs{End}_{\calS}(X))$
and $\rfs{Good}_{\calR}(X)$ is closed, it follows that
$\rfs{Good}_{\calS}(X) \subseteq \rfs{Good}_{\calR}(X)$.}
Hence by the induction hypothesis,
$$
\rfs{Good}_{\calS}(X) \subseteq \rfs{Good}_{\calP}(X).
$$
This concludes the inductive construction.

Since for every $n$,
$\rfs{Good}_{\calP_n}(X) \subseteq \rfs{Good}_{\calP}(X)$
and $\abs{\rfs{Good}_{\calP_n}(X)} \geq \beths_n(\mu)$,
it follows that 
$\abs{\rfs{Good}_{\calP}(X)} \geq \beths_{\omega}(\mu)$.
\medskip\hfill\proofend

\kern2mm

\noindent
We shall apply Corollary~\ref{tightly-haus-c1.5}
to the class of all $\itOmega$-demonstrative codes.
To this end we show that this class forms a P-system.

For a limit ordinal $\itOmega$ and $\alpha < \itOmega$
define
\vspace{1.5mm}
\newline
\rule{7pt}{0pt}
\renewcommand{\arraystretch}{1.5}
\addtolength{\arraycolsep}{-0pt}
$
\begin{array}{lll}
\whatcalX^{\itOmega}_{\alpha}
&
=
&
\setm{\pair{F}{e}}
{F \in K_{\srfs{TH}},
\mbox{ $F$ is an $\itOmega$-demonstrative $(\alpha + 1)$-code}
\\
&
&
\mbox{\rule{3pt}{0pt}and } e \in \rfs{Pend}(F)}.
\vspace{1.7mm}
\end{array}
$
\renewcommand{\arraystretch}{1.0}
\addtolength{\arraycolsep}{0pt}
\newline
Now define
$\calP^{\itOmega} =
\setm{\whatcalX^{\itOmega}_{\alpha}}{1 \leq \alpha <\itOmega}$.
\index{N@$\whatcalX^{\itOmega}_{\alpha}$.
       The class of all pointed spaces which are
       $\itOmega$-domnstrative\\\indent
       {\thickmuskip=2mu \medmuskip=1mu \thinmuskip=1mu 
       $(\alpha + 1)$-codes
       with a member of $\rfs{Pend}(X)$ as their distiguished point}}
\index{N@$\calP^{\itOmega}$. The P-system of
       $\itOmega$-demonstrative codes}

\begin{prop}\label{p2.32-a}
If $\itOmega$ is a limit ordinal,
then $\calP^{\itOmega}$ is a proliferation system.
\end{prop}

\noindent{\bf Proof }
By the definitions, $\calP^{\itOmega}$ fulfills (P1).
That (P2) holds is proved in Proposition~\ref{p2.21}.

We prove (P3).
Let $\pair{F}{e} \in \whatcalX^{\itOmega}_{\alpha}$.
Then $F$ is an $\itOmega$-demonstrative $(\alpha + 1)$-code
and $e \in \rfs{Pend}(F)$.
Let $V \in \rfs{Nbr}^F(e)$ and $\beta < \itOmega$.
By definition, $F \in K_{\srfs{TH}}$, $\rfs{ker}(F) \neq \emptyset$,
$\itOmega(F) = \itOmega$ and $\itOmega(F)$ is not attained.
By Clause~(D2) of~\ref{d2.22}, $e \in \rfs{Good}(F)$.
So the assumptions of Lemma~\ref{l2.24} are fulfilled by
$F$, $e$ and $V$.
So by~\ref{l2.24},
$V$ contains an $\itOmega$-demonstrative $(\beta + 1)$-code.
\smallskip\hfill\proofend

\kern3mm

\noindent
{\bf 6. The conclusion of the proof.}
\smallskip

We next see that if $\itOmega(X)$ is not attained,
then $\calP^{\itOmega(X)}$ occurs in $X$.

\begin{prop}\label{p2.32-b}
Let $X \in K_{\srfs{TH}}$
be such that $\rfs{ker}(X) \neq \emptyset$ and $\itOmega(X)$
is not attained.
Then
$\calP^{\itOmega(X)}$ occurs in $X$
and
$\rfs{Good}_{\calP^{\raise1pt\hbox{\tiny$\itOmega(X)$}}}(X)
\subseteq
\rfs{Good}(X)$.
\end{prop}

\noindent{\bf Proof }
The fact ``$\calP^{\itOmega(X)}$ occurs in $X$''
is part of Lemma~\ref{l2.24}.
Denote $P^{\itOmega(X)}$ by $\calP$.
We verify that
$\rfs{Good}_{\calP}(X)
\subseteq
\rfs{Good}(X)$.
Recall that by definition,\break
$\rfs{Good}_{\calP}(X) = \rfs{cl}(\rfs{End}_{\calP}(X))$
and that $\rfs{Good}(X)$ is closed.
So it suffices to show that
$\rfs{End}_{\calP}(X) \subseteq \rfs{Good}(X)$.
Let $e \in \rfs{End}_{\calP}(X)$.
This means that there is $F \subseteq X$ such that
$\pair{F}{e} \in \whatcalX_{\calP}$.
By the definition of $\calP$, we have that $e \in \rfs{Pend}(F)$,
and from $\itOmega(X)$-demonstrativeness of $F$
it follows that $e \in \rfs{Good}(F)$.
The fact $\itOmega(F) = \itOmega(X)$ implies that
$\rfs{Good}(F) \subseteq \rfs{Good}(X)$.
So $e \in \rfs{Good}(X)$.
That is, $\rfs{End}_{\calP}(X) \subseteq \rfs{Good}(X)$
and hence
$\rfs{Good}_{\calP}(X) \subseteq \rfs{Good}(X)$.
\smallskip\hfill\proofend

We need one last property of spaces which are a continuous image of
a compact interval space.

\begin{prop}\label{p2.31}
Let $X \in K_{\srfs{CII}}$.
Then for every infinite cardinal $\lambda$ and a closed subset
$F \subseteq X$:
If $\abs{F} \geq (2^{\lambda})^+$,
then $F$ contains a scattered subspace $H$
such that $\abs{H} = \lambda^+$.
\end{prop}

\noindent{\bf Proof }
(a) A closed subspace of an interval space is an interval space.
This implies that a closed subspace of a member of $K_{\srfs{CII}}$
is a member of $K_{\srfs{CII}}$.
So it suffices to show that if $X \in K_{\srfs{CII}}$
and $\abs{X} \geq (2^{\lambda})^+$,
then $X$ contains a scattered subspace $H$
such that $\abs{H} = \lambda^+$.

Let $\pair{L}{<}$ be a compact linear ordering and
$\fnn{g}{L}{X}$ be continuous and surjective.
There is $A \subseteq L$ such that $\abs{A} = (2^{\lambda})^+$
and $g \nrestriction A$ is $\onetoone$.
By \Erdos\ Rado Theorem, there is $B \subseteq A$
such that $B$ is order isomorphic to $\lambda^+$
or to the reverse ordering of $\lambda^+$. Let $C = \rfs{cl}(B)$
and $H = g(C)$.
Then $\abs{C} = \lambda^+$.
It is obvious that $C$ is homeomorphic to $\lambda^+ + 1$
with is order topology, so $C$ is scattered.
Hence $H = g(C)$ is scattered.
Since $g \nrestriction B$ is $\onetoone$ and $\abs{B} = \lambda^+$,
it follows that $\abs{H} = \lambda^+$.
So $H$ is as desired.
\medskip\hfill\proofend

Let $K$ be the class of all compact Hausdorff spaces $X$
such that
\begin{itemize}
\addtolength{\parskip}{-11pt}
\addtolength{\itemsep}{06pt}
\item[(1)] 
$X \in K_{\srfs{TH}}$.
\item[(2)] 
$X$ is sequentially compact.
\item[(3)] 
For every infinite cardinal $\lambda$ and a closed subset
$F \subseteq X$:
If $\abs{F} \geq (2^{\lambda})^+$,
then $F$ contains a scattered subspace $H$
such that $\abs{H} = \lambda^+$.
\vspace{-2mm}
\end{itemize}

The class $K_{\srfs{CII}}$ is contained in $K$,
and indeed, this has been already shown.
So  the following statement implies Theorem~\ref{t1.2}.

\begin{theorem}\label{t2.28}
For every $X \in K$: if $X$ is a CO space, then $X$ is scattered.
\end{theorem}

\noindent{\bf Proof }
Suppose by contradiction that $X \in K$, $X$ is a CO space
and $X$ is not scattered.
Since $X$ is tightly Hausdorff,
it is collectionwise Hausdorff
and strongly Hausdorff for convergent sequences.
$X$ is also sequentially compact.
Hence by Corollary~\ref{c2.9}, $\itOmega(X)$ is not attained in $X$.

By Proposition~\ref{p2.32-a}, $\calP^{\itOmega(X)}$ is a P-system,
and by \ref{p2.32-b}, $\calP^{\itOmega(X)}$ occurs in~$X$.

Denote $\calP^{\itOmega(X)}$ by $\calP$ and $\abs{\itOmega(X)}$
by $\mu$. Note that $\abs{T_{\calP}} = \mu$.
Then by Corollary~\ref{tightly-haus-c1.5},
$\abs{\rfs{Good}_{\calP}(X)} \geq \beths_4(\mu) \geq (\beths_3(\mu))^+$.
By \ref{p2.32-b}, $\rfs{Good}_{\calP}(X) \subseteq \rfs{Good}(X)$.
Hence $\abs{\rfs{Good}(X)} \geq (\beths_3(\mu))^+$.

By Clause (3) in the definition of $K$,
there is a scattered subspace $F \subseteq \rfs{Good}(X)$
such that $\abs{F} = (\beths_2(\mu))^+$.
From the scatteredness of $F$ it follows that $A \eqdf \rfs{Is}(F)$
is dense in $F$.
So if $\abs{A} \leq \mu$,
then $\abs{F} \leq \beths_2(\mu)$.
It follows that $\abs{A} \geq \mu^+$.

Since $A$ is relatively discrete and $X \in K_{\srfs{TH}}$,
there is a tight Hausdorff system for $A$.
Denote it by $\calU = \setm{U_a}{a \in A}$.
Let $\fnn{\gamma}{A}{\itOmega(X)}$ be a surjection.
Recall that $A \subseteq F \subseteq \rfs{Good}(X)$.
So for every $a \in A$ there is $F_a \in \calS(X)$
such that $F_a \subseteq U_a$ and $\rfs{rk}(F_a) = \gamma(a)$.
Let $H = \rfs{cl}(\bigcup_{a \in A} F_a)$.
Since $\bigcup_{a \in A} F_a \subseteq \rfs{ker}(X)$,
it follows that $H \subseteq \rfs{ker}(X)$.
Also, it is easy to see that $H$ is scattered and that
$\rfs{rk}(H) \geq \itOmega(X)$.
This contradicts the fact that $\itOmega(X)$ is not attained.
\smallskip\hfill\proofend

\noindent
{\bf Proof of Theorem \ref{t1.2} }
It suffices to check that $K_{\srfs{CII}} \subseteq K$.
Let $X \in K_{\srfs{CII}}$.
By Lemma~\ref{l2.3}, $X$ is tightly Hausdorff,
that is, Clause (TH1) of $K_{\srfs{TH}}$ is fulfilled by $X$.
Clause (TH2) is implied by Proposition~\ref{p2.16},
and Clause (TH3) is implied by Proposition~\ref{p2.14}.
By Proposition~\ref{p2.7}, $X$ is sequentially compact,
and Proposition~\ref{p2.31} implies that
$X$ fulfills Clause (3) in the definition of $K$.
So $X \in K$.
\hfill\proofend

\section{Orderability of continuous images of interval spaces}
\label{s3}

In this section we consider the following question.
Suppose that $X$ is a continuous image of a compact interval space,
and $X$ is scattered. Is $X$ an interval space?
The answer to this question is in terms of obstructions.
That is, we define a class of spaces $\calO$, and prove that $X$ is an
interval space iff $X$ has no subspace homeomorphic to a member of $\calO$.
 
For infinite cardinals $\kappa,\lambda$ and $\mu$ we define the
topological space $X_{\kappa,\lambda,\mu}$ as follows.
$X_{\kappa,\lambda,\mu}$ is the quotient of the disjoint union of
the interval spaces $\kappa + 1,\lambda + 1$ and $\mu + 1$,
where the points $\kappa,\lambda$ and $\mu$ are identified.
Say that $\trpl{\kappa}{\lambda}{\mu}$ is a {\it legal triple} if
$\kappa,\lambda$ and $\mu$ are regular cardinals and
$\lambda,\mu > \aleph_0$.
Let\break
$\calT =
\setm{X_{\kappa,\lambda,\mu}}
{\trpl{\kappa}{\lambda}{\mu} \mbox{ is a legal triple}}$.

Let $\lambda$ be an uncountable regular cardinal
and $S \subseteq \lambda$.
Let $\fnn{\vecmu}{S}{\rfs{\rfs{On}}}$
be such that
$\vecmu(\alpha)$ is an uncountable regular cardinal
for every $\alpha \in S$. Denote $\vecmu(\alpha)$ by $\mu_{\alpha}$.
Define the space $ X = X_{\lambda,\vecmu}$ as follows.
Let $\hatomega = \setm{-i - 1}{i \in \omega}$ be the set of negative integers.
The universe of $X$ is
$$
(\lambda + 1) \cup
(\bigcup_{\alpha \in S} \sngltn{\alpha} \times (\mu_{\alpha}
\cup \hatomega)).
$$
An open base $\calB$ of the topology of $X$ consists of the following
sets.
\begin{itemize}
\item[(1)] For every $\alpha \in S$, an open set
$U \subseteq \mu_{\alpha}$ and a subset $V$ of $\hatomega$,\break
$\sngltn{\alpha} \times U,
\sngltn{\alpha} \times V
\in \calB$.
\item[(2)]
Let $W$ be an open subset of $\lambda + 1$
and $\sigma \subseteq W \cap S$ be finite.
\hbox{For every $i \in \sigma$} let $F_i$ be a closed subset of
$\mu_{i} + 1$ not containing $\mu_{i}$ and $G_i$ be a finite subset of
$\hatomega$.
Then
\\
$W \cup
\bigcup_{\alpha \in W \cap S}
\sngltn{\alpha} \times (\mu_{\alpha} \cup \hatomega) -
\bigcup_{i \in \sigma} \sngltn{i} \times (F_i \cup G_i) \in \calB$.
\end{itemize}
Denote $\calB$ by $\calB_{\lambda,\vecmu}$.
\index{N@$\calB_{\lambda,\vecmu}$. A base of $X_{\lambda,\mu}$}
It is left to the reader to check that $X_{\lambda,\vecmu}$ is compact.
Let $\calS =
\setm{X_{\lambda,\vecmu}}
{\rfs{Dom}(\vecmu) \mbox{ is a stationary subset of } \lambda}$
and $\calO = \calT \cup \calS \cup \sngltn{X_{\aleph_1}}$.
(Recall that $X_{\aleph_1}$
is the one point compactification of a discrete space of cardinality $\aleph_1$).

\begin{theorem}\label{t2.1}
Let $X$ be a scattered continuous image of a compact interval space.
Then $X$ is an interval space iff no subset of $X$ is homeomorphic
to a member of $\calO$.
\end{theorem}

\begin{prop}\label{p1.2}
\num{a}
Let $X$ be a closed subspace of a scattered continuous image
of a compact interval space.
Then $X$ is a continuous image of a compact interval space.

\num{b}
If $X$ is a scattered continuous image of a compact interval space,
then there is a 0-dimensional compact interval space $Y$
such that $X$ is a continuous image of $Y$.
\end{prop}

\noindent
{\bf Proof } (a) Suppose that $X$ is a subspace of $Y$,
and $Y$ that is a continuous image of $Z$.
Then the preimage of $X$ in $Z$ is a closed subset of $Z$,
and thus it is an interval space.

(b) Let $\pair{L}{<}$ be a compact chain and $\fnn{f}{L}{X}$
be continuous and onto.
Let $L' = \dbltn{0}{1} \times L$ and $<'$ be the lexicographic order of $L'$.
Then $L'$ is compact and $0$-dimensional. Define $\fnn{f'}{L'}{X}$
by $f'(\pair{i}{a}) = f(a)$. Then $f'$ is continuous
and $\rfs{Rng}(f') = X$.
\hfill\proofend

We need the following theorem from \cite{B}.

\begin{theorem}\label{t1.3}
The following conditions are equivalent.
\begin{itemize}
\item[{\rm(1)}]
$X$ is a scattered continuous image of a compact interval space.
\item[{\rm(2)}]
There is a family $\setm{U_x}{x \in X}$ of clopen subsets of $X$
such that 
\begin{itemize}
\item[{\rm(2.1)}]
For every $x \in X$, $\sngltn{x} = D_{\srfs{rk}^X(x)}(U_x)$.
\item[{\rm(2.2)}] For every $x,y \in X$: if $y \in U_x$,
then $U_y \subseteq U_x$.
\item[{\rm(2.3)}] For every $x,y \in X$: if $y \not\in U_x$
and $x \not\in U_y$, then $U_x \cap U_y = \emptyset$.
\end{itemize}
\end{itemize}
\end{theorem}

\noindent
{\bf Proof }
In \cite{B} Theorem 1.5 it is proved that for every topological space $X$:\break
$X$ is a scattered continuous image of a compact $0$-dimensional interval space,
iff $X$ satisfies Clause (2).
But by Proposition~\ref{p1.2}(b),
$X$ is a scattered continuous image of a compact interval space
iff
$X$ is a scattered continuous image of a compact $0$-dimensional interval space.
So Theorem~\ref{t1.3} follows.
\hfill\proofend

\kern2mm

Let $\calU = \setm{U_x}{x \in X}$ be as in the above theorem.
Then we call $\calU$ a {\it  tree-like clopen system} for $X$.
\index{D@treelike clopen system@@tree-like clopen system}
It is easy to see that if $X$ is scattered and compact and
$\calU$ is a tree-like clopen system for $X$,
then $\calU \cup \setm{X - U}{U \in \calU}$
is a subbase for the topology of $X$.
Let $X$ be a scattered compact space. We say that $X$ is {\it unitary}
if for some $e \in X$, $D_{\srfs{rk}(X)}(X) = \sngltn{e}$.
If $X$ is unitary, then the above $e$ is denoted by $e^X$.
Every scattered compact space $X$ is a finite union
of pairwise disjoint clopen sets $U$ such that
$U$ is unitary and $\rfs{rk}(U) = \rfs{rk}(X)$.

It is clear that if $X$ is a finite union of pairwise disjoint
clopen sets which are interval spaces, then $X$ is an interval space.

Let $\pair{P}{<}$ be a poset and $x \in P$.
We define $P^{< x} = \setm{y \in P}{y < x}$.
The sets $P^{\leq x}$ etc.\ are defined analogously.
\index{N@$P^{< x} = \setm{y \in P}{y < x}$}
Suppose that $\pair{L}{<}$ is a linear ordering and $a \in L$.
We denote the cofinality of $L^{< a}$ by $\rfs{cf}^-_{\pair{L}{<}}(a)$
and the coinitiality of $L^{> a}$ by $\rfs{cf}^+_{\pair{L}{<}}(a)$.
\index{N@$\rfs{cf}^-_{\pair{L}{<}}(a)$ and $\rfs{cf}^+_{\pair{L}{<}}(a)$.
The cofinality of $a$ from the left and the\\
\indent cofinality of $a$ from the right}

We shall also need the following well-known facts.

\begin{prop}\label{p2.3.1}
Let $X$ be a compact Hausdorff space
and $\leq$ be a partial ordering of $X$
such that $\leq$ is a closed subset of $X \times X$.

\num{a} If $C \subseteq X$ is a chain,
then $C$ has a supremum and an infimum.

\num{b} Suppose that $X$ has the following property.
(H1) For every $x,y \in X$: if $x \not\geq y$,
then there are open sets $U$ and $V$ such that
U is an initial segement of $\pair{X}{\leq}$ and $x \in U$,
$V$ is a final segment of $\pair{X}{\leq}$ and $y \in V$,
and $U \cap V = \emptyset$.
\underline{Then} for every chain $C \subseteq X$,
$\sup(C),\inf(C) \in \rfs{cl}^X(C)$.

\num{c} Suppose that $X$ satisfies (H1)
and let $C \subseteq X$ be a chain
such that for every nonempty $A \subseteq C$,
$\sup(A),\inf(A) \in C$.
Then $C$ is closed in $X$,
and the order topology of $C$ coincides with the induced topology of $C$.
\end{prop}

\noindent
{\bf Proof of Theorem \ref{t2.1} }
We first prove that an interval space does not have a subspace homeomorphic
to a member of $\calO$.
Since a closed subspace of an interval space is necessarily an interval space,
it suffices to show that for every $X \in \calO$, $X$ is not an interval space.

The proofs that $X_{\aleph_1}$ is not an interval space,
and that for a legal triple $\trpl{\kappa}{\lambda}{\mu}$
the space $X_{\kappa,\lambda,\mu}$ is not an interval space
are left to the reader.

We show that a space of the type $X_{\lambda,\vecmu}$,
where $\rfs{Dom}(\vecmu)$ is stationary in $\lambda$
is not an interval space.

So let $X_{\lambda,\vecmu} \in \calS$ and suppose by way of contradiction that
$\prec$ is a linear ordering of $X_{\lambda,\vecmu}$
which induces the topology of $X_{\lambda,\vecmu}$.
Let $S = \rfs{Dom}(\vecmu)$ and denote $\vecmu(\alpha)$
by $\mu_{\alpha}$.
Also denote $X_{\lambda,\vecmu}$ by $X$.
Let $<^{\srfs{On}}$ denote the linear ordering of the ordinals.
For an ordinal $\alpha$ let $\tau^{\alpha}$ be the order topology of
$\pair{\alpha}{<^{\srfs{On}} \restriction \alpha}$.
Consider the sets $(\lambda + 1) \cap X^{\tpreceq \lambda}$
and $(\lambda + 1) \cap X^{\tsucceq \lambda}$.
One of them must be of cardinality $\lambda$ and the other of cardinality
$< \lambda$.
This is so since in the interval space  $\pair{\lambda + 1}{\tau^{\lambda + 1}}$
every two closed sets of cardinality $\lambda$
intersect in a set of cardinality $\lambda$.
So we may assume that $\abs{(\lambda + 1) \cap X^{\tpreceq \lambda}} = \lambda$
and $\abs{(\lambda + 1) \cap X^{\tsucceq \lambda}} < \lambda$.
Hence for some $\alpha_0 < \lambda$,
$[\alpha_0,\lambda]^{<^{\trfs{On}}} \subseteq X^{\preceq \lambda}$.

Let $X_{\kappa,\vecnu} \in \calS$
and $C \subseteq \kappa$ be a club.
Define $X_{\kappa,\vecnu} \restriction C$ as follows.
$$
X_{\kappa,\vecnu} \restriction C =
C \cup \sngltn{\kappa} \cup \bigcup_{\alpha \in S \cap C}
\sngltn{\alpha} \times (\vecnu(\alpha) \cup \hatomega).
$$
Then for some $\vecnu\fprime$,
$X_{\kappa,\vecnu} \restriction C \cong X_{\kappa,\vecnu\fprime}$
and $\rfs{Dom}(\vecnu\fprime)$ is stationary in $\kappa$.
So if $X_{\kappa,\vecnu}$ is counter-example,
then $X_{\kappa,\vecnu} \restriction C$ too is a counter-example.
We may thus replace $X_{\lambda,\vecmu}$
by $X_{\lambda,\vecmu} \restriction [\alpha_0,\lambda]^{<^{\trfs{On}}}$,
and assume that $\lambda \subseteq X^{\prec \lambda}$.

We say that $\alpha \in \lambda$ is a bad
if there is $\beta = \beta_{\alpha} <^{\srfs{On}} \alpha$
such that $\alpha \prec \beta$.
Suppose by way of contradiction that the set $B$
of bad points is stationary.
Then the function taking every $\alpha \in B$ to $\beta_{\alpha}$
is constant on an unbounded set.
Let $\gamma$ be this constant value.
Then $\gamma \succeq \lambda$. A contradiction, so $B$ is non-stationary.
Let $C$ be a club disjoint from $B$
and $X_0 = X \restriction C$.
For some $\vecmu\fprime$, $X_0 \cong X_{\lambda,\vecmu\fprime}$
and $\rfs{Dom}(\vecmu\fprime)$ is stationary in $\lambda$.
Replacing $X$ by $X_0$,
we may assume that
$<^{\srfs{On}} \restriction\kern-1.5pt \lambda =
\kern4pt\prec\kern1pt \restriction\kern-1.5pt \lambda$.

Let $\alpha \in S$ be a limit ordinal.
We show that there is $\gamma_{\alpha} \in \mu_{\alpha}$ such that
$\sngltn{\alpha} \times [\gamma_{\alpha},\mu_{\alpha}) \subseteq
X^{\succ \alpha}$.
The subspace $Y = \sngltn{\alpha} \cup (\sngltn{\alpha} \times \mu_{\alpha})$ of $X$
is homeomorphic to $\pair{\mu_{\alpha + 1}}{\tau^{\mu_{\alpha + 1}}}$,
and the subsets $A_1 \eqdf (\sngltn{\alpha} \times \mu_{\alpha}) \cap X^{\tpreceq \alpha}$
and $A_2 \eqdf (\sngltn{\alpha} \times \mu_{\alpha}) \cap X^{\tsucceq \alpha}$ of $Y$
have only one common point in their closures.\break
So one of these sets must have cardinality $\mu_{\alpha}$ and the other must have cardinality
$< \mu_{\alpha}$.
Suppose by contradiction that
$\abs{A_1} = \mu_{\alpha}$. Then $A_1 \cup \sngltn{\mu_{\alpha}}$ is a closed subset of $Y$
of cardinality $\mu_{\alpha}$. Hence $A_1 \cup \sngltn{\mu_{\alpha}} \cong Y$.
It follows that for every $B \subseteq A_1$: if $\abs{B} = \mu_{\alpha}$,
then $\mu_{\alpha} \in \rfs{cl}(B)$.
This implies that
$\rfs{cf}(\pair{X^{\prec \mu_{\alpha}}}{\prec \restriction X^{\prec \mu_{\alpha}}}) =
\mu_{\alpha} > \aleph_0$.
The sets $A_1$ and $\mu_{\alpha} \cap X^{\prec \mu_{\alpha}}$
are closed unbounded subsets of $X^{\prec \mu_{\alpha}}$, and they are disjoint.
A contradiction.
So $\abs{A_1} < \mu_{\alpha}$, and hence there is
$\gamma_{\alpha} \in \mu_{\alpha}$ such that
$\sngltn{\alpha} \times [\gamma_{\alpha},\mu_{\alpha}) \subseteq
X^{\succ \alpha}$.

Let $\alpha \in S$ be a limit ordinal.
It follows that
$\rfs{cf}^+_{\pair{X}{\prec}}(\alpha) = \mu_{\alpha} > \aleph_0$.
However, $\sngltn{\alpha} \times \hatomega$ is an $\omega$-sequence
converging to $\alpha$.
We conclude that
$(\sngltn{\alpha} \times \hatomega) \cap X^{\succ \alpha}$ is finite.
Since $\prec$ and $<^{\srfs{On}}$ coincide on $\lambda$,
the subset $\alpha$ of $X$ is cofinal in $X^{\prec \alpha}$.
So there are $\gamma_{\alpha} < \alpha$
and $a_{\alpha} \in \sngltn{\alpha} \times \hatomega$
such that $a_{\alpha} \prec \gamma_{\alpha}$.
Then $\alpha \mapsto \gamma_{\alpha}$ is a regressive function
defined on a stationary subset of $\lambda$.
Let $\gamma$ be such that $\gamma_{\alpha} = \gamma$ for an unbounded
set of $\alpha$'s. Denote this unbounded set by $D$.
Then by the definition of the topology of $X_{\lambda,\vecmu}$,
$\lambda \in \rfs{acc}(\setm{a_{\alpha}}{\alpha \in D})$.
But $a_\alpha \prec \gamma \prec \lambda$ for every $\alpha \in D$.
A contradiction. So $X_{\lambda,\vecmu}$ is not homeomorphic to
an interval space.

We prove the other direction of the theorem by induction on
$\rfs{rk}(X)$. The statement of the induction hypothesis requires
some prepartion.

Let $X$ be a scattered continuous image of a compact interval space
and $\calU = \setm{U_x}{x \in X}$ be a tree-like clopen system for $X$.
For $x,y \in X$ define $x \leq_{\calU} y$ if $x \in U_y$.
Clearly $\leq_{\calU}$ is a partial ordering of $X$.
We say that $\pair{X}{\calU}$ is simple if $\pair{X}{\leq_{\calU}}$
has a maximum $e^{\calU}$,
and there are an uncountable regular cardinal $\lambda$
and a strictly increasing sequence
$\setm{x_{\alpha}}{\alpha < \lambda}$
in $\pair{X}{\leq_{\calU}}$ such that
$X - \sngltn{e^{\calU}} = \bigcup_{\alpha < \lambda} U_{x_{\alpha}}$.

We shall prove by induction on $\alpha$ the following statement.
\begin{itemize}
\item[$(*)_{\alpha}$]
If $X$ is a scattered continuous image of a compact interval space,
$\rfs{rk}(X) = \alpha$
and no subspace of $X$ is homeomorphic to a member of $\calO$,
then $X$ is an interval space.
If in addition, $\pair{X}{\calU}$ is simple,
then there is a linear ordering $\leq_X$ of $X$ such that
the order topology of $\leq_X$ is the topology of $X$ and
$e^{\calU}$ is the maximum of $\pair{X}{\leq_X}$.
\end{itemize}

Denote $(*)_{< \alpha} \equiv \bigwedge_{\beta < \alpha} (*)_{\beta}$.
It is trivial that $(*)_0$ holds.
We shall prove that if $\alpha > 0$,
and for every $(*)_{< \alpha}$ holds,
then $(*)_{\alpha}$ holds.

A poset $\pair{P}{\leq}$ is called a {\it reverse tree}
if $P^{>x}$ is a chain for every $x \in X$.
A subset $D$ of a poset $P$ is {\it directed}, if for every $a,b \in D$
there is $c \in D$ such that $a,b \leq c$.
We say that $D$ is 
{\it principal} if for some $d \in P$, $D = P^{\leq d}$.
\index{D@principal}
We say that $D$ is {\it generated} by $A$ if
$D =
\setm{p \in P}{\mbox{ there is } a \in A \mbox{ such that } p \leq a}$.
For $x \in P$ we set
$\calD_x =
\setm{D}{D \mbox{ is a maximal directed subset of } P^{< x}}$.
We leave it to the reader to check that if $P$ is a reverse tree,
$x \in P$ and $D_1,D_2 \in \calD_x$ are distinct,
then $D_1 \cap D_2 = \emptyset$.
Also, if $P$ is a reverse tree, $x \in P$
and $D \in \calD_x$, then there is a chain which generates $D$.
Let $D \in \calD_x$. If $D$ has a maximum,
then the cofinality of $D$ is defined to be $1$.
Otherwise, there is unique regular cardinal $\nu$ such that $D$
is generated by a chain of type $\nu$.
We denote this $\nu$ by $\rfs{cf}(D)$.
A subset $A$ of a poset $P$ is {\it unbounded},
if there is no $p \in P$ such that $p \geq a$ for every $a \in A$.

Let $X$ be a scattered continuous image of a compact interval space.
Let $\calU = \setm{U_x}{x \in X}$ be a tree-like clopen system for $X$.
Clearly, $\pair{X}{\leq_{\calU}}$ is a reverse tree.
Note also that $\leq_{\calU}$ is a closed subset of $X \times X$.
This is so, since
$\not\leq_{\calU} \kern8pt=\kern8pt \bigcup_{x \in X} (X - U_x) \times U_x$.
It is trivial that $\trpl{X}{\tau^X}{\leq_{\calU}}$ satisfies
Property~(H1) from Proposition~\ref{p2.3.1}.
The set of maximal points in $\pair{X}{\leq_{\calU}}$ is finite.
Suppose otherwise. Let $y$ be an accumulation point of the set of
maximal points.
Since $U_y$ is a neighborhood of $y$ it contains a maximal point
$z \neq y$.
But then $z <_{\calU} y$. A contradiction.
It follows that the set of maximal points is finite.
Also, $X = \bigcup \setm{U_x}{x \mbox{ is a maximal point of } X}$.
If $V$ is a clopen subset of $X$, then $\setm{V \cap U_x}{x \in V}$
is a tree-like clopen system for $V$.
It follows that if $X$ is a scattered continuous image of
a compact interval space, then there are
$\pair{X_1}{\calU_1},\ldots,\pair{X_n}{\calU_n}$
such that $\fsetn{X_1}{X_n}$ is a partition of $X$ into clopen sets,
and $\pair{X_i}{\leq_{\calU_i}}$ has a maximum for every $i \leq n$.
It is trivial that if $X$ has a maximum $e^{\calU}$,
then $X$ is unitary and that $e^{\calU} = e^X$.

\noindent
{\bf Claim 1 }
Let $X$ be a scattered continuous image of a compact interval space
such that no subspace of $X$ is homeomorphic to a member of $\calO$.
Let $x \in X$.
Then the following facts hold.
\begin{itemize}
\item[(1)] $\abs{\calD_x} \leq \aleph_0$.
\item[(2)] 
If there are distinct $D_0,D_1 \in \calD_x$
such that $\rfs{cf}(D_0),\rfs{cf}(D_1) \geq \aleph_1$,
then $\calD_x$ is finite and every
member of $\calD_x$ other than $D_0$ and $D_1$ has a maximum.
\end{itemize}
{\bf Proof }
(1) The proof relies on the fact that $X_{\aleph_1}$ is not embeddable in~$X$.
Let $A \subseteq U_x$ be infinite, and assume that for every $D \in \calD_x$,
$\abs{A \cap D} \leq 1$. We show that $A$ is discrete and that
$A \cup \sngltn{x}$ is the one point compactification of $A$.
If $y,z \in A$ are distinct, then $z \not\leq_{\calU} y$. That is,
$z \not\in U_y$. So $U_y$ is a neighborhod of $y$ disjoint from
$A - \sngltn{y}$.
We show that every neighborhood $V$ of $x$ contains all
but finitely many members of $A$.
We may assume that $V$ has the form
$U_x - \bigcup_{y \in \sigma} U_y$,
where $\sigma$ is a finite subset of $U_x$.
If $y \in \sigma$, then $y \leq_{\calU} x$ and so
$\abs{U_y \cap A} \leq 1$.
So $A - V$ is finite.
Hence $x \in \rfs{acc}(A)$.
Let $y \in X - \sngltn{x}$.
We show that $y \not\in \rfs{acc}(A)$.
If $y \not\in U_x$, then $U_y \cap U_x = \emptyset$.
Hence $U_y \cap A = \emptyset$.
If $y \in U_x$, then $\abs{U_y \cap A} \leq 1$.
So $y \not\in \rfs{acc}(A)$.
We have shown that $\rfs{acc}(A) = \sngltn{x}$.
So $A \cup \sngltn{x}$ is the one point compactification of $A$.
Since $X_{\aleph_1}$ is not embeddable in $X$,
$\abs{A} \leq \aleph_0$, and since $\calD_x$ is a pairwise disjoint family,
$\abs{\calD_x} = \aleph_0$.

(2) The proof relies on the fact that no member of $\calT$ is embeddable in $X$.
Suppose that $D_0,D_1 \in \calD_x$, $D_0 \neq D_1$
and $\rfs{cf}(D_0),\rfs{cf}(D_1) \geq \aleph_1$.
For $i = 0,1$ let $E_i$ be a chain which generates $D_i$ and
such that the order type of $E_i$ is a regular cardinal $\lambda_i$.
We may assume that for every nonempty bounded $B \subseteq E_i$,
$\sup(B) \in E_i$.
Suppose by way of contradiction that $\calD_x$ is infinite.
Since $\calD_x$ is a pairwise disjoint family,
it follows that there is a countably infinite set
$A \subseteq U_x$ such that
$\abs{A \cap D} \leq 1$ for every $D \in \calD_x$.
We show that $Y \eqdf A \cup E_1 \cup E_2 \cup \sngltn{x}$
is homeomorphic to $X_{\aleph_0,\lambda_1,\lambda_2}$.
For $i = 1,2$, $E_i \cup \sngltn{x}$ is a chain in $X$
closed under infima and suprema.
So by Proposition~\ref{p2.3.1}(c),
its induced topology $\tau_i$ coincides with its order topology.
That is, $\pair{E_i \cup \sngltn{x}}{\tau_i} \cong \lambda_i + 1$.
It is also clear that $A \cup \sngltn{x}$
with its induced topology is homeomorphic to $\omega + 1$.
These facts imply that
$\pair{Y}{\tau^X \restriction Y} \cong
X_{\aleph_0,\lambda_1,\lambda_2} \in \calT$.
A contradiction.

Suppose by contradiction that
$D \in \calD_x - \dbltn{D_0}{D_1}$ and $D$ is nonprincipal.
Let $E$ be unbounded chain in $D$ such that the order type of $E$
is a regular cardinal $\mu$.
We may assume that for every nonempty bounded $B \subseteq E$,
$\sup(B) \in E$.
Just as in the previous argument we conclude that
$E \cup E_1 \cup E_2 \cup \sngltn{x} \cong
X_{\mu,\lambda_1,\lambda_2} \in \calT$.
A contradiction.
We have proved Claim 1.
\smallskip

Suppose that
$(*)_{< \alpha_0}$ holds,
and we prove $(*)_{\alpha_0}$.
Let $X$ be a scattered continuous image of a compact interval space,
suppose that $\rfs{rk}(X) = \alpha_0$,
and that no subspace of $X$ is homeomorphic to a member of $\calO$.
Let $\calU$ be a tree-like clopen system for $X$.
Since $X$ can be partitioned into finitely many clopen sets $X_i$
with tree-like systems $\calU_i$
such that each $\pair{X_i}{\leq_{\calU_i}}$ has a maximum,
we may assume that $X$ has a maximum.
We deal separately with three cases.

{\bf Case 1 }
Assume that $\pair{X - \sngltn{e}}{\leq_{\calU} \restriction (X - \sngltn{e})}$
contains a chain $C$ with uncountable cofinality,
such that $X - \sngltn{e} = \bigcup_{x \in C} U_x$.
In such a situation it is required that we prove that there is a linear
ordering $\leq_X$ of $X$ which induces the topology of $X$
and in which $e = \max(\pair{X}{\leq_X})$.

We may assume that $C$
is order isomorphic to an uncountable regular cardinal $\lambda$.
We may further assume that for every nonempty bounded $A \subseteq C$,
$\sup(A) \in C$.
Obviously, $C \cup \sngltn{e}$ is order isomorphic to $\lambda + 1$.
So by proposition~\ref{p2.3.1},
the induced topology on $C \cup \sngltn{e}$
coincides with the order topology of $C \cup \sngltn{e}$.
So the order isomorphism between $\lambda + 1$ and $C \cup \sngltn{e}$
is a homeomorphism.
Let
$$
\alpha \mapsto y_{\alpha}, \ \alpha \leq \lambda,
$$
be the isomorphism between $\lambda + 1$ and $C$.

Let $J_{\alpha} = \bigcup_{\beta < \alpha} U_{y_{\beta}}$.
We check that $J_{\alpha} \cup \sngltn{y_{\alpha}}$ is closed.
If $\alpha = \beta + 1$, then $J_{\alpha} = U_{\beta}$.
So $J_{\alpha} \cup \sngltn{y_{\alpha}}$ is closed.
Suppose $\alpha$ is a limit.
Since $U_{y_{\alpha}}$ is closed,
$\rfs{cl}(J_{\alpha}) \subseteq U_{y_{\alpha}}$.
Let
$x \in U_{y_{\alpha}} - J_{\alpha} - \sngltn{y_{\alpha}}$,
$U_x$ is an open neighborhood of $x$.
Suppose by contradiction that $U_x \cap J_{\alpha} \neq \emptyset$.
Then for some $\beta < \alpha$,
$U_x \cap U_{\beta} \neq \emptyset$,
and hence for every $\gamma \geq \beta$, $x$ and $y_{\gamma}$
are comparable in $\leq_{\calU}$.
But $x <_{\calU} y_{\alpha}$
and hence $x \not\geq_{\calU} y_{\alpha}$.
Since $y_{\alpha} = \sup^{<_{\calU}}(\setm{y_{\gamma}}{\beta \leq \gamma < \alpha})$,
there is $\delta < \alpha$ such that $x \not\geq_{\calU} y_{\delta}$.
So $x <_{\calU} y_{\delta}$.
So $x \in U_{\delta} \subseteq J_{\alpha}$. A contradiction.
Hence $U_x \cap U_{\beta} = \emptyset$.
So $J_{\alpha} \cup \sngltn{y_{\alpha}}$ is closed.

Let $\alpha < \lambda$. We say that $\alpha$ is inconvenient
(with respect to the sequence $\setm{y_{\alpha}}{\alpha < \lambda}$), if $\alpha$
is a limit ordinal, and there are an uncountable regular cardinal
$\mu_{\alpha}$ and disjoint sets
$Y_{\alpha},Z_{\alpha} \subseteq U_{y_{\alpha}} - J_{\alpha} - \sngltn{y_{\alpha}}$
such that
\begin{itemize}
\addtolength{\parskip}{-11pt}
\addtolength{\itemsep}{06pt}
\item[(1)] 
$Y_{\alpha}$ is discrete,
$Y_{\alpha} \cup \sngltn{y_{\alpha}}$
is homeomorphic to $\omega + 1$,
\item[(2)] 
$Z_{\alpha}$ is homeomorphic to $\mu_{\alpha}$
and $Z_{\alpha} \cup \sngltn{y_{\alpha}}$
is homeomorphic to $\mu_{\alpha} + 1$.
\vspace{-05.7pt}
\end{itemize}
Let $S$ be the set of inconvenient ordinals,
$\vecmu = \setm{\mu_{\alpha}}{\alpha \in S}$
and
$$
Y = \setm{y_{\alpha}}{\alpha \leq \lambda} \cup
\bigcup_{\alpha \in S} Y_{\alpha} \cup
\bigcup_{\alpha \in S} Z_{\alpha}.
$$

{\bf Claim 2 }
$Y \cong X_{\lambda,\vecmu}$.

{\bf Proof }
At first we verify that $Y$ is closed.
Let $x \in \rfs{cl}(Y)$. Let $\alpha$ be the first ordinal such that
$x \in U_{y_{\alpha}}$. So $x \not\in J_{\alpha}$.
Hence $U_x$ is a neighborhood of $x$ disjoint from
$J_{\alpha} \cup (X - U_{y_{\alpha}})$.
The complement of this set in $X$ is $U_{y_{\alpha}} - J_{\alpha}$
and
$Y \cap (U_{y_{\alpha}} - J_{\alpha}) = 
\sngltn{y_{\alpha}} \cup Y_{\alpha} \cup Z_{\alpha}$.
So
$x \in \rfs{cl}(\sngltn{y_{\alpha}} \cup Y_{\alpha} \cup Z_{\alpha})$.
But $\sngltn{y_{\alpha}} \cup Y_{\alpha} \cup Z_{\alpha}$
is closed. So
$x \in \sngltn{y_{\alpha}} \cup Y_{\alpha} \cup Z_{\alpha} \subseteq Y$.
Hence $Y$ is closed.

For $\alpha \in S$
let $\setm{z_{\alpha,i}}{i < \mu_{\alpha}}$ be an enumeration of
$Z_{\alpha}$ such that the function
$i \mapsto z_{\alpha,i}, \ i < \mu_{\alpha}$, is a homeomorphism
between $\mu_{\alpha}$ and $Z_{\alpha}$,
and let $\setm{y_{\alpha,i}}{i < \omega}$ be a $\onetoone$ enumeration
of $Y_{\alpha}$.
Define $\fnn{\psi}{X_{\lambda,\vecmu}}{Y}$ as follows:
\begin{itemize}
\item[(1)]
$\psi(\alpha) = y_{\alpha}, \ \alpha \leq \lambda$;
\item[(2)]
$\psi(\pair{\alpha}{i}) =
y_{\alpha,i}, \ \alpha \in S, \, i \in \mu_{\alpha}$;
\item[(3)]
$\psi(\pair{\alpha}{-i - 1}) =
z_{\alpha,i}, \ \alpha \in S, \, i \in \omega$.
\end{itemize}
Clearly, $\psi$ is a bijection.
We prove that $\psi$ is a homeomorphism between $X_{\lambda,\vecmu}$ and $Y$.
Since both $X_{\lambda,\vecmu}$ and $Y$ are compact,
it suffices to show
that for every\break
$B \in \calB_{\lambda,\vecmu}$, $\psi(B)$ is open in $Y$.

We  first show that $Y_{\alpha}$ and $Z_{\alpha}$ are open in $Y$.
Note that
$Y_{\alpha} \cup Z_{\alpha} =
Y \cap (U_{y_{\alpha}} - (J_{\alpha} \cup \sngltn{y_{\alpha}}))$.
Since
$U_{y_{\alpha}} - (J_{\alpha} \cup \sngltn{y_{\alpha}})$
is open in $X$,
it follows that $Y_{\alpha} \cup Z_{\alpha}$ is open in $Y$.
Both $Z_{\alpha} \cup \sngltn{y_{\alpha}}$
and $Y_{\alpha} \cup \sngltn{y_{\alpha}}$
are compact and hence closed in $X$.
So they are closed in $Y$.
Since
$Y_{\alpha} =
(Y_{\alpha} \cup Z_{\alpha}) - (Z_{\alpha} \cup \sngltn{y_{\alpha}})$,
it follows that $Y_{\alpha}$ is open in $Y$.
Similarly, $Z_{\alpha}$ is open in $Y$, because
$Z_{\alpha} = 
(Y_{\alpha} \cup Z_{\alpha}) - (Y_{\alpha} \cup \sngltn{y_{\alpha}})$.

Let $B = \sngltn{\alpha} \times V \in \calB_{\lambda,\vecmu}$,
where $V$ is an open subset of $\mu_{\alpha}$.
Since $\psi \nrestriction \sngltn{\alpha} \times \mu_{\alpha}$
is a homeomorphism onto $Z_{\alpha}$,
$\psi(B)$ is open in $Z_{\alpha}$. So it is open in $Y$.
Similarly,
if $B = \sngltn{\alpha} \times V \in \calB_{\lambda,\vecmu}$,
where $V \subseteq \hatomega$,
then $\psi(B)$ is open in $Y_{\alpha}$. So it is open in~$Y$.

Let $W$ be an open subset of $\lambda + 1$
and $\sigma \subseteq W \cap S$ be finite.
For every $i \in \sigma$ let $F_i$ be a closed subset of
$\mu_{i} + 1$ not containing $\mu_{i}$ and $G_i$ be a finite subset of
$\hatomega$.
Let $B = W \cup
\bigcup_{\alpha \in W \cap S}
\sngltn{\alpha} \times (\mu_{\alpha} \cup \hatomega) -
\bigcup_{i \in \sigma} \sngltn{i} \times (F_i \cup G_i)$.\break
It remains to show that when $B$ has this form,
then $\psi(B)$ is open in~$Y$.

For $i \in \sigma$, $F_i$ is compact in $\mu_i$.
So $\setm{z_{i,j}}{j \in F_i}$ is compact in $Z_i$. So it is closed
in $Y$.
But $\psi(\sngltn{i} \times F_i) = \setm{z_{i,j}}{j \in F_i}$.
Hence $\psi(\sngltn{i} \times F_i)$ is closed in $Y$.
Also, $\psi(\sngltn{i} \times G_i)$ is finite and hence closed in $Y$.
So $\psi(\bigcup_{i \in \sigma} (F_i \cup G_i)$ is closed in $Y$.

Set
$W' = W \cup
\bigcup_{\alpha \in W \cap S}
\sngltn{\alpha} \times (\mu_{\alpha} \cup \hatomega)$.
Then $B = W' - \bigcup_{i \in \sigma} (F_i \cup G_i)$.
We have already shown that $\psi(\bigcup_{i \in \sigma} (F_i \cup G_i))$ is closed
in $Y$.
So it remains to show that $\psi(W')$ is open in $Y$.

We may assume that $W$ is an open convex subset of $\lambda + 1$.
Let us first deal with the case that $W = (\beta,\gamma)$,
where $\beta,\gamma \in \lambda + 1$. Then
$$
\psi(W')
=
\setm{y_{\alpha}}{\alpha \in (\beta,\gamma)} \cup
\bigcup_{\alpha \in (\beta,\gamma) \cap S} (Y_{\alpha} \cup Z_{\alpha})
=
Y \cap (U_{y_{\gamma}} - U_{y_{\beta}}).
$$
So 
$\psi(W')$ is open in $Y$.

If $W$ has the form $(\beta,\lambda]$.
Then $\psi(W') = Y - U_{y_{\beta}}$.
If $W = [0,\gamma)$, then $\psi(W') = Y \cap U_{y_{\gamma}}$,
and finally, if $W = \lambda + 1$, then $\psi(W') = Y$.
In all cases $\psi(W')$ is open in $Y$. Hence $\psi(B)$ is open in $Y$.

Since $\psi$ takes all members of an open base of $X_{\lambda,\vecmu}$
to open subsets of $Y$
and both $X_{\lambda,\vecmu}$ and $Y$ are compact Hausdorff,
$\psi$ is a homeomorphism between $X_{\lambda,\vecmu}$ and $Y$.
This proves Claim 2.
\smallskip

We found that $Y$ is homeomorphic to $X_{\lambda,\vecmu}$.
However, no subspace of $X$ is homeomorphic to a member of $\calS$.
So $X_{\lambda,\vecmu} \not\in \calS$.
This implies that the set $S$ of inconvenient ordinals is non-stationary.
Let $A$ be a closed and unbounded subset of $\lambda$ disjoint from $S$,
and let $\setm{x_{\alpha}}{\alpha < \lambda}$ be a strictly
increasing enumeration of $\setm{y_{\beta}}{\beta \in A}$.
Denote $e$ by $x_{\lambda}$.
The function $\alpha \mapsto x_{\alpha}, \, \alpha \leq \lambda$
is again a homeomorphism.

We claim that
there are no inconvenient ordinals with respect to\break
$\setm{x_{\alpha}}{\alpha < \lambda}$.
Let $\alpha \in \lambda$ be a limit ordinal. There is $\beta$
such that $x_{\alpha} = y_{\beta}$.
Clearly, $\beta$ is a limit ordinal.
Hence
$U_{x_{\alpha}} - \bigcup_{\gamma < \alpha} U_{x_{\gamma}} =
U_{y_{\beta}} - \bigcup_{\gamma < \beta} U_{y_{\gamma}}$.
Since $\beta$ is not inconvenient with respect to
$\setm{y_{\gamma}}{\gamma < \lambda}$ it follows that
$\alpha$ is not inconvenient with respect to~$\setm{x_{\gamma}}{\gamma < \lambda}$.

We shall define by induction on $\alpha \leq \lambda$ linear orderings
$\leq_{\alpha}$ of $U_{x_{\alpha}}$.
Since $e = x_{\lambda}$ and $U_e = X$, the ordering $\leq_{\lambda}$
is an ordering of $X$.
This will be the ordering required in the theorem.

Denote $U_{x_{\alpha}}$ and $\calD_{x_{\alpha}}$
by $U_{\alpha}$ and $\calD_{\alpha}$ respectively,
and define $I_{\alpha} = \bigcup_{\beta < \alpha} U_{x_{\beta}}$.
We need the following facts.
\begin{itemize}
\item[(3)] Suppose that $\alpha$ is a limit ordinal.
Then $I_{\alpha} \in \calD_{\alpha}$.
{\thickmuskip=2mu \medmuskip=1mu \thinmuskip=1mu 
Also, there do not exist distinct
$D_1,D_2 \in \calD_{\alpha} - \sngltn{I_{\alpha}}$
such that $\rfs{cf}(D_1) \geq \aleph_0$
and $\rfs{cf}(D_2) \geq \aleph_1$.
}
\item[(4)] Suppose that $\alpha$ is a limit ordinal.
If $\calD_{\alpha}$ is infinite, then for every
$D \in \calD_{\alpha} - \sngltn{I_{\alpha}}$,
$\rfs{cf}(D) \leq \aleph_0$.
\end{itemize}

The proof of (3) relies on the fact that $\alpha$ is not inconvenient.
For suppose by contradiction that
$D_1,D_2 \in \calD_{\alpha} - \sngltn{I_{\alpha}}$,
$\rfs{cf}(D_1) \geq \aleph_0$, $\rfs{cf}(D_2) \geq \aleph_1$.
Let $E_1,E_2$ be chains which generate $D_1$ and $D_2$ respectively,
and such that the order types of $E_1$ and $E_2$ are regular cardinals.
Suppose further that $E_i \cup \sngltn{x_{\alpha}}$ is closed under suprema
in $\pair{X}{\leq_{\calU}}$ for $i = 1,2$.
Let $E$ be a closed and unbounded subset of $\setm{x_{\beta}}{\beta < \alpha}$
with order type which is a regular cardinal.
It follows that if $\rfs{cf}(D_1),\rfs{cf}(D_2) \geq \aleph_1$,
then
$E \cup E_1 \cup E_2 \cup \sngltn{x_{\alpha}} \in \calT$.
If
$\rfs{cf}(D_1) = \aleph_0$
and
$\rfs{cf}(D_2) \geq \aleph_1$,
then setting $Y = E_1$ and $Z = E_2$ shows that $\alpha$ is inconvenient.
We have proved (3)

The proof of (4) relies on the fact that $\alpha$ is not inconvenient.
Suppose by contradiction that $\calD_{\alpha}$ is infinite
and for some $D \in \calD_{\alpha} - \sngltn{I_{\alpha}}$,
$\rfs{cf}(D) \geq \aleph_1$.
Let $Y \subseteq U_{\alpha} - I_{\alpha} - D$ be a countably infinite set
such that for every $D' \in \calD_{\alpha}$,
$\abs{Y \cap D'} \leq 1$.
Let $Z$ be a chain which generates $D$
whose order type is a regular cardinal
and such that $Z \cup \sngltn{x_{\alpha}}$ is closed under suprema.
Then the pair $Y,Z$ is an evidence that $\alpha$ is inconvenient.
We have proved (4)

If $\pair{W}{\zeta}$ is a topological space and $A \subseteq W$,
denote the relative topology that $A$ inherits from $\pair{W}{\zeta}$
by $\zeta \nrestriction A$.
If $\kern2pt\preceq\kern2pt$ is a linear ordering of a set $A$,
denote by $\tau^{\preceq}$ the order topology of $\pair{A}{\preceq}$.
Let $\rho$ denote the topology $\tau^X$ of~$X$.

We now define by induction on $\alpha \leq \lambda$
the linear ordering $\leq_{\alpha}$ of $U_{\alpha}$.
We assume by induction that
\begin{enumerate}
\addtolength{\parskip}{-08pt}
\addtolength{\itemsep}{03pt}
\item[(I1)]
$\tau^{\leq_{\alpha}} = \rho \nrestriction U_{\alpha}$.
\item[(I2)]
If $\beta < \gamma$,
then $\leq_{\beta} \kern8pt\subseteq\kern8pt \leq_{\gamma}$.
\item[(I3)]
If $\beta < \gamma$,
then $U_{\beta}$ is an initial segment of
$\pair{U_{\gamma}}{\leq_{\gamma}}$.
\vspace{-05.7pt}
\end{enumerate}
Since $\rfs{rk}(x_0) < \rfs{rk}(e) = \alpha_0$,
we may apply the induction hypothesis
that $(*)_{< \alpha_0}$ holds to $x_0$.
Let $\leq_0$ be a linear ordering of $U_0$ which induces the
relative topology of $U_0$.
Suppose that $\leq_{\beta}$ has been defined.
Since
$\rfs{rk}(U_{\beta + 1} - U_{\beta}) < \rfs{rk}(e) = \alpha_0$,
by the induction hypothesis, there is a linear ordering $\leq'$ of
$U_{\beta + 1} - U_{\beta}$ which induces the relative
topology of $U_{\beta + 1} - U_{\beta}$.
Then $\kern8pt\leq_{\beta} \kern8pt\cup\kern8pt \leq'
\kern8pt\cup\kern8pt
U_{\beta} \times (U_{\beta + 1} - U_{\beta})$
is a linear ordering of $U_{\beta + 1}$ which satisfies the
induction hypotheses.

Suppose that $\delta$ is a limit ordinal and
$\leq_{\beta}$ has been defined for every $\beta < \delta$.
Let $\leq'_{\delta} \kern8pt=\kern8pt
\bigcup_{\beta < \delta} \leq_{\beta}$.
So $\leq'_{\delta}$ is a linear ordering of $I_{\delta}$.

{\bf Case 1.1 } There is
$D \in \calD_{\delta} - \sngltn{I_{\delta}}$ such that
$\rfs{cf}(D) \geq \aleph_1$.
By (4), $\calD_{\delta}$ is finite, and by (3) every
$D' \in \calD_{\delta} - \dbltn{I_{\delta}}{D}$ has a maximum.
Let
$\sigma =
\setm{\max(D')}{D' \in \calD_{\delta} - \dbltn{I_{\delta}}{D}}$.
For every $x \in \sigma$, $\rfs{rk}(x) < \alpha_0$, so by the
induction hypothesis, there is a linear ordering $\leq_x$ of
$U_x$ which induces the topology of $U_x$.
Also, let $\leq_{\sigma}$ be a linear ordering of $\sigma$.
We claim that $\delta \neq \lambda$.
This is so since $I_{\delta}$ and $D$
are distinct maximal directed sets in $\calD_{\delta}$,
so there cannot be a chain $I$ in $X^{<_{\calU}\kern1pt x_{\delta}}$ such that
$\bigcup_{x \in I} U_x = X^{<_{\calU}\kern1pt x_{\delta}}$,
and we assumed that such an $I$ exists for $x_{\lambda}$.
Let $Z = D \cup \sngltn{x_{\delta}}$.
It is easy to check that $Z$ is closed in $X$.
So $Z$ is a scattered continuous image of a compact interval space.
Since $Z \subseteq U_{\delta}$ and 
$\rfs{rk}(U_{\delta}) = \rfs{rk}(x_{\delta}) < \alpha_0$,
it follows that $\rfs{rk}(Z) < \alpha_0$. 
Let $\calU^Z = \setm{U^Z_x}{x \in Z}$,
where $U^Z_x$ is defined as follows: if $x \neq x_{\delta}$,
then $U^Z_x = U_x$ and $U^Z_{x_{\delta}} = Z$.
Then $\calU^Z$ is a tree-like clopen system for $Z$.
Let $J$ be an unbounded chain in $D$.
Then $Z - \sngltn{x_{\delta}} = D = \bigcup_{x \in J} U^Z_x$.
We assumed that $(*)_{< \alpha_0}$ holds.
So there is a linear ordering $\leq_Z$ of $Z$ such that
$\tau^{\leq_Z} = \rho \restriction Z$ and such that $x_{\delta} = \min(\pair{Z}{\leq_Z})$.
We define the required linear ordering $\leq_{\delta}$ of $U_{\delta}$ as follows.
\begin{itemize}
\item[(1)] $I_{\delta} \leq_{\delta} Z$.
(This means: for every $a \in I_{\delta}$ and $b \in Z$, $a \leq_{\delta} b$).
\item[(2)] For every $x \in \sigma$, $Z \leq_{\delta} U_x$.
\item[(3)] For every $x,y \in \sigma$: if $x \leq_{\sigma} y$,
then $U_x \leq_{\delta} U_y$.
\item[(4)]
$\leq_{\delta} \restriction\kern-3pt I_{\delta}
\kern3pt=\kern6pt \leq'_{\delta}$,\ \ %
$\leq_{\delta} \restriction\kern-3pt Z \kern3pt=\kern6pt \leq_Z$\ \ %
and\ \ %
$\leq_{\delta} \restriction\kern-3pt U_x
\kern3pt=\kern6pt \leq_x$ \ for every $x \in \sigma$.
\end{itemize}
Clearly, $\leq_{\delta}$ is a linear ordering of $U_{\delta}$.
Recall that $\tau^{\leq_{\delta}}$
denotes the order topology of $\pair{U_{\delta}}{\leq_{\delta}}$.
We show that $\tau^{\leq_{\delta}} = \rho \nrestriction U_{\delta}$.

Let 
$W = I_{\delta} \cup \sngltn{x_{\delta}}$
and
$\calX \eqdf \dbltn{W}{Z} \cup \setm{U_x}{x \in \sigma}$.
Clearly, $\calX$ is a finite cover of $U_{\delta}$.
We shall argue as follows. At first we check that every member
of $\calX$ is closed in both $\tau^{\leq_{\delta}}$ and $\rho \nrestriction U_{\delta}$.
Then we show
$(*)$ For every $T \in \calX$,
$\tau^{\leq_{\delta}} \nrestriction T = \rho \nrestriction T$.

Notice that if $\calF$ is a finite cover of a space
$\pair{S}{\eta}$ consisting
of closed sets, then for every $V \subseteq S$:
$V \in \eta$ iff $V \cap F \in \eta \nrestriction F$ for every
$F \in \calF$.
That is $\setm{\eta \nrestriction F}{F \in \calF}$ determines $\eta$.
Hence $(*)$ implies that $\tau^{\leq_{\delta}} = \rho$.

We show that for every $x \in \sigma$,
$\tau^{\leq_{\delta}} \restriction U_x = \rho \restriction U_x$.
Recall that $\leq_x$ is a linear ordering of $U_x$ such that
$\tau^{\leq_x} = \rho \restriction U_x$.
Also, $U_x$ is a closed interval in $\pair{U_{\delta}}{\leq_{\delta}}$
and $\leq_{\delta} \restriction U_x =\kern6pt \leq_x$.
So
$$
\tau^{\leq_{\delta}} \restriction U_x =
\tau^{\leq_{\delta} \restriction U_x} = \tau^{\leq_x} = \rho \restriction U_x.
$$

An identical argument shows that
$\tau^{\leq_{\delta}} \restriction Z = \rho \restriction Z$.

We show that $\tau^{\leq_{\delta}} \restriction W = \rho \restriction W$.
By (I3) and the definition of $\leq_{\delta}$,
for every $\alpha < \delta$,
$U_{\alpha}$ is an initial segment
of $\pair{I_{\delta}}{\leq'_{\delta}}$
and $\leq'_{\delta} \restriction U_{\alpha} =\kern6pt
\leq_{\alpha}$.
Since $I_{\delta}$ is an initial segment of
$\pair{U_{\delta}}{\leq_{\delta}}$
and $\leq_{\delta} \restriction I_{\delta} =\kern6pt
\leq'_{\delta}$, it follows that
$U_{\alpha}$ is an initial segment of
$\pair{U_{\delta}}{\leq_{\delta}}$
and $\leq_{\delta} \restriction U_{\alpha} =\kern6pt \leq_{\alpha}$
for every $\alpha < \delta$.
Since $U_{\alpha}$ is an initial segment of
$\pair{U_{\delta}}{\leq_{\delta}}$,
the order topology of
$\pair{U_{\alpha}}{\leq_{\delta} \restriction U_{\alpha}}$
is equal to the relative topology it inherits from
$\pair{U_{\delta}}{\tau^{\leq_{\delta}}}$.
And hence $\tau^{\leq_{\alpha}} = \tau^{\leq_{\delta}} \nrestriction U_{\alpha}$.
By the induction hypothesis,
$\tau^{\leq_{\alpha}} = \rho \nrestriction U_{\alpha}$.
So $\tau^{\leq_{\delta}} \nrestriction U_{\alpha} = \rho \nrestriction U_{\alpha}$.
Hence for every $\alpha < \delta$, $U_{\alpha}$ is compact in the
topology $\tau^{\leq_{\delta}}$.
$\setm{U_{\alpha}}{\alpha < \delta}$ is an increasing sequence of
initial segments of $\pair{U_{\delta}}{\leq_{\delta}}$
and $\sup^{\leq_{\delta}}(\bigcup_{\alpha < \delta} U_{\alpha}) = x_{\delta}$.
So $\bigcup_{\alpha < \delta} U_{\alpha} \cup \sngltn{x_{\delta}}$
is compact in $\pair{U_{\delta}}{\tau}$.
Recall that
$\bigcup_{\alpha < \delta} U_{\alpha} \cup \sngltn{x_{\delta}} = W$.
So $W$ is compact in $\pair{U_{\delta}}{\tau^{\leq_{\delta}}}$.
It is easy to see that
$\rfs{cl}^{\pair{X}{\rho}}(\bigcup_{\alpha < \delta} U_{\alpha}) = W$.
So $W$ is compact in $\pair{X}{\rho}$.

In order to show that $\rho \nrestriction W = \tau^{\leq_{\delta}} \nrestriction W$
it thus suffices to prove that
$\tau^{\leq_{\delta}} \nrestriction W \subseteq \rho \nrestriction W$.
Let $V$ be open in $\pair{W}{\tau^{\leq_{\delta}} \nrestriction W}$.
Suppose that $x_{\delta} \not\in V$.
Then $V = \bigcup_{\alpha < \delta}(V \cap U_{\alpha})$.
Take $\alpha < \delta$.
Then $V \cap U_{\alpha}$ is open in
$\pair{U_{\alpha}}{\tau^{\leq_{\delta}} \restriction U_{\alpha}}$.
Since $U_{\alpha}$ is an initial segment of $\pair{U_{\delta}}{\leq_{\delta}}$,
$\tau^{\leq_{\delta} \restriction U_{\alpha}} =
\tau^{\leq_{\delta}} \restriction U_{\alpha}$.
So $V \cap U_{\alpha}$ is open in
$\pair{U_{\alpha}}{\leq_{\delta} \restriction U_{\alpha}}$.
But $\leq_{\delta} \restriction U_{\alpha} =\kern6pt \leq_{\alpha}$.
So $V \cap U_{\alpha}$ is open in $\pair{U_{\alpha}}{\leq_{\alpha}}$.
By (I1), $V \cap U_{\alpha}$
is open in $\pair{U_{\alpha}}{\rho \restriction U_{\alpha}}$.
Since $U_{\alpha} \in \rho$, $V \cap U_{\alpha} \in \rho$. So $V \in \rho$.

Suppose next that $x_{\delta} \in V$.
Then $V - \sngltn{x_{\delta}}$ is open in $\pair{W}{\leq_{\delta} \restriction W}$.
By the previous paragraph,
$V - \sngltn{x_{\delta}}$ is open in $\pair{W}{\rho \restriction W}$.
It remains to show that $V$ contains a
$(\rho \restriction W)$-neighborhood of $x_{\delta}$.
Clearly, $V$ contains an open final segment of
$\pair{W}{\leq_{\delta} \restriction W}$.
Hence for some $\alpha < \delta$, $V \supseteq W - U_{\alpha}$.
But $W - U_{\alpha} = (U_{\delta} - U_{\alpha}) \cap W$.
Obviously, $(U_{\delta} - U_{\alpha}) \cap W$
is a $(\rho \restriction W)$-neighborhood of $x_{\delta}$.
So $V$ is open in $\pair{W}{\rho \restriction W}$.
This implies that 
$\tau^{\leq_{\delta}} \nrestriction W = \rho \nrestriction W$.

It follows that $\rho \restriction U_{\delta} = \tau^{\leq_{\delta}}$.

{\bf Case 1.2 } $\calD_{\delta} - \sngltn{I_{\delta}} \neq \emptyset$
and there is no $D \in \calD_{\delta} - \sngltn{I_{\delta}}$
such that $\rfs{cf}(D) \geq \aleph_1$.
For every $D \in \calD_{\delta} - \sngltn{I_{\delta}}$
let $\beta_D \in \dbltn{1}{\omega}$
and $\setm{x_{D,i}}{i < \beta_D}$ be a strictly increasing unbounded
sequence in $D$.
For every $D \in \calD_{\delta} - \sngltn{I_{\delta}}$
let $V_{D,0} = U_{x_{D,0}}$ 
and for $0 < i < \beta_D$
let $V_{D,i} = U_{x_{D,i}} - U_{x_{D,i - 1}}$.
Let $\gamma \leq \omega$ and $\setm{V_i}{i < \gamma}$
be a $\onetoone$ enumeration of
$\setm{V_{D,i}}{D \in \calD_{\delta} - \sngltn{I_{\delta}}
\mbox{ and } i < \beta_D}$.
Then $V_i$ is a scattered continuous image of a compact interval space
with rank $< \alpha_0$.
By the induction hypothesis there is a linear ordering
$\leq_i$ of $V_i$ such that $\tau^{\leq_i} = \rho \nrestriction V_i$.
Let $\leq_{\delta}$ be defined as follows.
\begin{itemize}
\item[(1)] $I_{\delta} \kern3pt\leq_{\delta}\kern3pt
x_{\delta} \leq \ldots
\kern3pt\leq_{\delta}\kern3pt V_n
\kern3pt\leq_{\delta}\kern3pt
\ldots \kern3pt\leq_{\delta}\kern3pt V_1 \kern3pt\leq_{\delta}\kern3pt
V_0$.
\item[(2)] $\leq_{\delta} \restriction\kern-3pt I_{\delta}
\kern3pt=\kern8pt \leq'_{\delta}$,
\ and \ $\leq_{\delta} \restriction \kern-2ptV_i
\kern6pt=\kern8pt \leq_i\ $
for every $i \in \omega$.
\end{itemize}
It is left to the reader to check that
$\tau^{\leq_{\delta}} = \rho \nrestriction U_{\delta}$.

{\bf Case 1.3 } $\calD_{\delta} = \sngltn{I_{\delta}}$.
Define $\leq_{\delta}$ as follows:
$I_{\delta} \leq_{\delta} x_{\delta}$
and
$\leq_{\delta} \restriction\kern-3pt I_{\delta} \kern3pt=\kern3pt 
\leq'_{\delta}$.
Note that in this case
$x_{\delta} = \max(\pair{U_{\delta}}{\leq_{\delta}})$.
So if $\delta = \lambda$,
then the second part of $(*)_{\alpha_0}$ is fulfilled.
It is left to the reader to check that
$\tau^{\leq_{\delta}} = \rho \nrestriction U_{\delta}$.

{\bf Case 2 }
Assume that $\pair{X - \sngltn{e}}{\leq_{\calU} \restriction (X - \sngltn{e})}$
contains an unbounded chain $I$ with uncountable cofinality,
and that $X - \sngltn{e} - \bigcup_{x \in I} U_x \neq \emptyset$.
Let $\lambda$ be an uncountable cardinal
and $\setm{x_{\alpha}}{\alpha < \lambda}$ be an unbounded strictly
increasing sequence in $X^{\hbox{\tiny$<$}_{\calU}\kern1.5pt e}$.
Let $D_0 = \bigcup_{\alpha < \lambda} U_{x_{\alpha}}$
and $W = D_0 \cup \sngltn{e}$.
Then $D_0 \in \calD_e$.
By Facts (1) and (2), there are two possiblities.

\begin{itemize}
\item[(1)]
There is $D_1 \in \calD_e - \sngltn{D_0}$
such that $\rfs{cf}(D_1) \geq \aleph_1$,
$\calD_e$ is finite and every member of $\calD_e - \dbltn{D_0}{D_1}$
has a maximum.

\item[(2)]
$\abs{\calD_e} \leq \aleph_0$
and for every $D \in \calD_e - \sngltn{D_0}$,
$\rfs{cf}(D) \leq \aleph_0$.
\end{itemize}

{\bf Case 2.1 } (1) happens.
Let $\sigma = \setm{\max(D)}{D \in \calD_e - \dbltn{D_0}{D_1}}$.
For\break
every $x \in \sigma$, $\rfs{rk}(U_x) < \alpha_0$.
So by the induction hypothesis $U_x$ is homeomorphic to
an interval space.
Since $U_x$ is clopen for every $x \in \sigma$,
it suffices to show
that $Z \eqdf X - \bigcup_{x \in \sigma} U_x$
is homeomorphic to an interval
space.
For $i = 0,1$ let $Z_i = D_i \cup \sngltn{e}$.
Then $Z_i$ fulfill the assumptions of Case 1.
Let $\leq_0$ be a linear ordering of $Z_0$
such that $e = \max(\pair{Z_0}{\leq_0})$
and $\tau^{\leq_0} = \rho \nrestriction Z_0$.
Let $\leq_1$ be a linear ordering of $Z_1$
such that $e = \min(\pair{Z_1}{\leq_1})$
and $\tau^{\leq_1} = \rho \nrestriction Z_1$.
Clearly, $Z = Z_0 \cup Z_1$.
Define the relation $\leq$ on $Z$ as follows:
\begin{itemize}
\item[(1)] $\leq \restriction\kern-3pt Z_0 \kern3pt=\kern6pt \leq_0$
and
$\leq \restriction\kern-3pt Z_0 \kern3pt=\kern6pt \leq_0$.
\item[(2)] $Z_0 \leq Z_1$.
\end{itemize}
It is left to the reader to check that $\leq$ is a linear ordering of
$Z_0 \cup Z_1$ and that 
$\tau^{\leq} = \rho \nrestriction Z$.

{\bf Case 2.2 } (2) happens.
This case is similar to Case 2.1.
For every $D \in \calD_e - \sngltn{D_0}$
let $\beta_D \in \dbltn{1}{\omega}$
and $\setm{x_{D,i}}{i < \beta_D}$ be a strictly increasing unbounded
sequence in $D$.
For every $D \in \calD_e - \sngltn{D_0}$
let $V_{D,0} = U_{x_{D,0}}$ 
and for $0 < i < \beta_D$
let $V_{D,i} = U_{x_{D,i}} - U_{x_{D,i - 1}}$.
Let $\gamma \leq \omega$ and $\setm{V_i}{i < \gamma}$
be a $\onetoone$ enumeration of
$\setm{V_{D,i}}{D \in \calD_e - \sngltn{D_0}
\mbox{ and } i < \beta_D}$.
Then $V_i$ is a scattered continuous image of a compact interval space
with rank $< \alpha_0$.
By the induction hypothesis there is a linear ordering
$\leq_i$ of $V_i$ such that $\tau^{\leq_i} = \rho \nrestriction V_i$.
Let $Z = D_0 \cup \sngltn{e}$. Then $Z$ fulfills the assumptions of
Case 1. So there is a linear ordering $\leq'$ of $Z$ such that
$e = \max(\pair{Z}{\leq'})$ and $\tau^{\leq'} = \rho \nrestriction Z$.
Let $\leq$ be defined as follows.
\begin{itemize}
\item[(1)] $D_0 \kern3pt\leq\kern3pt
e \leq \ldots
\kern3pt\leq\kern3pt V_n
\kern3pt\leq\kern3pt
\ldots \kern3pt\leq\kern3pt V_1 \kern3pt\leq\kern3pt
V_0$.
\item[(2)] $\leq \restriction\kern-3pt Z
\kern3pt=\kern8pt \leq'$,
\ and \ $\leq \restriction \kern-2ptV_i
\kern6pt=\kern8pt \leq_i\ $
for every $i \in \omega$.
\end{itemize}
It is left to the reader to check that
$\tau^{\leq} = \rho$.

{\bf Case 3 } 
Assume that $X - \sngltn{e}$
does not contains an unbounded chain with uncountable cofinality.
This case too is similar to Case 1.2.
For every $D \in \calD_e$
let $\beta_D \in \dbltn{1}{\omega}$
and $\setm{x_{D,i}}{i < \beta_D}$ be a strictly increasing unbounded
sequence in $D$.
For every $D \in \calD_e$
let $V_{D,0} = U_{x_{D,0}}$ 
and for $0 < i < \beta_D$
let $V_{D,i} = U_{x_{D,i}} - U_{x_{D,i - 1}}$.
Let $\gamma \leq \omega$ and $\setm{V_i}{i < \gamma}$
be a $\onetoone$ enumeration of
$\setm{V_{D,i}}{D \in \calD_e \mbox{ and } i < \beta_D}$.
Then $V_i$ is a scattered continuous image of a compact interval space
with rank $< \alpha_0$.
By the induction hypothesis there is a linear ordering
$\leq_i$ of $V_i$ such that $\tau^{\leq_i} = \rho \nrestriction V_i$.
Let $\leq$ be defined as follows.
\begin{itemize}
\item[(1)]
$e \leq \ldots
\kern3pt\leq\kern3pt V_n
\kern3pt\leq\kern3pt
\ldots \kern3pt\leq\kern3pt V_1 \kern3pt\leq\kern3pt
V_0$.
\item[(2)]
$\leq \restriction \kern-2ptV_i
\kern6pt=\kern8pt \leq_i\ $
for every $i \in \omega$.
\end{itemize}
It is left to the reader to check that
$\tau^{\leq} = \rho$.

\section{A lemma about CO spaces}
\label{s4}
\label{s-c}
\begin{definition}
\label{d3.1}
\label{c-dfn-1.6-24-01-04}
\begin{rm}
Let $K$ and $L$ be unitary scattered compact spaces.

(a) $K$ and $L$ are {\it almost homeomorphic} ($K \approx L$)
if there are clopen neighborhoods $U$ and $V$
of $e^K$ and $e^L$
respectively such that $U \cong V$.

(b) We define the relation
$K \prec^{\rm w} L$ as follows.
$K \prec^{\rm w} L$
if for some $K' \approx K$, $K' \subseteq L$,
$e^{K'} = e^{L}$ and $K \not\approx L$.

We also define the relation $K \prec L$.
Say that $K \prec L$
if for some $K' \approx K$, $K' \subseteq L$,
$\rfs{rk}(K) = \rfs{rk}(L)$ and $K \not\approx L$.
Note that this implies that $e^{K'} = e^L$.

(c) Let $X$ be a compact space and $D \subseteq X$.
For every $d \in D$ let $V_d$ be an open neighborhood of $d$.
The family $\calV \eqdf \setm{V_d}{d \in D}$ is called
a {\it strong Hausdorff system} for $D$
if for every distinct $d,e \in D$, $V_d \cap V_e = \emptyset$
and
$$\rfs{cl}(\bigcup\setm{V_d}{d \in D}) =
\bigcup \setm{\rfs{cl}(V_d)}{d \in D} \cup \rfs{cl}(D).$$
$\calV$ is called a {\it clopen strong Hausdorff system} for $D$
if every $V_d$ is clopen.
\end{rm}
\end{definition}

\begin{theorem}\label{bthm-24-05-04}\label{t3.2}
Let $X$ be a scattered compact space,
and assume that for every subset $S$ of $X$ with regular cardinality
there is $D \subseteq S$ such that $\abs{D} = \abs{S}$,
and $D$ has a clopen strong Hausdorff system.

\num{a} Suppose that there are unitary scattered compact spaces
$L$ and $M$ and a family $\setm{L_i}{i \in \omega}$
of subsets of $X$
such that $M \prec L$ and for every\break
$i < j <\omega$, $L_i \approx L$
and $e^{L_i} \neq e^{L_j}$.
Then $X$ is not a CO space.

\num{b} Suppose that there are unitary scattered compact spaces
$K$, $L$ and $M$ such that $M \prec L \prec^{\rm w} K \subseteq X$,
then $X$ is not a CO space.
\end{theorem}

\begin{definition}
\begin{rm}
Let $Y$ be a scattered compact space.

(a) For an ordinal $\theta$ define
$R_{\theta}(Y) \eqdf \setm{z \in Y}{\rfs{rk}^Y(z) = \theta}$.

(b) Let $K$ be a unitary space with rank~$\theta$.
We say that $Y$ is {\it $K$-based}\break
if $\rfs{rk}(Y) \geq \theta + 1$,
and there are $\calU,\calV \subseteq \rfs{Clop}(Y)$
such that the following holds.
\begin{itemize}
\item[(1)] $\calU$ is a pairwise disjoint family,
and $\calV$ is a pairwise disjoint family.
\item[(2)] For every $U \in \calU$,
$U \approx K$,
and for every $V \in \calV$,
$V$ is unitary and $\rfs{rk}(V) = \theta + 1$.
\item[(3)] $R_{\theta}(Y) \subseteq \bigcup \calU$
and $R_{\theta + 1}(Y) \subseteq \bigcup \calV$.
\end{itemize}

(c) Suppose that $K,L$ are unitary spaces with the same
rank~$\theta$.
We say that $Y$ is {\it $\dbltn{K}{L}$-based}
if $\rfs{rk}(Y) \geq \theta + 1$,
and there are $\calU,\calV \subseteq \rfs{Clop}(Y)$
such that the following holds.
\begin{itemize}
\item[(1)] $\calU$ is a pairwise disjoint family,
and $\calV$ is a pairwise disjoint family.
\item[(2)] For every $U \in \calU$,
$U \approx K$ or $U \approx L$.
For every $V \in \calV$,
$V$ is unitary and $\rfs{rk}(V) = \theta + 1$.
\item[(3)] $R_{\theta}(Y) \subseteq \bigcup \calU$
and $R_{\theta + 1}(Y) \subseteq \bigcup \calV$.
\item[(4)] For every $W \in \rfs{Clop}(Y)$,
if $\rfs{rk}(W) \geq \theta + 1$,
then there are $U,V \in \rfs{Clop}(Y)$
such that $U,V \subseteq W$, $U \approx K$ and $V \approx L$.
\end{itemize}

Note that a space $Y$ is $K$-based
iff it is $\dbltn{K}{K}$-based.
Suppose that $Y$ is $\dbltn{K}{L}$-based and $\calU$, $\calV$ are
families assured by the $\dbltn{K}{L}$-basedness of $Y$.
We denote $\calU,\calV$ by $\calU_Y$ and $\calV_Y$ respectively.
\hfill\myqed
\end{rm}
\end{definition}

The trivial proof of the following proposition is left to the
reader.

\begin{prop}\label{p3.4}\label{lemma-1.7}
\label{lemma-1}
Let $Y,Z,K,L$ be a compact scattered spaces,
and assume that $K$ and $L$ are unitary with the
same rank $\theta$.

\num{a}
Suppose that $Y$ is $\dbltn{K}{L}$-based
and $U \in \rfs{Clop}(Y)$. If $\rfs{rk}(U) \geq \theta + 1$,
then $U$ is $\dbltn{K}{L}$-based.

\num{b}
Assume that $K \not\approx L$
and that $Y$ is $K$-based
and $Z$ is $\dbltn{K}{L}$-based.
Then $Y \not\approx Z$.

\num{c} 
Assume that $K \not\approx L$
and that $Y$ is $K$-based
and $Z$ is $\dbltn{K}{L}$-based.
Assume further that $Y,Z$
are clopen unitary subspaces of $X$.
Then $e^Y \neq e^Z$.

\num{d}
Assume that $K \not\approx L$
and that $Y$ is $\dbltn{K}{L}$-based.
Let\break
$Z = Y - \bigcup \setm{U \in \calU_Y}{U \approx L}$.
Then $Z$ is $K$-based, $\rfs{rk}(Z) = \rfs{rk}(Y)$,
$e^Z = e^Y$,
$\calU_Z = \setm{U \in \calU_Y}{U \approx K}$
and $\calV_Z = \setm{V \cap Z}{V \in \calV_Y}$.
\end{prop}

\begin{prop}\label{p3.5}
Let $X$ be a scattered compact space,
and assume that for every infinite subset $S$ of $X$
there is an infinite subset $D \subseteq S$
such that $D$ has a clopen strong Hausdorff system.
Let $L$ be a unitary scattered compact space,
and $\setm{L_i}{i \in \omega}$ be a family
of subsets of $X$ such that
\begin{itemize}
\item[$(1)$]
$L_i \approx L$ for every $i \in \omega$
and for $i < j < \omega$, $e^{L_i} \neq e^{L_j}$,
\item[]
or
\item[$(2)$]
for every $\,i$, \,$L_i$ is unitary and $L$-based,
for every $i,j$, $\rfs{rk}(L_i) = \rfs{rk}(L_j)$,
and for $i < j < \omega$, $e^{L_i} \neq e^{L_j}$.
\end{itemize}
Then $X$ has a unitary $L$-based subset $F$ such that
$\rfs{rk}(F) = \rfs{rk}(L_i) + 1$.
\end{prop}

\noindent
{\bf Proof }
Let $X$, $L$ and $\setm{L_i}{i \in \omega}$ be as in the hypotheses
of the proposition.
Denote $\rfs{rk}(L)$ by $\theta$
and $\rfs{rk}(L_i)$ by $\alpha$.
Let $A = \setm{e^{L_i}}{i \in \omega}$.
For every $a \in A$, if $a = e^{L_i}$, denote $L_i$ by $L_a$.
Let $x$ be an isolated point of $\rfs{acc}(A)$
and let $U \in \rfs{Nbr}_{\srfs{clp}}^X(x)$ be such that
$U \cap \rfs{acc}(A) = \sngltn{x}$.
Let $B = A \cap U$.
Hence $\rfs{cl}(B) = B \cup \sngltn{x}$.
There are an infinite subset $C \subseteq B$
and a family $\calT = \setm{T_c}{c \in C}$
such that $\calT$ is a clopen strong Hausdorff system for $C$.
So $\rfs{cl}^X(\bigcup_{c \in C} T_c) = (\bigcup_{c \in C} T_c) \cup \sngltn{x}$.
For every $c \in C$ let $L'_c = L_c \cap U \cap T_c$,
and let $F = \rfs{cl}^X(\bigcup_{c \in C} L'_c)$.
Then
$F = (\bigcup_{c \in C} L'_c) \cup \sngltn{x}$.
Clearly, $R_{\alpha}(F) \cap (\bigcup_{c \in C} L'_c) = C$
and $\rfs{acc}(C) = \sngltn{x}$.
So $R_{\alpha + 1}(F) = \sngltn{x}$, $R_{\alpha}(F) = C$
and $\rfs{rk}(F) = \alpha + 1$.
Hence $F$ is unitary.
We now distinguish between the two cases.

{\bf Case 1 } For every $i \in \omega$, $L_i \approx L$.
\\
Hence $\alpha = \theta$. Let $\calU = \setm{L'_c}{c \in C}$.
Then $\calU \subseteq \rfs{Clop}(F)$,
$\calU$ is a pairwise disjoint family
and $R_{\theta}(F) = C \subseteq \bigcup \calU$.
Also, for every $c \in C$, $L'_c \approx L_c \approx L$.
We now distinguish between the two cases.
Define $\calU_F$ to be $\calU$ and $\calV_F$ to be $\sngltn{F}$.
Then $\calU_F$ and $\calV_F$ demonstrate that $F$ is $L$-based.

{\bf Case 2 } For every $i \in \omega$, $L_i$ is $L$-based.
\newline\noindent
Recall that $L'_c = L_c \cap U \cap T_c$. So $L'_c \in \rfs{Clop}(L_c)$.
So by Proposition~\ref{p3.4}(a), $L'_c$ is $L$-based.
Let $\calU = \bigcup \setm{\calU_{L'_c}}{c \in C}$
and $\calV = \bigcup \setm{\calV_{L'_c}}{c \in C}$.
Clearly, $\calU$ and $\calV$ are pairwise disjoint families.
(This is so because $\setm{L'_c}{c \in C}$ is a pairwise disjoint family).
Also, for every $U \in \calU$, $U \approx L$ and for every $V \in \calV$,
$V$ is unitary and $\rfs{rk}(V) = \theta + 1$.
Since $\rfs{rk}^F(x) = \alpha + 1 > \theta + 1$
and $F - \sngltn{x} = \bigcup_{c \in C} L'_c$,
it follows that
$R_{\theta}(F) = \bigcup_{c \in C} R_{\theta}(L'_c)$
and
$R_{\theta + 1}(F) = \bigcup_{c \in C} R_{\theta + 1}(L'_c)$.
Recall that $R_{\theta}(L'_c) \subseteq \bigcup \calU_{L'_c}$
and
$R_{\theta + 1}(L'_c) \subseteq \bigcup \calV_{L'_c}$.
Hence $R_{\theta}(F) \subseteq \bigcup \calU$
and
$R_{\theta + 1}(F) \subseteq \bigcup \calV$.
So $\calU$ and $\calV$ demonstrate that $F$ is $L$-based.
\hfill\proofend

\begin{prop}\label{p3.6}
Let $X$ be a compact CO space and $L \prec^{\rm w} K \subseteq X$.
Then there is a family $\setm{L_i}{i \in \omega}$ of subsets of $X$
such that for every $i$, $L_i \approx L$
and for every $i \neq j$, $e^{L_i} \neq e^{L_j}$.
\end{prop}

{\bf Proof }
We may assume that $K$ is clopen in $X$.
Let $L' \subseteq K$ be such that $L' \approx L$
and $e^{L'} = e^K$.
We define by induction clopen sets
$L_n$ and $L_{n,i}$, $i \leq n$.
We assume by induction on $n$ that for every $i \leq n$,
$L_{n,i} \cong L'$,
$L_n = \bigcup_{i \leq n} L_{n,i}$
and that for every $i \neq j$, $e^{L_{n,i}} \neq e^{L_{n,j}}$.

Let $L_0 = L_{0,0} \in \rfs{Clop}(X)$ be homeomorphic to $L'$.
Then the the induction hypotheses hold for $n = 0$.
Suppose that $L_n$ and $L_{n,i}$, \,$i \leq n$ have been defined.
Let $L_{n + 1}$ be a clopen set homeomorphic to $L' \cup L_n$
and $\iso{\psi}{L' \cup L_n}{L_{n + 1}}$.
For $i \leq n$ define $L_{n + 1,i} = \psi(L_{n,i})$
and let $L_{n + 1,n + 1} = \psi(L')$.
We check
that the induction hypotheses hold.
The only fact that needs to be verified is that for every distinct $i,j \leq n + 1$,
$e^{L_{n + 1,i}} \neq e^{L_{n + 1,j}}$.
If $i,j \leq n$ then $e^{L_{n + 1,i}} = \psi(e^{L_{n,i}})$
and
$e^{L_{n + 1,j}} = \psi(e^{L_{n,j}})$.
So since $e^{L_{n,i}} \neq e^{L_{n,i}}$ and $\psi$ is $\onetoone$,
it follows that $e^{L_{n + 1,i}} \neq e^{L_{n + 1,j}}$.
Suppose that $i \leq n$ and $j = n + 1$.
Note that $L_{n,i} \in \rfs{Nbr}^X_{\srfs{clp}}(e^{L_{n,i}})$.
However, there is no $U \in \rfs{Nbr}^X_{\srfs{clp}}(e^{L'})$
such that $U \cong L_{n,i}$.
Suppose by contradiction that such a $U$ exists.
Since $K \in \rfs{Nbr}^X_{\srfs{clp}}(e^{L'})$,
it follows that $K \approx K \cap U \approx U \approx L' \approx L$.
A contradiction. So $U$ does not exist, and hence $e^{L'} \neq e^{L_{n,i}}$.
It follows that
$$
e^{L_{n + 1,n +1}} = \psi(e^{L'}) \neq \psi(e^{L_{n,i}}) = e^{L_{n + 1,i}}.
$$
So the induction hypotheses hold for $n + 1$.
This completes the inductive construction.

It follows that $\setm{e^{L_{n,i}}}{n \in \omega, \ i \leq n}$ is infinite.
\hfill\proofend

\kern2mm

\noindent
{\bf Proof of Theorem \ref{bthm-24-05-04}}
(a) 
Assume by way of contradiction that $X$ is a CO space
and that $L$, $M$ and $\setm{L_i}{i \in \omega}$
are as in the hypotheses of~\ref{bthm-24-05-04}(a).
Denote $\rfs{rk}(L)$ by $\theta$.

We prove by induction on $\alpha \geq \theta + 1$
that there is $X_{\alpha} \subseteq X$ such that
$X_{\alpha}$ is $L$-based, $X_{\alpha}$ is unitary
and $\rfs{rk}(X_{\alpha}) = \alpha$.
By the first case of Proposition~\ref{p3.5},
there is a subspace $F \subseteq X$ such that $F$ is unitary of rank
$\theta + 1$ and $F$ is $L$-based. That is, $X_{\theta + 1}$ exists.
Suppose that $X_{\alpha}$ exists.
We may assume that $X_{\alpha}$ is clopen in $X$.
Denote $\calU_{X_{\alpha}}, \calV_{X_{\alpha}}$ by $\calU$ and $\calV$
respectively.
For every $x \in R_{\theta}(X_{\alpha})$ there is a unique
$L_x \in \calU$ such that $x \in L_x$.
So $L_x \approx L$.
Let $x \in R_{\theta + 1}(X_{\alpha})$. Then there is a unique
$V_x \in \calV$ such that $x \in V_x$.
Choose $A_x \subseteq R_{\theta}(V_x)$ such that
$A_x$ and $R_{\theta}(V_x) - A_x$ are infinite.
Note that
$\rfs{acc}(A_x) = \rfs{acc}(R_{\theta}(V_x) - A_x) = \sngltn{x}$.
For every $y \in A_x$ choose $M_y \subseteq L_y$ such that
$M_y \approx M$.
Define 
\medskip

$\begin{array}{rcl}
\calU\nprime &=& \setm{M_y}{x \in R_{\theta + 1}(X_{\alpha}) \mbox{ and }
y \in A_x} \\
\rule{0pt}{16pt}
&&\cup \setm{L_y}{x \in R_{\theta + 1}(X_{\alpha}) \mbox{ and }
y \in R_{\theta}(V_x) - A_x} \\
\rule{0pt}{16pt}
&&\cup
\setm{L_y}{ y \in R_{\theta}(X_{\alpha}) - \bigcup \calV}.
\end{array}$
\medskip

\noindent
Let $Y = \rfs{cl}(\bigcup \calU\nprime)$ and
$\calV\nprime = \setm{V \cap Y}{V \in \calV}$.

We shall see that $Y$ is $\dbltn{L}{M}$-based,
and that $\calsubU_Y$ and $\calsubV_Y$ can be taken to be
$\calU\nprime$ and $\calV\nprime$.
Let $y \in R_{\theta}(X_{\alpha})$. Define $N_y = M_y$
if for some $x \in R_{\theta + 1}(X_{\alpha})$, $y \in A_x$,
and otherwise let $N_y = L_y$.
Then for every $y \in R_{\theta}(X_{\alpha})$, $y \in N_y \subseteq Y$
and $N_y \in \rfs{Clop}(Y)$.
So $R_{\theta}(X_{\alpha}) \subseteq Y$
and for every  $y \in R_{\theta}(X_{\alpha})$,
$\rfs{rk}^Y(y) = \rfs{rk}^{N_y}(y) = \theta$.
That is, $R_{\theta}(X_{\alpha}) \subseteq R_{\theta}(Y)$.
Suppose by way of contradiction that
$R_{\theta}(Y) - R_{\theta}(X_{\alpha}) \neq \emptyset$,
and let $y \in R_{\theta}(Y) - R_{\theta}(X_{\alpha})$.
Since $\rfs{rk}^Y(y) \leq \rfs{rk}^{X_{\alpha}}(y)$,
$\rfs{rk}^{X_{\alpha}}(y) > \theta$.
Hence every neighborhood of $y$ intersects  $R_{\theta}(X_{\alpha})$.
Since $R_{\theta}(X_{\alpha}) \subseteq R_{\theta}(Y)$,
every neighborhood of $y$ intersects $R_{\theta}(Y)$.
This contradicts the fact that $\rfs{rk}^Y(y)  = \theta$.
So $R_{\theta}(X_{\alpha}) = R_{\theta}(Y)$.
Clearly, $\rfs{cl}(R_{\theta}(X_{\alpha})) = D_{\theta}(X_{\alpha})$ and
$\rfs{cl}(R_{\theta}(Y)) = D_{\theta}(Y)$.
So $D_{\theta}(X_{\alpha}) = D_{\theta}(Y)$.
So for every $\eta \geq \theta$, $D_{\eta}(X_{\alpha}) = D_{\eta}(Y)$.
It follows that $\rfs{rk}(Y) = \rfs{rk}(X_{\alpha})$, $Y$ is unitary
and $e^Y = e^{X_{\alpha}}$.

We next show that one can take $\calsubU_Y$ to be $\calU\nprime$.
Clearly, $\calU\nprime$ is a pairwise disjoint family, every member of
$\calU\nprime$ is almost homeomorphic to either $L$ or $M$.
For every space $Z$,
$R_{\theta + 1}(Z) = D_{\theta + 1}(Z) - D_{\theta + 2}(Z)$.
So since
$D_{\eta}(X_{\alpha}) = D_{\eta}(Y)$ for every $\eta \geq \theta$,
$R_{\theta + 1}(X_{\alpha}) = R_{\theta + 1}(Y)$.
By the construction, for every $x \in R_{\theta + 1}(X_{\alpha})$
and a neighborhood $W$ of $x$ there are $U,V \in \rfs{Clop}(Y)$
such that $U,V \subseteq W$, $U \approx L$ and $V \approx M$.
So for every $x \in R_{\theta + 1}(Y)$
and a neighborhood $W$ of $x$ there are $U,V \in \rfs{Clop}(Y)$
such that $U,V \subseteq W$, $U \approx L$ and $V \approx M$.
By the construction, $R_{\theta}(X_{\alpha}) \subseteq \bigcup \calU\nprime$.
So since\break
$R_{\theta}(Y) = R_{\theta}(X_{\alpha})$,
it follows that $R_{\theta}(Y)\subseteq \bigcup \calU\nprime$.

We check that $\calsubV_Y$ can be taken to be $\calV\nprime$.
Clearly, $\calV\nprime$ is a pairwise disjoint family.
Since $R_{\theta + 1}(X_{\alpha}) \subseteq \bigcup \calV$,
it follows that
$R_{\theta + 1}(X_{\alpha}) \cap Y \subseteq
(\bigcup \calV) \cap Y$.\break
But
$R_{\theta + 1}(X_{\alpha}) = R_{\theta + 1}(Y) =
R_{\theta + 1}(Y) \cap Y$
and $(\bigcup \calV) \cap Y = \bigcup \calV\nprime$.
So\break
$R_{\theta + 1}(Y) \subseteq \bigcup \calV\nprime$.
Let $V \in \calV$. Denote $e^V$ by $x$. So $V = V_x$.
Also, $x \in R_{\theta + 1}(X_{\alpha}) = R_{\theta + 1}(Y)$.
So $x \in V \cap Y$ and indeed
$(V \cap Y) \cap R_{\theta + 1}(Y) = \sngltn{x}$.
So $V \cap Y$ is unitary and $\rfs{rk}^Y(V) = \theta + 1$.
We have shown everything that is required in order to conclude that
$Y$ is $\dbltn{L}{M}$-based.

We verify that $Y \prec X_{\alpha}$.
We have already seen that $\rfs{rk}(Y) = \rfs{rk}(X_{\alpha})$
and that $Y$ is unitary. Also, $Y \subseteq X_{\alpha}$ and $X_{\alpha}$ is unitary.
It remains to show that $Y \not\approx X_{\alpha}$.
Recall that $X_{\alpha}$ is $L$-based and that $Y$ is $\dbltn{L}{M}$-based.
Also, $M \prec L$ and hence $M \not\approx L$.
Then by Proposition~\ref{p3.4}(b), $Y \not\approx X_{\alpha}$.
So $Y \prec X_{\alpha}$.

It follows that $Y \prec^{\srfs{w}} X_{\alpha}$.
By Proposition~\ref{p3.6},
there is a family $\setm{Y_i}{i \in \omega}$ of subsets of $X$
such that for every $i$, $Y_i \approx Y$
and for every $i \neq j$, $e^{Y_i} \neq e^{Y_j}$.

For every $i \in \omega$
let $Z_i = Y_i - \bigcup \setm{U \in \calsubU_{Y_i}}{U \approx M}$.
By Proposition~\ref{p3.4}(d), $Z_i$ is $L$-based, $e^{Z_i} = e^{Y_i}$
and $\rfs{rk}(Z_i) = \alpha$.
By the second case of Proposition~\ref{p3.5},
there is $Z \subseteq X$ such that
$\rfs{rk}(Z) = \alpha + 1$ and $Z$ is $L$-based.
So $X_{\alpha + 1} \eqdf Z$ is as required.

\kern2mm

Let $\delta$ be a limit ordinal,
and suppose that for every $\alpha < \delta$, $X_{\alpha}$ exists.
Set $\theta = \rfs{rk}(L)$.
Let $\lambda = \rfs{cf}(\delta)$
and $\setm{\alpha_i}{i < \lambda}$ be a strictly increasing sequence
converging to $\delta$ such that $\alpha_0 > \theta + 1$.
For $i < \lambda$ let $Y_i \subseteq X$ be a clopen unitary $L$-based set with rank
$\alpha_i$.
Hence for every $i < \lambda$, $\rfs{rk}^X(e^{Y_i}) = \alpha_i$.
It follows that $e^{Y_i} \neq e^{Y_j}$ for every $i \neq j$.
So $\abs{\setm{e^{Y_i}}{i < \lambda}} = \lambda$.
Hence there are $A \subseteq \setm{e^{Y_i}}{i < \lambda}$ and
$\calW = \setm{W_a}{a \in A}$
such that $\abs{\calA} = \lambda$
and $\calW$ is a clopen strong Hausdorff system for $A$.
That is, $\calW$ is a pairwise disjoint family consisting of clopen sets,
$a \in W_a$ for every $a \in A$,
and $\rfs{cl}^X(\bigcup \calW) = (\bigcup \calW) \cup \rfs{acc}^X(A)$.
For $a = e^{Y_i} \in A$ set $Y^a = Y_i$.

Note the following fact.
$(*)$ If $\calV = \setm{V_b}{b \in B}$ is a clopen strong Hausdorff system for $B$,
and $\calF = \setm{F_b}{b \in B}$ is a family of closed sets such that
$b \in F_b \subseteq V_b$ for every $b \in B$.
Then $\rfs{cl}^X(\bigcup \calF) = (\bigcup \calF) \cup \rfs{acc}^X(B)$.

For every $a \in A$ let $Y_0^a = Y^a \cap W_a$,
and let $\whatY = \rfs{cl}^X(\bigcup \setm{Y_0^a}{a \in A})$.
So by $(*)$, $\whatY = (\bigcup \setm{Y_0^a}{a \in A}) \cup \rfs{acc}^X(A)$.
Note that for every $a \in A$, $Y_0^a$ is clopen in $\whatY$,
$Y_0^a \approx Y^a$ and $a = e^{Y_0^a}$.
Hence $\rfs{rk}^{\whatY}(a) = \rfs{rk}^{Y_0^a}(a) = \rfs{rk}^{Y^a}(a)$.
It follows that $\sup_{a \in A} \rfs{rk}^{\whatY}(a) \geq \delta$.
So $\rfs{rk}(\whatY) \geq \delta$.

{\bf Claim 1 } For every $y \in \whatY - \bigcup_{a \in A} Y_0^a$,
$\rfs{rk}^{\whatY}(y) > \theta + 1$.

{\bf Proof }
Let $y \in \whatY - \bigcup_{a \in A} Y_0^a$.
Recall that
$\whatY = (\bigcup_{a \in A} Y_0^a) \cup \rfs{acc}(A)$.
So $y \in \rfs{acc}^X(A)$.
We show that for every $V \in \rfs{Nbr}^X(y)$
there is $z \in V \cap \whatY$ such that $\rfs{rk}^{\whatY}(z) > \theta + 1$.
We may assume that $V$ is clopen.
Since $y \in \rfs{acc}^X(A)$, it follows that $V \cap A \neq \emptyset$.
Let $a \in V \cap A$.
Then
$\rfs{rk}^{Y_0^a \cap V}(a) = \rfs{rk}^{Y_0^a}(a) = \rfs{rk}^{Y^a}(a) > \theta + 1$.
Clearly, $Y_0^a \cap V \subseteq \whatY$.
So $\rfs{rk}^{\whatY}(a) > \theta + 1$.
It follows that $\rfs{rk}^{\whatY}(y) > \theta + 1$.
So Claim 1 is proved.

{\bf Claim 2 } $\whatY$ is $L$-based.

{\bf Proof }
Recall that for every $a \in A$, $Y_0^a = Y^a \cap W_a$
and $\rfs{rk}(Y_0^a) = \rfs{rk}(Y^a) > \theta + 1$.
Hence by Proposition~\ref{p3.4}(a), $Y_0^a$ is $L$-based.
Let $\calsubU_a,\calsubV_a$ demonstrate that $Y_0^a$ is $L$-based.
Set $\calU = \bigcup_{a \in A} \calsubU_a$
and $\calV = \bigcup_{a \in A} \calsubV_a$.
We show that $\calU,\calV$ demonstrate that $\whatY$ is $L$-based.
Since $\setm{Y_0^a}{a \in A}$ is a pairwise disjoint family,
and for every $a \in A$, $\calsubU_a$ is a pairwise disjoint family,
it follows that $\calU$ is a pairwise disjoint family.
Similarly, $\calV$ is a pairwise disjoint family.
It is also trivial that for every $U \in \calU$, $U \approx L$
and that for every $V \in \calV$, $\rfs{rk}(V) = \theta + 1$.

Let $y \in R_{\theta}(\whatY)$.
By Claim 1, there is $a \in A$ such that $y \in Y_0^a$.
Since $Y_0^a$ is clopen in $\whatY$, $\rfs{rk}^{Y_0^a} = \rfs{rk}^{\whatY} = \theta$.
So there is $U \in \calsubU_a$ such that $y \in U$. But $U \in \calU$.
It follows that $y \in \bigcup \calU$.
That is, $R_{\theta}(\whatY) \subseteq \bigcup \calU$.

An identical argument shows that
$R_{\theta + 1}(\whatY) \subseteq \bigcup \calV$.
We have shown that $\whatY$ is $L$-based.
So Claim 2 is proved.

Let $x \in \whatY$ be such that $\rfs{rk}^{\whatY}(x) = \delta$,
and let $T \in \rfs{Clop}(X)$
be such that $T \cap D_{\delta}(\whatY) = \sngltn{x}$.
Set $Y = T \cap \whatY$.
Hence $Y$ is unitary of rank $\delta$.
By Proposition~\ref{p3.4}(a), $Y$ is $L$-based.
Define $X_{\delta} = Y$. Then $X_{\delta}$ is as required.

We have proved that for every ordinal $\alpha$,
$X$ contains a subset with rank~$\alpha$. A contradiction.
So $X$ is not a CO-space.

(b) Let $M \prec L \prec^{\rm w} K \subseteq X$ be as in Part (b).
By Proposition \ref{p3.6}, there is a family
$\setm{L_i}{i \in \omega}$ such for every $i \in \omega$,
$L_i \cong L$ and for every $i \neq j$, $e^{L_i} \neq e^{L_j}$.
By Part (a) of this theorem, and since $M \prec L$, $X$ is not
a CO space.
\hfill\proofend

\section{CO spaces must omit the obstructions}
\label{s5}

The existence of strong Hausdorff systems is used in this section.
However, the full strength of Lemma \ref{l2.3} is not needed
and only the following fact is used.

\begin{cor}\label{c5.1}
Let $X$ be a continuous image of a compact interval space.
Let $A \subseteq X$ be such that $\abs{A}$ is regular
and $\rfs{cl}(A)$ is scattered.
Then there is $B \subseteq A$ such that $\abs{B} = \abs{A}$
and $B$ has a strong Hausdorff system.
\end{cor}

\begin{prop}\label{p5.1}
Let $X$ be a scattered continuous image of an interval space,
and assume that $X$ is a CO space.

\num{a}
If $\trpl{\kappa}{\lambda}{\mu}$ is a legal triple,
then $X_{\kappa,\lambda,\mu}$ is not embeddable in $X$.

\num{b}
The set $\setm{e^F}{F \subseteq X \mbox{ and } F \cong X_{\aleph_1}}$
is finite.
%
\end{prop}

\noindent
{\bf Proof }
(a) Suppose by way of contradiction that
$\trpl{\kappa}{\lambda}{\mu}$ is a legal triple and $X_{\kappa,\lambda,\mu}$
is embeddable in $X$.
We may assume that $\kappa \leq \lambda \leq \mu$.
It is then obvious that
$\mu + 1 \prec \mu + 1 + \lambda^* \prec X_{\kappa,\lambda,\mu}$.
So by Theorem \ref{t3.2}(b) and Corollary~\ref{c5.1},
$X$ is not a CO space. A contradiction,
so $X_{\kappa,\lambda,\mu}$ is not embeddable in $X$.

(b) Assume by contradiction that
$\setm{e^F}{F \subseteq X \mbox{ and } F \cong X_{\aleph_1}}$
is infinite.
Clearly,
$X_{\aleph_0} \prec^{\rm w} X_{\aleph_1}$,
so by Theorem~\ref{t3.2}(a) and Corollary~\ref{c5.1},
$X$ is not a CO space. A contradiction, so
$\setm{e^F}{F \subseteq X \mbox{ and } F \cong X_{\aleph_1}}$ is finite.
%
\hfill\proofend

\kern2mm

We also have to prove that obstructions of the type $X_{\lambda,\vecmu}$
are not embeddable in $X$.
In order to show this we consider the following space.
Let $\lambda$ be a cardinal and $S \subseteq \lambda$.
For $\alpha \in S$ let $L_{\alpha} = 1 + \omega^*$
and for every $\alpha \in \lambda - S$ let $L_{\alpha} = 1$.
Define $X_{\lambda,S}$ to be the topological space with universe
$\sum_{\alpha < \lambda} L_{\alpha} + 1$ and with the order topology
as its topology.
\index{N@$X_{\lambda,S}$}

\begin{prop}\label{p5.2}
\num{a}
Let $\lambda$ be an uncontable regular cardinal
and $S_1,S_2 \subseteq \lambda$ be subsets of $\lambda$
such that $S_1 - S_2$ is stationary.
Then $X_{\lambda,S_1} \not\cong X_{\lambda,S_2}$.
Also, $X_{\lambda,S_1}$ is not homeomorphic to an ordinal.

\num{b}
Let $X$ be a scattered continuous image of a compact interval space,
and assume that $X$ is a CO space.
Let $Y = X_{\lambda,\vecmu}$ and assume that $\rfs{Dom}(\vecmu)$
is a stationary set in $\lambda$. Then $Y$ is not embeddable in $X$.
\end{prop}

\noindent
{\bf Proof }
(a) Suppose by way of contradiction that
$\iso{f}{X_{\lambda,S_1}}{X_{\lambda,S_2}}$.
For $S \subseteq \lambda$ represent $X_{\lambda,S}$ in the following way.
$X_{\lambda,S} = (\lambda + 1) \cup (S \times \omega)$,
where for $\beta \in S$,
$L_{\beta} = \sngltn{\beta} \cup (\sngltn{\beta} \times \omega)$
and for $\beta \not\in S$,
$L_{\beta} = \sngltn{\beta}$.
For $\alpha < \lambda$ denote
$X_{\lambda,S} \restriction [\alpha,\lambda] = 
[\alpha,\lambda] \cup ((S \cap [\alpha,\lambda]) \times \omega)$.

Note that
$\rfs{acc}(\lambda + 1) \subseteq \rfs{acc}(X_{\lambda,S_i}) \subseteq
\lambda + 1$.
This implies that
$f[\rfs{acc}(\lambda)] \subseteq \lambda + 1$.
It follows that there is a club $C$ in $\lambda$
such that $f \nrestriction C = \rfs{Id}$.
For every $\beta \in S_1$,
$\beta \in \rfs{acc}^{X_{\lambda,S_1}}(\sngltn{\beta} \times \omega)$
and
for every $\beta \in S_2$,
$\beta \not\in
\rfs{acc}^{X_{\lambda,S_2}}(X_{\lambda,S_2} \restriction [\beta + 1,\lambda])$.
Hence for every $\beta \in (S_1 - S_2) \cap C$,
there are $\gamma_{\beta} < \beta$ and $m_{\beta},n_{\beta} \in \omega$
such that $f(\pair{\beta}{m_{\beta}}) = \gamma_{\beta}$
or $f(\pair{\beta}{m_{\beta}}) = \pair{\gamma_{\beta}}{n_{\beta}}$.
There is a stationary subset $S \subseteq (S_1 - S_2) \cap C$
such that for every $\alpha,\beta \in S$,
$m_{\alpha} = m_{\beta}$ and $n_{\alpha} = n_{\beta}$.
By Fodor's Lemma, for some distinct $\alpha,\beta \in S$,
$\gamma_{\alpha} = \gamma_{\beta}$.
So $f(\pair{\beta}{n_{\beta}}) = f(\pair{\gamma}{n_{\gamma}})$.
So $f$ is not $\onetoone$. A contradiction, so $f$ does not exist.

We show that $X_{\lambda,S_1}$ is not homeomorphic to an ordinal space.
Since $\lambda + 1 = X_{\lambda,\emptyset}$,
$X_{\lambda,S_1} \not\cong \lambda + 1$.
But $X_{\lambda,S_1}$ is a unitary space with rank $\lambda$,
and upto a homeomorphism,
the only ordinal space which is unitary with rank $\lambda$
is $\lambda + 1$.
So $X_{\lambda,S_1}$ is not homeomorphic to an ordinal space.

(b) Suppose by way of contradiction that
$Y = X_{\lambda,\vecmu} \subseteq X$.
Let\break
$S = \rfs{Dom}(\vecmu)$.
Let $Y_0$ be a clopen unitary subspace of $Y$
such that\break
$e^{Y_0} = \lambda + 1$.
Then $Y_0$ contains a clopen subspace $Y_1$ homeomorphic to
$X_{\lambda,\vecmu_1}$, where $\rfs{Dom}(\vecmu_1)$
is a final segment of $S$.
Clearly $e^{Y_1} = \lambda$.
Now, replace $Y$ by $Y_1$. Then we may assume that 
$e^{Y} = \lambda$.

By the easy direction of Theorem \ref{t2.1}, $X_{\lambda,\vecmu}$
is not homeomorphic to an interval space.
Moreover, if $U \in \rfs{Nbr}^Y_{\srfs{clp}}(\lambda)$,
then $U$ contains a space homeomorphic to $X_{\lambda,\vecmu'}$,
where $\rfs{Dom}(\vecmu')$ is a final segment of $S$.
So $U$ is not homeomorphic to an interval space.
Let $Z$ be the subspace of $Y$ whose universe is
$(\lambda + 1) \cup
(S \times \hatomega)$.
Then $Z$ is homeomorphic to $X_{\lambda,S}$ which is an interval space.
So $Y \not\approx Z$.
Also, $e^Z = \lambda$.
It follows that $Z \prec^{\rm w} Y$.

Clearly, $\lambda + 1$ is a closed unitary subspace of $Z$,
$e^{\lambda + 1} = \lambda = e^Z$ and
$\rfs{rk}^Z(\lambda) = \rfs{rk}^{\lambda + 1}(\lambda)$.
For every $U \in \rfs{Nbr}^Z_{\srfs{clp}}(\lambda)$
there is a final segment $S'$ of $S$ such that $U$ contains a subspace
homeomorphic to $X_{\lambda,S'}$.
By Part (a), $U$ is not homeomorphic to an ordinal.
So $Z \not \approx \lambda + 1$.
It follows that $\lambda + 1 \prec Z$.
We thus have $\lambda + 1 \prec Z \prec^{\rm w} Y$.
By Corollary~\ref{c5.1} and Theorem~\ref{t3.2}(b),
$X$ is not a CO space.
A contradiction so $X_{\lambda,\vecmu}$ is not embeddable in $X$.
\rule{0pt}{0pt}\hfill\proofend

\section{The characterization}
\label{s6}

In this section we prove Theorem \ref{t1.1}.
\begin{theorem}\label{t6.1}
Let $X$ be a scattered continuous image of a compact interval space,
and assume that $X$ is a CO space.
Then there is
a finite family of pairwise disjoint spaces
$\setm{Y_i}{i \in I}$
and an ordinal $\alpha + 1$ disjoint from the $Y_i$'s
such that
$X \cong (\alpha + 1) \cup \bigcup_{i \in I} Y_i$ and
\begin{itemize}
\item[\num{1}] For every $i \in I$ either
$Y_i \cong X_{\lambda,\mu}$, where $\lambda,\mu$ are infinite
regular cardinals and  $\mu > \aleph_0$ or
$Y_i \cong X_{\aleph_1}$.
\item[\num{2}] $\alpha \geq \alpha(Y_i)$ for every $i \in I$.
\end{itemize}
\end{theorem}

We quote the following Theorem from \cite{BBR}
\begin{theorem}\label{t6.2}
Let $X$ be a compact interval space,
and assume that $X$ is a CO space.
Then there is
a finite family of pairwise disjoint spaces
$\setm{Y_i}{i \in I}$
and an ordinal $\alpha + 1$ disjoint from the $Y_i$'s
such that
$X \cong (\alpha + 1) \cup \bigcup_{i \in I} Y_i$ and
\begin{itemize}
\item[\num{1}] For every $i \in I$,
$Y_i \cong X_{\lambda,\mu}$, where $\lambda,\mu$ are infinite
regular cardinals and  $\mu > \aleph_0$.
\item[\num{2}] $\alpha \geq \alpha(Y_i)$ for every $i \in I$.
\end{itemize}
\end{theorem}

\noindent
{\bf Proof of Theorem \ref{t6.1} }
Let
$\sigma = \setm{e^F}{F \subseteq X \mbox{ and } F \cong X_{\aleph_1}}$.
Then by Proposition \ref{p5.1}(b), $\sigma$ is finite.
For every $x \in \sigma$ let $F_x \subseteq X$ be such that
$F_x \cong X_{\aleph_1}$ and $x = e^{F_x}$.
Let $\whatF = \bigcup_{x \in \sigma} F_x$,
$Z$ be a clopen subset of $X$ homeomorphic to $\whatF$
and $\psi$ be a homeomorphism between $\whatF$ and $Z$.
Clearly, $\psi[\sigma] \subseteq \sigma$ and 
$\abs{\psi[\sigma]} = \abs{\sigma}$.
So $\psi[\sigma] = \sigma$. That is, $\sigma \subseteq Z$.

Let $Y = X - Z$.
Since $Y \cap \sigma = \emptyset$, it follows that $Y$ does not contain a subspace
homeomorphic to $X_{\aleph_1}$.
Since $X$ is a CO space
and by Propositions~\ref{p5.1}(a) and \ref{p5.2}(b),
$X$ does not contain a subspace homeomorphic to
$X_{\kappa,\lambda,\mu}$, where $\trpl{\kappa}{\lambda}{\mu}$
is a legal triple, and 
$X$ does not contain a subspace homeomorphic to
$X_{\lambda,\vecmu}$, where $\rfs{Dom}(\vecmu)$ is stationary in
$\lambda$.
So the same holds for $Y$.
By Theorem \ref{t2.1}, $Y$ is homeomorphic to an interval space.
We claim that $Y$ is a CO space.
If $Y$ is countable, then $Y$ is homeomorphic to an ordinal.
So $Y$ is CO.
Assume that $Y$ is uncountable.
Let $F$ be a closed subset of $Y$. There is $U \in \rfs{Clop}(X)$
such that $U \cong F$.
Let $V = U \cap Z$ and $W = U \cap Y$.
If $V = \emptyset$, then $U \in \rfs{Clop}(Y)$,
so there is nothing more to do.
Suppose that $V \neq \emptyset$.
Recall that $V \subseteq U \cong F \subseteq Y$. 
So $V$ is a closed subspace of an interval space,
and hence $V$ too is a compact interval space.
The only compact interval spaces embeddable in $Z$ are finite spaces
and spaces which are a disjoint union
of finitely many copies of $X_{\aleph_0}$.
So for some $n \in \omega$, $V \cong \omega \mcdot n + 1$,
or $V$ is finite.
Since $Y$ is uncountable and $Y$ is a scattered compact interval space,
it contains a clopen set homeomorphic to
$\omega^2 + 1$. So for every $n \in \omega$,
$Y \cong Y \cup (\omega \mcdot n + 1)$, where the union is disjoint.
It thus suffices to find a clopen subset of
$Y \cup (\omega \mcdot n + 1)$
which is homeomorphic to $U = V \cup W$.
$W$ is a clopen subset of $Y$ and either $V$ is finite or it is
homeomorphic to $\omega \mcdot n + 1$.
In either case $V$ is homeomorphic to a clopen subset of
$\omega \mcdot n + 1$.
So $Y$ is a CO space.

By Theorem \ref{t6.2},
there is a finite family of pairwise disjoint spaces\break
$\setm{Y_i}{i \in I}$
and an ordinal $\alpha + 1$ disjoint from the $Y_i$'s
such that:\break
(1) $Y \cong (\alpha + 1) \cup \bigcup_{i \in I} Y_i$;
(2) for every $i \in I$,
$Y_i \cong X_{\lambda,\mu}$, where $\lambda,\mu$ are infinite
regular cardinals and  $\mu > \aleph_0$;
(3) $\alpha \geq \alpha(Y_i)$ for every $i \in I$.

If $\sigma = \emptyset$, then the above description of $Y$
fulfills the requirements of the theorem.
Suppose that $\sigma \neq \emptyset$.
Then it remains to show that $\alpha \geq \omega^2$.
This is certainly true if $I \neq \emptyset$.
So suppose that $I = \emptyset$.
Note that $X_{\aleph_0} \prec X_{\aleph_1}$.
So by Proposition \ref{p3.6},
there is a family $\setm{L_i}{i \in \omega}$ of subsets of $X$
such that for every $i$, $L_i \approx X_{\aleph_0}$,
and for every $i \neq j$, $e^{L_i} \neq e^{L_j}$.
$Z \cap \setm{e^{L_i}}{i \in \omega}$ is finite.
So $(\alpha + 1) \cap \setm{e^{L_i}}{i \in \omega}$ is infinite.
That is, $R_1(\alpha + 1)$ is infinite. This implies that
$\alpha \geq \omega^2$.
\smallskip\hfill\proofend

\noindent
{\bf Proof of Theorem~\ref{t1.1} }
Combine Theorems~\ref{t1.2} and \ref{t6.1}.
\hfill\proofend

\section{Characterization of CO compact interval spaces.}
\label{s7}
In the previous section we quoted without proof Theorem~\ref{t6.2} from~\cite{BBR}.
However, Theorem~\ref{t6.2} follows easily from the previous sections.
So for completeness, we include a proof of Theorem~\ref{t6.2}.

The following proposition is an addition to Theorem~\ref{t3.2}(a).
\begin{prop}\label{p7.1}
Let $X$ be a scattered compact space,
and assume that for every subset $S$ of $X$ with regular cardinality
there is $D \subseteq S$ such that $\abs{D} = \abs{S}$,
and $D$ has a clopen strong Hausdorff system.

Suppose that there is a family $\setm{L_i}{i \in \omega}$
of compact subsets of $X$ such that
$(1)$ for every $i \in \omega$, $L_i$ is unitary,
and $(2)$ for every distinct $i,j \in \omega$,
$U \in \rfs{Nbr}^{L_i}_{\srfs{clp}}(e^{L_i})$
and $V \in \rfs{Clop}(L_j)$,
$U \not\cong V$.
\underline{Then} $X$ is not a CO space.
\end{prop}

\noindent
{\bf Proof }
Suppose by contradiction that $X$ is a CO space.
Denote $e^{L_i}$ by $e_i$.
We may assume that for every $i \in \omega$, $L_i$ is clopen in $X$.
Then by (2), for every distinct $i,j \in \omega$,
$e_i \neq e_j$.
We may further assume that $\rfs{acc}(\setm{e_i}{i \in \omega})$ is a singleton.
Denote it by $x$.
We may also assume that for every $i < j \in \omega$,
$\rfs{rk}^X(e_i) \leq \rfs{rk}^X(e_j)$.
There is an infinite subset $\sigma \subseteq \omega$ such that
$\setm{e_i}{i \in \sigma}$ has a clopen strong Hausdorff system.
We may assume $\sigma = \omega$.
Let $\setm{U_i}{i \in \omega}$ be a clopen strong Hausdorff system for
$\setm{e_i}{i \in \omega}$.
Let $L'_i = L_i \cap U_i$
and
$K = (\bigcup_{i \in \omega} L'_i) \cup \sngltn{x}$.
Then $K$ is closed and unitary, $e^K = x$
and $\rfs{rk}^K(x) = \rfs{Sup}(\setm{\rfs{rk}^X(e_i)}{i \in \omega})$.
We define $M,L$ such that $M \subseteq L \subseteq K$.
Let $\tau \subseteq \sigma \subseteq \omega$
be such that $\tau$, $\sigma - \tau$ and $\omega - \sigma$ are infinite.
Let
$L = (\bigcup_{i \in \sigma} L'_i) \cup \sngltn{x}$
and $M = (\bigcup_{i \in \tau} L'_i) \cup \sngltn{x}$.
It is obvious that $M,L$ are closed and unitary,
that $e^L = e^M = x$ and that $\rfs{rk}(L) = \rfs{rk}(M) = \rfs{rk}(K)$.

We show that $L \prec K$.
We already know that $L,K$ are unitary, $L \subseteq K$
and that $\rfs{rk}(L) = \rfs{rk}(K)$.
So it remains to show that $L \not\approx K$.
Suppose by contradiction that $U \in \rfs{Nbr}^K_{\srfs{clp}}(x)$,
$V \in \rfs{Nbr}^L_{\srfs{clp}}(x)$ and $U \cong V$.
Let $\iso{f}{U}{V}$.
Since $x$ is the only accumulation point of $\setm{e_i}{i \in \omega}$,
it follows that $\setm{e_i}{i \in \omega} - U$ is finite.
So there is $i \in \omega - \sigma$ such that $e_i \in U$.
Clearly, $f(e_i) \neq x$, so there is $j \in \sigma$ such that $f(e_i) \in F'_j$.

We consider the sets $S = f\inverse(L'_j \cap V) \cap (L'_i \cap U)$
and $T = (L'_j \cap V) \cap f[L'_i \cap U]$.
Then $ e_i \in S \subseteq L_i$, $T \subseteq L_j$ and $\iso{(f \restriction S)}{S}{T}$.
We check that $S \in \rfs{Clop}(L_i)$ and $T \in \rfs{Clop}(L_j)$.
Since $L'_i,U \in \rfs{Clop}(K)$,
it follows that $L'_i \cap U \in \rfs{Clop}(K)$.
Also, since $L'_j \cap V \in \rfs{Clop}(V)$,
it follows that $f\inverse(L'_j \cap V) \in \rfs{Clop}(U)$.
Hence $f\inverse(L'_j \cap V) \in \rfs{Clop}(K)$.
So $S \in \rfs{Clop}(K)$. Now, $S \subseteq K \cap L_i \subseteq L_i$
and $K \cap L_i \in \rfs{Clop}(L_i)$. Hence $S \in \rfs{Clop}(L_i)$.
A similar calculation shows that $T \in \rfs{Clop}(L_j)$.
So $S \in \rfs{Nbr}^{L_i}_{\srfs{clp}}(e_i)$
and $f[S] \in \rfs{Clop}(L_j)$. That is,
$S \in \rfs{Nbr}^{L_i}_{\srfs{clp}}(e_i)$,
$T \in \rfs{Clop}(L_j)$ and $S \cong T$. But $i \in \omega - \sigma$ and $j \in \sigma$.
So $i \neq j$. These facts contradict (2).
Hence $K \not\approx L$. We have shown that $L \prec K$.

A similar argument shows that $M \prec L$.

We have shown that $M \prec L \prec K$.
So by Theorem~\ref{t3.2}(b), $X$ is not a CO space.
A contradiction. So $X$ is not a CO space.
\smallskip\hfill\proofend

Let $X$ be a space and $x \in X$.
Then $x$ is called a {\it double-limit point} of $X$,
if there are infinite cardinals $\lambda,\mu$ and an embedding
$\fnn{f}{X_{\lambda,\mu}}{X}$
such that $\rfs{cf}(\mu) \geq \aleph_1$
and $f(e^{X_{\lambda,\mu}}) = x$.

We represent $X_{\lambda,\mu}$ as $(\lambda + 1) \cup (\sngltn{0} \times \mu)$.
The subspace $\lambda + 1$ of $X_{\lambda,\mu}$ is denoted by $X^0_{\lambda,\mu}$
and the subspace $\sngltn{\lambda} \cup (\sngltn{0} \times \mu)$
of $X_{\lambda,\mu}$ is denoted by $X^1_{\lambda,\mu}$.

\begin{prop}\label{p7.2}
Let $X$ be a CO compact interval space.

\num{a}
Let $x \in X$ be a double-limit point.
Then there are regular cardinals $\lambda,\mu$
and $U \in \rfs{Nbr}_{\srfs{clp}}^X(x)$ such that $\mu \geq \aleph_1$
and $U \cong X_{\lambda,\mu}$.

\num{b}
The set of double-limit points of $X$ is finite.
\end{prop}

\noindent
{\bf Proof }
(a)
Note that the following facts.
\begin{itemize}
\addtolength{\parskip}{-11pt}
\addtolength{\itemsep}{06pt}
\item[(1)] 
There is a subset $F \subseteq X_{\lambda,\mu}$ such that
$F \cong X_{\srfs{cf}(\lambda),\srfs{cf}(\mu)}$
and
$e^F = e^{X_{\lambda,\mu}}$.
\item[(2)] 
If $U \in \rfs{Nbr}_{\srfs{clp}}(e^{X_{\lambda,\mu}})$,
then $U \cong X_{\lambda,\mu}$.
\vspace{-05.7pt}
\end{itemize}
Let $x$ be a double-limit point of $X$,
and let $f$, $\lambda$ and $\mu$ be as in the
definition of a double-limit point.
By Fact (1), we may assume that $\lambda,\mu$ are regular cardinals.
We may also assume that $\mu \geq \lambda$.
Let $F = \rfs{Rnf}(f)$,
and let $V \in \rfs{Nbr}^X_{\srfs{clp}}(x)$ be a unitary subspace
such that $e^V = x$.
Then either $F \approx V$ or $F \prec^{\srfs{w}} V$.
Suppose by contradiction that $F \prec^{\srfs{w}} V$.
Clearly, $f[X^1_{\lambda,\mu}] \prec F$.
By Corollary~\ref{c5.1}, Theorem~\ref{t3.2} applies to $X$.
So since $f[X^1_{\lambda,\mu}] \prec F \prec^{\srfs{w}} V$,
it follows from \ref{t3.2}(b) that $X$ is not a CO space.
A contradiction, so $V \approx F$.
By Fact (2), there is $W \in \rfs{Nbr}^X_{\srfs{clp}}(x)$
such that $W \cong X_{\lambda,\mu}$.

(b)
It follows from Part (a) that if $X$ is a CO compact interval space
and $x \in X$ is a double-limit point of $X$,
then there is a unique pair
$\pair{\lambda}{\mu} = \pair{\lambda_x}{\mu_x}$ which satisies:
\begin{itemize}
\addtolength{\parskip}{-11pt}
\addtolength{\itemsep}{06pt}
\item[(1)] 
$\mu \geq \aleph_1$ and $\mu \geq \lambda$.
\item[(2)] 
There is an embedding
$\fnn{f}{X_{\lambda,\mu}}{X}$
such that $f(e^{X_{\lambda,\mu}}) = x$.
\vspace{-05.7pt}
\end{itemize}
Also,
$\lambda_x,\mu_x$ are regular cardinals.

Suppose by contradiction that $X$ contains
infinitely many double-limit points.

{\bf Case 1 }
There are $\lambda,\mu$
and an infinite set $A$ of double-limit points of $X$
such that for every $x \in A$,
$\pair{\lambda_x}{\mu_x} = \pair{\lambda}{\mu}$.
Let $L = X_{\lambda,\mu}$ and $M = X^1_{\lambda,\mu}$.
Note that $X^1_{\lambda,\mu} \cong \mu + 1$.
So since $\rfs{rk}(L) = \rfs{rk}(M) = \mu$,
it follows that $M \prec L$.
Then there is a family $\setm{L_i}{i \in \omega}$
such that for every $i \in \omega$,
$L_i \subseteq X$ and $L_i \cong L$,
and for every distinct $i,j \in \omega$, $e^{L_i} \neq e^{L_j}$.
By Theorem~\ref{t3.2}(a), $X$ is not a CO space.
A contradiction.

{\bf Case 2 }
The set $\setm{\pair{\lambda_x}{\mu_x}}
{x \mbox{ is a double-limit point of } X}$ is infinite.
Note that if $\pair{\lambda}{\mu} \neq \pair{\kappa}{\nu}$,
then for every $U \in \rfs{Nbr}^{X_{\lambda,\mu}}_{\srfs{clp}}(e^{X_{\lambda,\mu}})$
and $V \in \rfs{Clop}(X_{\kappa,\nu})$,
$U \not\cong V$.
So $X$ satisfies the conditions of Proposition~\ref{p7.1}.
Hence $X$ is not a CO space. A contradiction.

It follows that the set of double-limit points of $X$ is finite.
\smallskip\hfill\proofend

\begin{prop}\label{p7.3}
Let $X$ be a CO compact interval space.
Then there are no cardinal with uncountable cofinality $\lambda$
and a stationary subset $S \subseteq \lambda$ such that $X_{\lambda,S}$
is embeddable in $X$.
\end{prop}

\noindent
{\bf Proof }
Suppose by contradiction that $X_{\lambda,S}$ is embeddable in $X$.
Without loss of generality $\lambda$ is regular.
Let $T \subseteq S$ be a stationary subset of $\lambda$ such that $S - T$ is stationary
in $\lambda$.
By Proposition~\ref{p5.2}(a),
$\lambda + 1 \prec X_{\lambda,T} \prec X_{\lambda,S}$.
So by Theorem~\ref{t3.2}(b), $X$ is not a CO space.
A contradiction.
\hfill\proofend

\begin{theorem}\label{t7.4}
Let $X$ be a compact scattered interval space.
Suppose that $X$ does not have double-limit points,
and there are no cardinal $\lambda$ with uncountable cofinality
and a stationary subset $S \subseteq \lambda$ such that $X_{\lambda,S}$
is embeddable in $X$.
Then $X$ is homeomorphic to an ordinal space.
\end{theorem}

\noindent
{\bf Proof }
The proof is by induction on the rank of $X$.
If $X$ is a scattered compact interval space with countable rank,
then $X$ is countable. Hence $X$ is homeomorphic
to the interval space of a countable ordinal.
So the claim is true for every space with countable rank.

Suppose that the claim is true for every space with rank $< \alpha$.
Let $X$ be a compact scattered interval space with rank $\alpha$,
and suppose that $X$ does not have double-limit points,
and there are no uncountable cardinal $\lambda$
and a stationary subset $S \subseteq \lambda$ such that $X_{\lambda,S}$
is embeddable in $X$.
We may assume that $X$ is unitary. Let $<$ be a linear ordering of $X$
such that $\tau^X = \tau^{<}$.
We may assume that $e^X \in \rfs{acc}(X^{< e^X})$.

{\bf Case 1 } $\rfs{cf}^-_{\pair{X}{<}}(e^X) = \omega$.
Assume first that $e^X \in \rfs{acc}(X^{> e^X})$.
Then since $e^X$ is not a double-limit point of $X$,
it follows that
$\rfs{cf}^+_{\pair{X}{<}}(e^X) = \omega$.
Let $\setm{x_i}{i \in \omega}$ be a strictly increasing sequence
converging to $e^X$ such that for every $i \in \omega$,
$x_i$ has a successor in $\pair{X}{<}$.
Similarly let $\setm{y_i}{i \in \omega}$ be a strictly deccreasing sequence
converging to $e^X$ such that for every $i \in \omega$,
$y_i$ has a predecessor in $\pair{X}{<}$.
Let $U_0 = X^{\leq x_0}$,
and for every $i > 0$ let $U_i = (x_{i - 1},x_i]$.
Similarly,
let $V_0 = X^{\geq y_0}$,
and for every $i > 0$ let $V_i = [y_i,y_{i - 1})$.
So for every $i \in \omega$, $U_i$ and $V_i$ are clopen subsets of $X$
and $\rfs{rk}(U_i), \rfs{rk}(V_i) < \alpha$.
By the induction hypothesis, for every $i \in \omega$ there are
well orderings $<_{U_i}$ of $U_i$ and
$<_{V_i}$ of $V_i$ which induce the topologies of $U_i$ and of $V_i$.
Define a new linear ordering $<'$ on $X$.
$$
U_0 <' V_0 <' U_1 <' V_1 <' \ldots <' e^X
$$
and for every $x,y \in X$: if for some $i$, $x,y \in U_i$,
then $x <' y$ iff $x <_{U_i} y$,
and if for some $i$, $x,y \in V_i$,
then $x <' y$ iff $x <_{V_i} y$.

It is obvious that $<'$ is a well ordering of $X$ and that
$\tau^{<'} = \tau^X$.

The case that $e^X \not\in \rfs{acc}(X^{> e^X})$
is similar but simpler.

{\bf Case 2 } $\rfs{cf}^-_{\pair{X}{<}}(e^X) > \omega$.
Denote $\rfs{cf}^-_{\pair{X}{<}}(e^X)$ by $\lambda$.
Let $\setm{x_i}{i \in \lambda}$
be a strictly increasing continuous sequence converging to $e^X$.
Let
$$
S = \setm{i \in \lambda}{i \mbox{ is a limit ordinal, and }
x_i \in \rfs{acc}^X(X^{> x_i})}.
$$
For every $i \in S$,
$\rfs{cf}^-_{\pair{X}{<}}(x_i) = \omega$,
for otherwise $x_i$ is a double-limit point.
It follows that $X_{\lambda,S}$ is embeddable in $X$.
So $S$ is not stationary.
Let $C$ be a club of $\lambda$
such that every point of $C$ is a limit point
and such that $C \cap S = \emptyset$.
So for every $i \in C$, $x_i$ has a successor in $X$.
For $i \in C$ let $i \kern2pt+^C\kern2pt 1$ be the successor of $i$ in $C$.
Hence for every $i \in C$,
$X_i \eqdf
(x_i,x_{i \kern2.3pt+^{\kern-1.5pt\tfs{C}}\kern1.1pt 1}]^{\pair{X}{<}}$
is a clopen subset of $X$.
Also, $X_0 \eqdf [\min(X),\min(C)]^{\pair{X}{<}}$ is clopen in $X$.
Clearly, $e^X = \max(X)$, for otherwise,
$e^X$ is a double-limit point in $X$.
Hence $X = (\bigcup_{i \in C \cup \sngltn{0}} X_i) \cup \sngltn{e^X}$.
For every $i \in C \cup \sngltn{0}$,
$\rfs{rk}(X_i) < \rfs{rk}(X)$.
So by the induction hypothesis there is a well-ordering $<_i$
of $X_i$ such that $\tau^{<_i} = \tau^X \restriction X_i$.
Define the linear ordering $<'$ of $X$ such that
$\pair{X}{<'}$ is the lexicographic sum
$\sum_{i \in C \cup {0}} \pair{X_i}{<_i} + \pair{\sngltn{e^X}}{\emptyset}$.
It is easy and left to the reader to check that $\tau^{<'} = \tau^X$.
\smallskip\hfill\proofend

\noindent
{\bf Proof of Theorem~\ref{t6.2}}
Let $X$ be a CO compact interval space.
By Theorem~\ref{t1.2}, $X$ is scattered.

By Proposition~\ref{p7.3}, there is no cardinal $\lambda$
with uncountable cofinality and a stationary set $S \subseteq \lambda$
such that $X_{\lambda,S}$ is embeddable in $X$.

By Proposition~\ref{p7.2}(a) and (b), there are $k \in \omega$,
clopen sets $U_i \subseteq X$, and regular infinite cardinals
$\lambda_i,\mu_i$ \,$i < k$,
such that 
\begin{itemize}
\addtolength{\parskip}{-11pt}
\addtolength{\itemsep}{06pt}
\item[(1)] 
$\mu_i \geq \aleph_1$ and $U_i \cong X_{\lambda_i,\mu_i}$.
\item[(2)] 
$X - \bigcup_{i < k} U_i$ has no double-limit points.
\vspace{-05.7pt}
\end{itemize}
By Proposition~\ref{t7.4}, $X - \bigcup_{i < k} U_i$ is homeomorphic
to an ordinal space $\alpha + 1$.
Let $\mu = \max(\mu_0,\ldots,\mu_{k - 1})$.
Then $\mu + 1 \prec X_{\lambda,\mu}$.
Hence by Proposition~\ref{p3.6},
there is a family $\setm{L_i}{i \in \omega}$ of subsets of $X$
such that for every $i$, $L_i \approx \mu + 1$
and for every $i \neq j$, $e^{L_i} \neq e^{L_j}$.
Clearly, for every $i \in \omega$, $e^{L_i} \not\in \bigcup_{i < k} U_i$.
This implies that $\mu \cdot \omega + 1$ is embeddable in $\alpha + 1$.
Hence $\mu \cdot \omega \leq \alpha$.
\smallskip\hfill\proofend


\newpage
\noindent
{\large\bf Notation index}\vspace{4mm}
\newline
\indexentry{$X_{\lambda,\mu}$. The interval space of $\lambda + 1 + \mu^*$}{1}
\indexentry{$\alpha(X_{\lambda,\mu}) = \max(\lambda,\mu) \cdot \omega$.}{1}
\indexentry{$K_{\srfs{CII}}$. The class of all Hausdorff spaces which are        a continuous\\\indent        image of a compact interval space}{2}
\indexentry{$\rfs{cl}^X(A)$. Closure of $A$ in $X$}{4}
\indexentry{$\rfs{int}^X(A)$. Interior of $A$ in $X$}{4}
\indexentry{$\rfs{acc}^X(A)$. Set of accumulation points of $A$ in $X$}{4}
\indexentry{$\rfs{Nbr}^X(x)$. Set of open neighborhoods of $x$ in $X$}{4}
\indexentry{$\rfs{Nbr}^X_{\srfs{cl}}(x)$.        Set of closed neighborhoods of $x$ in $X$}{4}
\indexentry{$\rfs{Nbr}^X_{\srfs{clp}}(x)$.        Set of clopen neighborhoods of $x$ in $X$}{4}
\indexentry{$\rfs{acc}(\calA)$. The set of accumulation points of a        family of sets $\calA$}{4}
\indexentry{$\calI(N)$. The family of convex components of $N$        in a linear ordering $L$}{5}
\indexentry{$\rfs{Is}(X)$. The set of isolated points of $X$}{8}
\indexentry{$\rfs{D}(X) = X - \rfs{Is}(X)$}{8}
\indexentry{$\rfs{D}_{\alpha}(X)$. The $\alpha$'s derivative of $X$}{8}
\indexentry{$\rfs{rk}(X)$. The rank of $X$}{8}
\indexentry{$\rfs{ker}(X) = \rfs{D}_{\srfs{rk}(X)}(X)$. The maximal perfect        subset of $X$}{8}
\indexentry{$\rfs{Clop}(X)$. The set of clopen subsets of $X$}{8}
\indexentry{$\rfs{Clsd}(X)$. The set of closed subsets of $X$}{8}
\indexentry{$\rfs{Po}(X) = \setm{x \in X}{\mbox{there is }        U \in \rfs{Nbr}(x) \mbox{ such that $\rfs{cl}(U)$ is perfect}}$}{8}
\indexentry{$\calS(X) =        \setm{F \in \rfs{Clsd}(\rfs{ker}(X))}{F \mbox{ is scattered}}$}{8}
\indexentry{$\itOmega(X) = \rfs{sup}(\setm{\rfs{rk}(F)}{F \in \calS(X)})$}{8}
\indexentry{$\rfs{Good}(X)$. The set of good points of $X$}{15}
\indexentry{$\rfs{rk}^{X}(x) =        \rfs{max}(\setm{\alpha}{x \in \rfs{D}_{\alpha}(X)})$}{15}
\indexentry{$\rfs{PD}(X) = X - \rfs{Is}(X) - \rfs{Po}(X)$.        The perfect derivative of $X$}{18}
\indexentry{$\rfs{PD}_{\alpha}(X)$. The $\alpha$'s perfect derivative        of $X$}{18}
\indexentry{$\rfs{prk}(X) =        \rfs{max}(\setm{\alpha}{\rfs{PD}_{\alpha}(X) \neq \emptyset})$.        The perfect rank of $X$}{18}
\indexentry{$\rfs{Pend}(X) = \rfs{PD}_{\srfs{prk}(X)}(X)$.        The perfect end of $X$.}{18}
\indexentry{$K_{\srfs{TH}}$. The class of all compact Hausdorff spaces        that have\\\indent Properties~(TH1)\,-\,(TH3)}{20}
\indexentry{$\whatcalX_t$. The family of pointed spaces associated with $t$        in a P-system $\calP$}{23}
\indexentry{$T_{\calP}$. The index set of a type system $\calP$}{23}
\indexentry{$\calX_t$. The family of spaces associated with $t$        in a P-system $\calP$}{23}
\indexentry{$\whatcalX_{\calP}$. The class of all pointed spaces of        a type system $\calP$}{23}
\indexentry{$\calX_{\calP}$. The class of all spaces of        a type system $\calP$}{23}
\indexentry{$\bfP(A)$. Powerset of $A$}{24}
\indexentry{$\calQ_{\bfGamma}^{\calP}$}{24}
\indexentry{$\whatcalM_{\itGamma}^{\calP} = \setm{\pair{X}{d}}{\pair{X}{d}        \mbox{ is a $\itGamma$-marker}}$}{24}
\indexentry{$\calQ_{\bfGamma}^{\calP} =        \setm{\whatcalM_{\itGamma}^{\calP}}        {\itGamma \in \mathbf \Gamma}$}{24}
\indexentry{$\rfs{acc}_{\mu}(A)$.        The set of $\mu$-accumulation points of $A$}{26}
\indexentry{$\rfs{End}_{\calP}(X)$}{29}
\indexentry{$\rfs{Good}_{\calP}(X)$}{29}
\indexentry{$\whatcalX^{\itOmega}_{\alpha}$.        The class of all pointed spaces which are        $\itOmega$-domnstrative\\\indent        {\thickmuskip=2mu \medmuskip=1mu \thinmuskip=1mu        $(\alpha + 1)$-codes        with a member of $\rfs{Pend}(X)$ as their distiguished point}}{30}
\indexentry{$\calP^{\itOmega}$. The P-system of        $\itOmega$-demonstrative codes}{30}
\indexentry{$\calB_{\lambda,\vecmu}$. A base of $X_{\lambda,\mu}$}{34}
\indexentry{$P^{< x} = \setm{y \in P}{y < x}$}{36}
\indexentry{$\rfs{cf}^-_{\pair{L}{<}}(a)$ and $\rfs{cf}^+_{\pair{L}{<}}(a)$. The cofinality of $a$ from the left and the\\ \indent cofinality of $a$ from the right}{36}
\indexentry{$X_{\lambda,S}$}{58}
\indexentry{$V_a$. For a member $a$ of a BA $B$,        $V_a = \setm{x \in \rfs{Ult}(B)}{a \in x}$}{68}

\newpage
\noindent
{\large\bf Definition index}\vspace{4mm}
\newline
\indexentry{accumulation point of $\calA$}{4}
\indexentry{accumulation point: $\lambda$-accumulation point.        A point $x$ is a\\\indent        $\lambda$-accumulation point of $A$ if $\abs{U \cap A} = \lambda$        for every $U \in \rfs{Nbr}(x)$}{13}
\indexentry{attained: $\itOmega(X)$ is not attained in $X$}{12}
\indexentry{code: $\alpha$-code}{18}
\indexentry{code}{18}
\indexentry{collectionwise Hausdorff space.        $X$ is collectionwise Hausdorff if every\\\indent        relatively discrete subset of $X$ has a Hausdorff system}{5}
\indexentry{demonstrative set: $\itOmega$-demonstrative set}{19}
\indexentry{dense: $\lambda$-dense linear ordering}{13}
\indexentry{derivative.}{8}
\indexentry{filler}{24}
\indexentry{good point}{12}
\indexentry{Hausdorff system}{4}
\indexentry{marker: $\itGamma$-marker}{24}
\indexentry{marker: $\mu$-special $\sngltn{t}$-marker}{26}
\indexentry{occurs: $\calX$ occurs in $Y$}{23}
\indexentry{occurs: A P-system occurs in $Y$}{23}
\indexentry{P-system. Abbreviation of a proliferation system}{23}
\indexentry{pairwise disjoint family of subsets of $X$}{4}
\indexentry{pairwise disjoint set of elements        of a Boolean algebra}{68}
\indexentry{perfect derivative}{18}
\indexentry{perfect end}{18}
\indexentry{perfect kernel}{8}
\indexentry{perfect rank}{18}
\indexentry{perfect set. A set which does not have isolated points in its        \\\indent relative topology}{8}
\indexentry{pointed space. A pair $\pair{X}{x}$,        where $X$ is a topological space and $x \in X$}{23}
\indexentry{principal}{38}
\indexentry{proliferation system}{23}
\indexentry{rank of $x$ in $X$}{15}
\indexentry{rank. The rank of $X$, the first ordinal $\alpha$ such that        $\rfs{D}_{\alpha}(X)$ is finite\\\indent        or perfect}{8}
\indexentry{relatively discrete.        $A$ is relatively discrete if $A \cap \rfs{acc}(A) = \emptyset$}{4}
\indexentry{scattered space}{2}
\indexentry{strong Hausdorff system}{4}
\indexentry{strongly Hausdorff for convergent sequences}{11}
\indexentry{strongly Hausdorff space.        $X$ is strongly Hausdorff if every relatively\\\indent        discrete subset of $X$ has a strong Hausdorff system}{5}
\indexentry{tight family of subsets of $X$}{4}
\indexentry{tightly Hausdorff space.        $X$ is tightly Hausdorff if every relatively\\\indent        discrete subset of $X$ has a tight Hausdorff system}{5}
\indexentry{tree-like clopen system}{35}
\indexentry{type system}{23}

\newpage
\begin{thebibliography}{DGMM}
\medskip
\smallskip

\bibitem[B]{B} Bonnet R.:
On superatomic Boolean algebras,
{\it Finite and infinite combinatorics in sets and logic},
(Banff, AB, 1991), NATO Adv. Sci. Inst. Ser. C Math. Phys. Sci., {\bf 411},
Kluwer Acad. Publ., Dordrecht, 1993, pp.~31--62.

\bibitem[BS]{BS} Bonnet R. and Shelah S.\,:
On HCO spaces. An uncountable compact $T\sb 2$ space, different from $\aleph\sb 1+1$, which is homeomorphic to each of its uncountable closed subspaces,
{\it Israel J. Math.} {\bf 84} (1993), no.3, pp.~289--332.


\bibitem[BBR]{BBR} Bekkal M., Bonnet R. and Rubin M.\,
Compact interval spaces in which all closed subsets are homeomorphic to clopen ones,
{\it Order} {\bf 9} (1992), no.1, pp.~69--95 and {\it Order} {\bf 9} (1992), no. 2, pp.~177--200.

\bibitem[H]{H} Heindorff L.\,:
On subalgebras of Boolean interval algebras.
{\it Proc. Amer. Math. Soc.} {\bf 125} (1997), no. 8, pp.~2265--2274

\bibitem[S]{S} Shelah S\,:
Factor = quotient, uncountable Boolean algebras,
number of endomorphism and width, {\it Math. Japonica} {\bf 37} (1992),
pp.~385-400.

\end{thebibliography}
\end{document}